\documentclass[10pt]{amsart}

\usepackage{amsfonts,amssymb,amsmath,mathrsfs,graphicx}
\usepackage{float}
\usepackage{epsfig}
\usepackage{CJK,color}

\definecolor{red}{rgb}{1.00,0.00,0.00}
\definecolor{blue}{rgb}{0.00,0.00,0.63}
\definecolor{black}{rgb}{0.00,0.00,0.00}

\newtheorem{theorem}{Theorem}[section]
\newtheorem{lemma}{Lemma}[section]

\newtheorem{proposition}{Proposition}[section]
\newtheorem{remark}{Remark}[section]

\def\charf {\mbox{{\text 1}\kern-.30em {\text l}}}

\def\mb{\mathbf}

\def\f{\frac}
\def\t{\theta}
\def\d{\delta}
\def\a{\alpha}
\def\b{\beta}
\def\i{\infty}
\def\z{\zeta}

\def\pa{\partial}

\def\di{\displaystyle}
\def\r{\rho}
\def\wt{\widetilde}
\def\k{\kappa}

\begin{document}
%

\title[Stability of wave patterns to bipolar VPB]
{Stability of the Superposition of a Viscous Contact Wave with two Rarefaction Waves to the bipolar Vlasov-Poisson-Boltzmann System}

\author[Li]{Hailiang Li}
\address[Hailiang Li]{\newline School of Mathematical Sciences, Capital Normal University, Beijing 100048, P. R. China
\newline and Beijing Center of Mathematics
and Information Sciences, Beijing 100048, P. R. China}
\email{Hailiang-Li@mail.cnu.edu.cn}

\author[Wang]{Teng Wang}
\address[Teng Wang]{\newline Department of Mathematics, School of Science, Beijing Forestry University, Beijing 100083, P. R. China}
\email{tengwang@amss.ac.cn}

\author[Wang]{Yi Wang}
\address[Yi Wang]{\newline Institute of Applied Mathematics, AMSS, CAS, Beijing 100190, P. R. China,
\newline and School of Mathematical Sciences, University of Chinese Academy of Sciences, Beijing 100049, P. R. China,
\newline and Beijing Center of Mathematics
and Information Sciences, Beijing 100048, P. R. China}
\email{wangyi@amss.ac.cn}

\date{\today}



\thanks{\textbf{Acknowledgment.} The work of H. L. Li is partially supported by  the NNSFC grants No. 11671384, 11231006,
and 11461161007, and by the Beijing New Century Baiqianwan Talent Project.
The work of T. Wang is partially supported by the Fundamental Research Funds for the Central Universities (No. 2015ZCQ-LY-01 and No. BLX2015-27),
and the NNSFC  grant No. 11601031.
The work of Y. Wang is partially supported by NNSFC grant No. 11671385 and Youth Innovation Promotion Association of CAS}

\maketitle

\begin{abstract}
We investigate the nonlinear stability of the superposition of a viscous contact wave and two rarefaction waves for one-dimensional bipolar Vlasov-Poisson-Boltzmann (VPB) system, which can be used to describe the transportation of charged particles under the additional electrostatic potential force. Based on a new micro-macro type decomposition around the local Maxwellian related to the bipolar VPB system in our previous work \cite{LWYZ}, we prove that the superposition  of  a viscous contact wave and two rarefaction waves is time-asymptotically stable to 1D bipolar VPB system under some smallness conditions on the initial perturbations and wave strength, which implies that this typical composite wave pattern is nonlinearly stable under the combined effects of the binary collisions, the electrostatic potential force, and the mutual interactions of different charged particles.  Note that this is the first result about the nonlinear stability of the combination of two different wave patterns for the Vlasov-Poisson-Boltzmann system.
\end{abstract}

\maketitle \centerline{\date}
%
%

\section{Introduction}
\setcounter{equation}{0}
In this paper, we investigate the time-asymptotic stability of the superposition of  a viscous contact wave and two rarefaction waves to the 1D bipolar
Vlasov-Poisson-Boltzmann system, which reads
\begin{equation}\label{VPB0}
\left\{
\begin{array}{ll}
\di
F_{At} +\xi_1\partial_x F_A+\partial_x\Phi\partial_{\xi_1} F_A= Q(F_A,F_A+F_B), \\
\di
F_{Bt} +\xi_1\partial_x F_B-\partial_x\Phi\partial_{\xi_1} F_B= Q(F_B,F_A+F_B), \\
\di \partial_{xx} \Phi=\int (F_A-F_B)d\xi,
\end{array}
\right.
\end{equation}
where $x\in \mb{R}^1, t\in \mb{R}^+, \xi=(\xi_1,\xi_2,\xi_3)^t\in \mb{R}^3$, and $F_A(t,x,\xi) $ and $ F_B(t,x,\xi)$ denote the
density distribution function of two-species particles (e.g. ions and electrons) at time-space point $(t,x)$ with the velocity $\xi$ respectively. For the hard sphere model, the collision operator
$Q(f,g)$ takes the following  bilinear form
$$
 Q(f,g)(\xi) \equiv \f{1}{2}
\int_{{\mb R}^3}\!\!\int_{{\mb S}_+^2} \Big(f(\xi')g(\xi_*')-f(\xi) g(\xi_*) \Big)
|(\xi-\xi_*)\cdot\Omega)|
 \; d \xi_* d\Omega,
$$
where $\xi',\xi_*'$ are the velocities after an elastic collision of
two particles with the pre-collision velocities  $\xi,\xi_*$, respectively,
and the unit vector $\Omega \in {\mb S}^2_+=\{\Omega\in
{\mb S}^2:\ (\xi-\xi_*)\cdot \Omega\geq 0\}$. The conservation of
momentum and energy of the collided particles yields the following relations between the
velocities before and after the collision:
$$
 \xi'= \xi -[(\xi-\xi_*)\cdot \Omega] \; \Omega, \qquad
 \xi_*'= \xi_* + [(\xi-\xi_*)\cdot \Omega] \; \Omega.
$$
The initial values and the far-field states to the system \eqref{VPB0} are given by
\begin{equation}\label{VPB0-i}
\left\{\begin{aligned}
&\di F_A(t=0,x,\xi)=F_{A0}(x,\xi)\rightarrow \mathbf{M}_{[\rho_\pm,u_\pm,\t_\pm]}(\xi),\quad {\rm as}~~x\rightarrow\pm\infty,\\
&\di F_B(t=0,x,\xi)=F_{B0}(x,\xi)\rightarrow \mathbf{M}_{[\rho_\pm,u_\pm,\t_\pm]}(\xi),\quad {\rm as}~~x\rightarrow\pm\infty,\\
&\di \Phi_x\rightarrow 0,\quad {\rm as}~~x\rightarrow\pm\infty,
\end{aligned}
\right.
\end{equation}
where $u_\pm=(u_{1\pm},0,0)^t$ and $\rho_\pm>0, u_{1\pm}, \t_\pm>0$ are the prescribed constant states such that the two states
$(\rho_\pm,u_\pm,\t_\pm)$ are connected by the Riemann solution with the combination of  1-rarefaction wave, 2-contact discontinuity and 3-rarefaction wave to the corresponding 1D Euler system
\begin{equation}\label{euler'}
\left\{
\begin{array}{l}
\di\rho_t+(\rho u_1)_x=0,\\
\di(\rho u_1)_t+(\rho u_1^2+p)_x=0,\\
\di(\rho u_i)_t+(\rho u_1u_i)_x=0,~i=2,3,\\
\di[\rho(e+\f{|u|^2}{2})]_t+[\rho u_1(e+\f{|u|^2}{2})+pu_1]_x=0.
\end{array}
\right.
\end{equation}
Note that the above fluid variables $(\rho, u,\theta)$ and the equilibrium Maxwellian $\mathbf{M}_{[\rho,u,\theta]}$ related to the density distribution functions $F_A$ and $F_B$ are defined in \eqref{macro} and \eqref{maxwellian}, respectively.
As in \cite{LWYZ}, set the total mass by $F_1$ and the neutrality by $F_2$ of the distribution functions as follows
$$
F_1=\f{F_A+F_B}{2},\qquad F_2=\f{F_A-F_B}{2}.
$$
Then one has the following reformulated system
\begin{equation}\label{VPB}
\left\{
\begin{array}{ll}
\di
F_{1t} +\xi_1\partial_x F_1+\partial_x\Phi\partial_{\xi_1} F_2= 2Q(F_1,F_1), \\
\di
F_{2t} +\xi_1\partial_x F_2+\partial_x\Phi\partial_{\xi_1} F_1= 2Q(F_2,F_1), \\
\di \partial_{xx} \Phi=2\int F_2d\xi:=2n_2,
\end{array}
\right.
\end{equation}
with the initial values $(F_{10}, F_{20})$ and the far-fields states satisfying
\begin{equation}
\begin{array}{ll}
\di F_1(t=0,x,\xi)=F_{10}(x,\xi)\rightarrow \mathbf{M}_{[\rho_\pm,u_\pm,\t_\pm]}(\xi),\quad {\rm as}~~x\rightarrow\pm\infty,\\
\di F_2(t=0,x,\xi)=F_{20}(x,\xi)\rightarrow 0,\quad {\rm as}~~x\rightarrow\pm\infty,\\
\di \Phi_x\rightarrow 0,\quad {\rm as}~~x\rightarrow\pm\infty.
\end{array}
\end{equation}

The bipolar Vlasov-Poisson-Boltzmann (abbreviated as VPB for simplicity below) system \eqref{VPB0} can be viewed as two Boltzmann equations coupled with each other through the self-consistent electrostatic potential determined by the Poisson equation. Moreover, if the effect of the electrons is neglected, the system~\eqref{VPB0} can reduce to the unipolar VPB equations to simulate the motion of one species. However, these bipolar and unipolar VPB systems are by no means the simple extension of the Boltzmann equation in the sense that the global solutions to these VPB system may  possibly admit, due to the influence of electric field, some different dynamical behaviors from the Boltzmann equation. Indeed, as analyzed in \cite{DS,LWYZ,LYZ1,LYZ}, the long time asymptotic behaviors of global solutions to VPB system are rather complicated under the combined effects of the binary collision among the particles of same charge, the electrostatic potential force, and/or the mutual interactions between different charged particles. In particular, it was proved for the unipolar VPB system in \cite{LYZ1} that the appearance of electric field affects the spectrum structure of the linearized unipolar VPB system and causes the slower but optimal (compared with the Boltzmann equation) time convergence rates of global solution to the equilibrium state, and there does not exist wave propagation any more, the readers can refer to \cite{DS, LYZ, Zhong} and references therein for details.

Yet, a completely different dynamical phenomena and long time behaviors of global solution are observed  in \cite{LWYZ,LYZ} for the bipolar VPB system~\eqref{VPB0}, namely, the global solution to the bipolar VPB system~\eqref{VPB0} may have  similar time asymptotical behaviors to that of the Boltzmann equation. Indeed, it has been shown  in \cite{LYZ}  that the linearized bipolar VPB system around the global Maxwellian consists of a decoupled system: one is the linearized Boltzmann equation for the distribution function $F_1$ which admits the wave modes at lower frequency and the algebraic time decay rates of the function $F_1$, the other is the equation of unipolar VPB type for the neutral function $F_2$ of the particles with different charge which admits the spectral gap at lower frequency and causes the strong neutrality in the sense that the neutral function $F_2$  and the electric field $\partial_x\Phi$ decay exponentially in time.
On the other hand, it is well-known that the Boltzmann equation is asymptotically equivalent to the compressible Euler equations as illustrated by the Hilbert's expansions~\cite{CC}, and there are three basic waves to the compressible Euler equations: two nonlinear waves (shock wave and rarefaction wave), and one linearly degenerate wave (contact discontinuity). As shown in~\cite{Huang-Xin-Yang,Huang-Yang, Liu-Yu-1, Wang-Wang, Yu}, the global solution to Boltzmann equation can have rich wave phenomena and the time-asymptotic stability of above basic wave patterns is proved for one-dimensional Boltzmann equation.
Therefore, due to the appearance of wave modes and the spectral gap of the linearized bipolar VPB system~\eqref{VPB0} around the global Maxellian as made in \cite{LYZ}, it is natural and interesting to investigate the nonlinear wave phenomena of the bipolar VPB system~\eqref{VPB0}  and understand the combined influence of electric field and mutual interactions between different charged particles.  Moreover, we have proved in \cite{LWYZ} that the global solutions the bipolar VPB system shall tend time-asymptotically to either viscous shock profile or rarefaction wave-fan if the initial data is near the corresponding local Maxwellian related to the shock or rarefaction wave profiles.

The main purpose of the present paper is to prove the nonlinear stability of the superposition of a viscous contact wave and two rarefaction waves for the bipolar VPB system~\eqref{VPB0}, as a continuation of the previous work \cite{LWYZ}. In \cite{LWYZ}  we first set up a new micro-macro type decomposition around the local Maxwellian related to the bipolar VPB system \eqref{VPB0} (e.g. \eqref{VPB}) and established a generic framework to prove the nonlinear stability of shock profile and rarefaction wave respectively to the 1D bipolar VPB system \eqref{VPB0} (or \eqref{VPB}). Our main result in the present paper can be stated roughly as follows: it can be proved that the composite waves with one viscous contact wave and two rarefaction waves is nonlinearly stable for the bipolar VPB system as the time $t$ tends to infinity, if the initial perturbations and total wave strength is suitably small (refer to Theorem~\ref{thm} for details).

\bigskip

There have been important progress on the existence and asymptotical behaviors of solutions to the VPB system. For instance,  Mischler first proved the global existence of renormalized (weak) solution for large initial data~\cite{Mi}, and Guo obtained the first global existence result on strong solution in torus when the initial data is near a global Maxwellian~\cite{Guo1}. And the global existence of classical solution in $\mathbf{R}^3$ was given in \cite{YYZ, YZ} . The case
with general stationary background density function was studied in \cite{DY}, and the perturbation of vacuum was investigated in \cite{DYZ, DS}. For the case case of soft potential, one can refer to\cite{DYZh, WXXZ, XXZ} and references therein.  Recently, Li-Yang-Zhong \cite{LYZ}  analyze the
spectrum of the linearized VPB system (unipolar and bipolar) and
obtain the optimal decay rate of solutions to the nonlinear system near global Maxwellian. See also the works on the stability of global Maxwellian and the optimal time decay rate in \cite{YY, W}. For the stability of single basic wave patterns, Duan-Liu proved the stability of rarefaction wave to a unipolar VPB system~\cite{Duan-Liu}, which can be viewed as an approximation of bipolar VPB system  \eqref{VPB0} when the electron density is very rarefied and reaches a local equilibrium state with small electron mass compared with the ion. Recently, Li-Wang-Yang-Zhong \cite{LWYZ} gave a unified micro-macro decomposition to the bipolar VPB system and proved the stability of rarefaction wave and viscous shock profile to the bipolar VPB system by this new decomposition. Then, Duan-Liu \cite{Duan-liu2} generalized the result in \cite{Duan-Liu} to the bipolar VPB system with the disparate mass and Huang-Liu \cite {Huang-liu} showed the stability of a single viscous contact wave case.

\bigskip

We should mention that the stability of elementary wave patterns of the Boltzmann equation are well-analyzed. For instance, the pioneering result on the stability of viscous shock wave was first proved by Liu-Yu \cite{Liu-Yu} with the zero total macroscopic mass condition by the energy method based on the micro-macro decomposition while the existence of viscous shock profile to the Boltzmann equations is given by Caflish-Nicolaenko \cite{Caflish-Nicolaenko} and Liu-Yu \cite{Liu-Yu-1}. Then stability of rarefaction wave is proved by  Liu-Yang-Yu-Zhao \cite{Liu-Yang-Yu-Zhao} and the stability of viscous contact wave, which is the viscous version of contact discontinuity, by Huang-Yang \cite{Huang-Yang} with the zero mass condition and Huang-Xin-Yang \cite{Huang-Xin-Yang} without the zero mass condition.  Furthermore, Yu \cite{Yu} proved the stability of single viscous shock profiles without zero mass condition by the point-wise method based on the Green function around the viscous shock profile. Recently, Wang-Wang \cite{Wang-Wang} proved the stability of superposition of two viscous shock profiles to the Boltzmann equation without the zero mass condition by the weighted characteristic energy method. Moreover, the hydrodynamic limit of Boltzmann equation to the Riemann solution to Euler system is investigated extensively in \cite{Huang-Wang-Yang, Huang-Wang-Yang-2, Huang-Wang-Yang-3, Xin-Zeng, Yu1} and finally justified in \cite{hwwy} to the generic Riemann solutions.

\bigskip

The rest part of the paper is arranged as follows. In section~\ref{decompisition} we present the new micro-macro decomposition around the local Maxellian for the bipolar VPB system \eqref{VPB0}. Then, the main result  Theorem \ref{thm} on the stability of the superposition
of a viscous contact wave and two rarefaction waves are stated in section~\ref{result}.
The lower and higher order a-priori estimates will be established in section 4 and 5, respectively. Finally, section 6 is devoted to the proof of the main result.

%
%
%
%

\section{New Micro-Macro Decompositions}
\label{decompisition}
\setcounter{equation}{0}
First we recall the new micro-macro decompositions around the local Maxwellian to the bipolar VPB system \eqref{VPB} in \cite{LWYZ}.
The equation $\eqref{VPB}_1$ can be viewed as the Boltzmann equation with additional force, we make use of the micro-macro decomposition as introduced by Liu-Yu \cite{Liu-Yu} and Liu-Yang-Yu \cite{Liu-Yang-Yu}. For any solution $F_1(t,x,\xi)$ to Eq.~$\eqref{VPB}_1$, there are five macroscopic (fluid) quantities: the mass density $\r(t,x)$, the momentum $m(t,x)=\r u(t,x)$,
and the total energy $E(t,x)=\r\big(e+\frac12|u|^2\big)(t,x)$ defined by
\begin{equation}
\left\{
\begin{array}{l}
\di\rho(t,x)=\int_{\mb{R}^3}\varphi_0(\xi)F_1(t,x,\xi)d\xi,\\
\di\rho u_i(t,x)=\int_{\mb{R}^3}\varphi_i(\xi)F_1(t,x,\xi)d\xi,~i=1,2,3,\\
\di\rho(e+\f{|u|^2}{2})(t,x)=\int_{\mb{R}^3}\varphi_4(\xi)F_1(t,x,\xi)d\xi,
\end{array}
\right. \label{macro}
\end{equation}
where  $\varphi_i(\xi)$ $(i=0,1,2,3,4)$ are the collision
invariants given by
\begin{equation}
 \varphi_0(\xi) = 1,~~~
 \varphi_i(\xi) = \xi_i~  (i=1,2,3),~~~
 \varphi_4(\xi) = \f{1}{2} |\xi|^2,
\label{collision-invar}
\end{equation}
and satisfy
$$
\int_{{\mb R}^3} \varphi_i(\xi) Q(g_1,g_2) d \xi =0,\quad {\textrm
{for} } \ \  i=0,1,2,3,4.
$$
Define the local Maxwellian $\mb{M}$ associated to the solution $F_1(t,x,\xi)$  to Eq.~$\eqref{VPB}_1$  in terms of the fluid quantities by
\begin{equation}
\mb{M}:=\mb{M}_{[\rho,u,\t]} (t,x,\xi) = \f{\rho(t,x)}{\sqrt{ (2 \pi
R \t(t,x))^3}} e^{-\f{|\xi-u(t,x)|^2}{2R\t(t,x)}},  \label{maxwellian}
\end{equation}
where $\t(t,x)$ is the temperature which is related to the internal
energy $e(t,x)$ by $e=\frac{3}{2}R\t$ with $R>0$ the gas constant,  and $u(t,x)=\big(u_1(t,x),u_2(t,x),u_3(t,x)\big)^t$
is the fluid velocity. Then, the collision operator of $Q(f,f)$ can be linearized to be $\mb{L}_\mb{M}$ with respect to the local Maxwellian $\mb{M}$ by
\begin{equation}
\mb{L}_\mb{M} g=2Q(\mb{M}, g)+ 2Q(g,\mb{M}).  \label{L_M}
\end{equation}
The null space $\mathfrak{N}_1$ of $\mb{L}_\mb{M}$ is spanned by the collision invariants $\varphi_i(\xi)~(i=0,1,2,3,4)$  in \eqref{collision-invar}.

Define an inner product $ \langle g_1,g_2\rangle_{\widetilde{\mb{M}}}$ for $g_i\in L^2(\mb{R}^3_\xi)$ with respect to the given local Maxwellian $\wt{\mb{M}}$ as:
\begin{equation}
 \langle g_1,g_2\rangle_{\widetilde{\mb{M}}}\equiv \int_{{\mb R}^3}
 \f{1}{\widetilde{\mb{M}}}g_1(\xi)g_2(\xi)d \xi.\label{product}
\end{equation}
For simplicity, if $\widetilde{\mb{M}}$ is the local Maxwellian $\mb{M}$ in \eqref{maxwellian}, we shall use the notation $\langle\cdot,\cdot\rangle$ instead of $\langle\cdot,\cdot\rangle_{\mb{M}}$.
Furthermore, there exists a positive constant $\wt\sigma_1>0$ such
that it holds for any function $g(\xi)\in \mathfrak{N}_1^\bot$ (cf. \cite{CC,Grad}) that
\begin{equation}
\langle g,\mb{L}_\mb{M}g\rangle
\le
 -\wt\sigma_1\langle \nu(|\xi|)g,g\rangle ,  \label{H-thm}
\end{equation}
where $\nu(|\xi|)\sim (1+|\xi|)$ is the collision frequency for the
hard sphere collision. The dissipative property of the linearized operator $\mb{L}_\mb{M}$ in \eqref{H-thm} comes from the celebrated H-theorem for the Boltzmann equation, which is crucial to estimates the microscopic parts for the Boltzmann equation.
With respect to the inner product $\langle\cdot,\cdot\rangle$,  the following pairwise orthogonal bases span the macroscopic space  $\mathfrak{N}_1$
\begin{equation}
\left\{
\begin{array}{l}
 \chi_0(\xi) \equiv {\di\f1{\sqrt{\rho}}\mb{M}}, \quad
 \chi_i(\xi) \equiv {\di\f{\xi_i-u_i}{\sqrt{R\t\rho}}\mb{M}} \ \ {\textrm {for} }\ \  i=1,2,3, \\[2mm]
 \chi_4(\xi) \equiv
 {\di\f{1}{\sqrt{6\rho}}(\f{|\xi-u|^2}{R\t}-3)\mb{M}},\quad \langle\chi_i,\chi_j\rangle=\delta_{ij}, ~i,j=0,1,2,3,4.
 \end{array}
\right.\label{orthogonal-base}
\end{equation}
In terms of the above five orthogonal  bases, the macroscopic projection $\mb{P}_0$ from $L^2(\mb{R}^3_\xi)$ to  $\mathfrak{N}_1$ and the microscopic projection $\mb{P}_1$  from $L^2({\mb R}^3_\xi)$ to $\mathfrak{N}_1^\bot$ can be defined as
\begin{equation*}
 \mb{P}_0g = {\di\sum_{j=0}^4\langle g,\chi_j\rangle\chi_j},\qquad
 \mb{P}_1g= g-\mb{P}_0g.
\end{equation*}
A function $g(\xi)$ is said to be microscopic or non-fluid, if it holds
$$
\int g(\xi)\varphi_i(\xi)d\xi=0,~i=0,1,2,3,4,
$$
where $\varphi_i(\xi)$ are the collision invariants defined in \eqref{collision-invar}.

Based on the above preparation, the solution $F_1(t,x,\xi)$  to Eq.~$\eqref{VPB}_1$ can be decomposed into
the macroscopic (fluid) part, i.e., the local Maxwellian $\mb{M}=\mb{M}(t,x,\xi)$
defined in \eqref{maxwellian}, and the microscopic (non-fluid) part, i.e.
$\mb{G}=\mb{G}(t,x,\xi)$ (cf. \cite{Liu-Yu, Liu-Yang-Yu}):
$$
F_1(t,x,\xi)=\mb{M}(t,x,\xi)+\mb{G}(t,x,\xi),\quad
\mb{P}_0F_1=\mb{M},~~~\mb{P}_1F_1=\mb{G},
$$
By the decomposition for $F_1=\mb{M}+\mb{G}$, one can derive from \eqref{VPB} and \eqref{macro}  the following
fluid-part system
\begin{equation}\label{F1-f-euler}
\left\{
\begin{array}{l}
\di \rho_{t}+(\rho u_1)_x=0,\\
\di (\rho u_1)_t+(\rho u_1^2 +p)_x-(\f{\Phi_x ^2}{4})_x=\f{4}{3}(\mu(\t)
u_{1x})_x-\int \xi_1^2\Pi_xd \xi,  \\
\di (\rho u_i)_t+(\rho u_1u_i)_x=(\mu(\t)
u_{ix})_x-\int \xi_1 \xi_i\Pi_xd \xi,~ i=2,3,\\
\di [\rho(\t+\f{|u|^2}{2})+\f{\Phi_x ^2}{4}]_t+[\rho
u_1(\t+\f{|u|^2}{2})+pu_1]_x=(\k(\t)\t_x)_x+\f{4}{3}(\mu(\t)u_1u_{1x})_x\\
\di\qquad\qquad +\sum_{i=2}^3(\mu(\t)u_iu_{ix})_x
-\int\f12 \xi_1| \xi|^2\Pi_xd \xi,
\end{array}
\right.
\end{equation}
where the viscosity coefficient $\mu(\t)>0$ and the heat
conductivity coefficient $\k(\t)>0$ are smooth functions of the
temperature $\t$. Here,  we renormalize the gas constant $R$ to be
$\f{2}{3}$ so that $e=\t$ and $p=\f23\rho\t$. The non-fluid equation for $F_1$ can be written as:
\begin{equation}
\mb{G}_t+\mb{P}_1( \xi_1\mb{M}_x)+\mb{P}_1(
\xi_1\mb{G}_x)+\mb{P}_1(\Phi_x \partial_{\xi_1}F_2)
=\mb{L}_\mb{M}\mb{G}+2Q(\mb{G}, \mb{G}). \label{F1-G-euler}
\end{equation}
Therefore,  $\mb{G}$ can be expressed explicitly by
\begin{equation}
\mb{G}= \mb{L}_\mb{M}^{-1}[\mb{P}_1(\xi_1\mb{M}_x)] +\Pi
\label{G-euler}
\end{equation}
with
\begin{equation}\label{Pi}
\Pi=\mb{L}_\mb{M}^{-1}[\mb{G}_t+\mb {P}_1(\xi_1\mb{G}_x)+\mb
{P}_1(\Phi_x\partial_{\xi_1} F_2)-2Q(\mb{G}, \mb{G})].
\end{equation}
Remark that the fluid-system \eqref{F1-f-euler} is the compressible Navier-Stokes type system strongly coupled with microscopic terms determined by \eqref{F1-G-euler} and electric field terms in $\eqref{VPB}_3$, and the system \eqref{F1-f-euler}, \eqref{F1-G-euler},  $\eqref{VPB}_3$ is not self-closed and the new decomposition to $F_2$-equation $\eqref{VPB}_2$ around the local Maxellian is needed (cf. \cite{LWYZ}):
$$
F_2=\frac{n_2}{\rho}\mb{M}+\mb{P}_c F_2.
$$
with the linearized collision operator $\mb{N}_\mb{M}$ (around the local Maxellian $\mb{M}$ to $F_1$ in \eqref{maxwellian}) defined by
$$
\mb{N}_{\mb{M}} h=2Q(h,\mb{M}).
$$
The null space $\mathfrak{N}_2$ of the operator $\mb{N}_\mb{M}$ is
spanned by the single macroscopic variable:
$$
\chi_0(\xi)=\f{\mb{M}}{\sqrt\rho},
$$
due to the quite different collision structure from the linearized operator $\mb{L}_{\mb{M}}$, whose null space $\mathfrak{N}_1$ is spanned by five macroscopic variables $\chi_j(\xi)~(j=0,1,2,3,4).$
Furthermore, there exists a positive constant $\wt\sigma_2>0$ such
that it holds for any function $g(\xi)\in \mathfrak{N}_2^\bot$ (cf. \cite{AY, LWYZ}) that
$$
\langle g,\mb{N}_\mb{M}g\rangle \le -\wt\sigma_2\langle
\nu(|\xi|)g,g\rangle ,
$$
where $\nu(|\xi|)\sim(1+|\xi|)$ is the collision frequency for the hard sphere collision.
Consequently, the linearized collision operator $\mb{N}_\mb{M}$ is dissipative on $\mathfrak{N}_2^\bot$, and its inverse
$\mb{N}_\mb{M}^{-1}$ exists and is a bounded operator on $\mathfrak{N}_2^\bot$.
Taking the inner product of the equation \eqref{VPB}$_2$ and the collision invariants $\varphi_0(\xi)=1$ with respect to $\xi$ over $\mb{R}^3$,
the macroscopic part $n_2$ satisfies the following equation
\begin{equation}\label{n21}
n_{2t}+(\int \xi_1F_2d\xi)_x=0,
\end{equation}
or equivalently,
\begin{equation}\label{n2-euler}
\begin{array}{ll}
\di n_{2t}+(u_1n_2)_x+ \Big(\f{\k_1(\t)}{R\t}\Phi_x \Big)_x- \Big(\k_1(\t)(\f{n_2}{\r})_x\Big)_x\\
\di
 =-\Big(\f{n_2}{\r}\int \xi_1 \mb{N}_{\mb{M}}^{-1}\big[\mb{P}_c(\xi _1\mb{M}_x) )\big]d\xi \Big)_x-\Big(\int \xi_1 \mb{N}_{\mb{M}}^{-1}\big[\mb{P}_c(\xi_1 (\mb{P}_cF_2)_x)\big]d\xi\Big)_x\\
\di
 - \Big(\int \xi_1 \mb{N}_{\mb{M}}^{-1}\big[\Phi_x \mb{G}_{\xi_1}\big]d\xi \Big)_x -\Big(\int \xi_1 \mb{N}_{\mb{M}}^{-1}\Big[\partial_t(\mb{P}_c F_2)+(\f{\mb{M}}{\rho})_t~n_2
 -2Q(F_2,\mb{G})\Big]d\xi\Big)_x,
\end{array}
\end{equation}
and the microscopic part $\mb{P}_c F_2$ satisfies the equation
\begin{equation}\label{F2-pc-1}
\partial_t(\mb{P}_c F_2)-\mb{N}_{\mb{M}}(\mb{P}_c F_2)+\mb{P}_c(\xi_1 F_{2x})+\mb{P}_c (\Phi_x \partial_{\xi_1}F_1)+(\f{\mb{M}}{\rho})_t~n_2=2Q(F_2,\mb{G}).
\end{equation}
The detailed derivation of the equations \eqref{n2-euler} and \eqref{F2-pc-1} can be found in \cite{LWYZ}. Note that the equation \eqref{n2-euler} is a diffusive equation with the damping term and higher order source terms if both the density $\r$ and the temperature $\t$ have the lower and upper positive bound, which is the main observation of the previous paper \cite{LWYZ} such that the electric potential terms can be estimated and thus the stability of basic wave patterns to the bipolar VPB system \eqref{VPB} could be expected. It should be remarked that this decomposition of the system \eqref{VPB} into the equations \eqref{F1-f-euler}-\eqref{F2-pc-1} is quite universal and give a unified framework for the stability analysis towards the wave patterns.

In the following sections, as an application of the new decomposition,  we prove the
nonlinear stability of the combination of a viscous contact wave and two rarefaction waves to the Cauchy problem of the 1D bipolar VPB system \eqref{VPB0} or \eqref{VPB}. Due to the different structure condition around the viscous contact wave in Eulerian and Lagrangian coordinates and one-dimensional regime,
here we conveniently use the Lagrangian coordinate transformation
\begin{equation}
(t,x)\Longrightarrow\left(t,\int_{(0,0)}^{(t,x)}\rho(\tau,y)dy
-(\rho u_1)(\tau,y)d\tau\right),
\end{equation}
where$\int_A^B fdy+g d\tau$ represents a line integration from point $A$ to point $B$ on the half-plane $\mb{R}^+\times\mb{R}$. Here, the line integration is independent of the path due to the mass conservation equation $\eqref{F1-f-euler}_1$.

Under this Lagrangian transformation, \eqref{VPB0} becomes
\begin{equation}
\left\{
\begin{array}{ll}
\di
F_{At}-\frac{u_1}{v}F_{Ax} +\frac{\xi_1}{v} F_{Ax}+\frac{\Phi_x}{v}\partial_{\xi_1}F_A= Q(F_A, F_A)+Q(F_A,F_B), \\[2mm]
\di
F_{Bt}-\frac{u_1}{v}F_{Bx} +\frac{\xi_1}{v} F_{Bx}-\frac{\Phi_x}{v}\partial_{\xi_1}F_B= Q(F_B, F_A)+Q(F_B,F_B), \\
\di \frac{1}{v}\left(\frac{\Phi_x}{v}\right)_x=\int (F_A-F_B)d\xi
\end{array}
\right.
\end{equation}
with the initial values and the far-field states given by
\begin{equation}\label{VPB10-i}
\begin{array}{ll}
\di F_A(t=0,x,\xi)=F_{A0}(x,\xi)\rightarrow \mathbf{M}_{[v_\pm,u_\pm,\t_\pm]}(\xi),\quad {\rm as}~~x\rightarrow\pm\infty,\\
\di F_B(t=0,x,\xi)=F_{B0}(x,\xi)\rightarrow \mathbf{M}_{[v_\pm,u_\pm,\t_\pm]}(\xi),\quad {\rm as}~~x\rightarrow\pm\infty,\\
\di \Phi_x\rightarrow 0,\quad {\rm as}~~x\rightarrow\pm\infty,
\end{array}
\end{equation}
where the two states $(v_\pm,u_\pm,\t_\pm)$ with $u_\pm=(u_{1\pm},0,0)^t$ and $v_\pm>0, u_{1\pm}, \t_\pm>0$ are connected by the Riemann solution with the combination of two rarefaction waves in the first and third characteristic fields and a contact discontinuity in the second characteristic field to the corresponding 1D Euler system
\begin{equation}\label{euler}
\left\{
\begin{array}{l}
\di v_t-u_{1x}=0,\\
\di u_{1t}+p_x=0, \quad  u_{it}=0,~i=2,3,\\
\di \left(e+\f{|u|^2}{2}\right)_t+(pu_1)_x=0
\end{array}
\right.
\end{equation}
with the Riemann initial data
\begin{equation}\label{R-in}
(v_0,u_0,\t_0)(x)=\left\{
\begin{array}{ll}
\di (v_+,u_+,\t_+), &\di x>0, \\
\di (v_-,u_-,\t_-), &\di x<0.
\end{array}
 \right.
\end{equation}
Furthermore, the transformed system \eqref{VPB} becomes
\begin{equation}\label{vpb-l}
\left\{
\begin{array}{ll}
\di
F_{1t}-\frac{u_1}{v}F_{1x} +\frac{\xi_1}{v} F_{1x}+\frac{\Phi_x}{v}\partial_{\xi_1}F_2= 2Q(F_1, F_1), \\[2mm]
\di
F_{2t}-\frac{u_1}{v}F_{2x} +\frac{\xi_1}{v} F_{2x}+\frac{\Phi_x}{v}\partial_{\xi_1}F_1= 2Q(F_2, F_1), \\[2mm]
\di \frac{1}{v}\left(\frac{\Phi_x}{v}\right)_x=2\int F_2d\xi=2n_2.
\end{array}
\right.
\end{equation}
with the initial values and the far-field states given by
\begin{equation}  \label{vpb-l-data}
\left\{\begin{aligned}
&\di F_1(t=0,x,\xi)=F_{10}(x,\xi)\rightarrow \mathbf{M}_{[v_\pm,u_\pm,\t_\pm]}(\xi),\quad {\rm as}~~x\rightarrow\pm\infty,\\
&\di F_2(t=0,x,\xi)=F_{20}(x,\xi)\rightarrow 0,\quad {\rm as}~~x\rightarrow\pm\infty,\\
&\di \Phi_x\rightarrow 0,\quad {\rm as}~~x\rightarrow\pm\infty.
\end{aligned}
\right.
\end{equation}
By the previous micro-macro decomposition for $F_1=\mb{M}+\mb{G}$, one has the following
fluid-part system for $F_1$:
\begin{equation}\label{ns-l}
\left\{
\begin{array}{lll}
\di v_{t}-u_{1x}=0, \\
\di u_t+p_x-\Phi_x n_2=-\int \xi_1^2\mb{G}_xd\xi,  \\
\di u_{it}=-\int \xi_1\xi_i\mb{G}_xd\xi,  \quad i=2,3,\\[2mm]
\di \left(e+\f{|u|^2}{2}\right)_t+(pu_1)_x-\Phi_x\int\xi_1 F_2 d\xi=-\int\f12\xi_1|\xi|^2\mb{G}_xd\xi
\end{array}
\right.
\end{equation}
and
\begin{equation}\label{F1-f}
\left\{
\begin{array}{l}
\di v_{t}-u_{1x}=0,\\
\di u_{1t}+p_x-\Phi_x n_2 =\f{4}{3}\left(\frac{\mu(\t)}{v}u_{1x}\right)_x-\int \xi_1^2\Pi_{1x}d \xi,  \\
\di u_{it}=\left(\frac{\mu(\t)}{v}u_{ix}\right)_x-\int \xi_1\xi_i\Pi_{1x}d \xi, \quad i=2,3, \\
\di \left(e+\f{|u|^2}{2}\right)_t+(pu_1)_x-\Phi_x\int \xi_1F_2d\xi=\left(\frac{\k(\t)}{v}\t_x\right)_x+\f{4}{3}\left(\frac{\mu(\t)}{v}u_1u_{1x}\right)_x\\
\di\qquad\qquad +\sum_{i=2}^3\left(\frac{\mu(\t)}{v}u_iu_{ix}\right)_x-\int\f12 \xi_1 |\xi|^2\Pi_{1x}d \xi,
\end{array}
\right.
\end{equation}
and the non-fluid equation for $F_1$:
\begin{equation}
\mb{G}_t-\frac{u_1}{v}\mb{G}_x+\f1v\mb{P}_1( \xi_1\mb{M}_x)+\f1v\mb{P}_1(\xi_1\mb{G}_x)
+\f1v\mb{P}_1(\Phi_x\partial_{\xi_1}F_2)
=\mb{L}_\mb{M}\mb{G}+2Q(\mb{G}, \mb{G}), \label{F1-G}
\end{equation}
with
\begin{equation}
\mb{G}= \mb{L}_\mb{M}^{-1}[\f1v\mb{P}_1( \xi_1\mb{M}_x)] +\Pi_1, \label{G1}
\end{equation}
and
\begin{equation}
\Pi_1=\mb{L}_\mb{M}^{-1}[\mb{G}_t-\f{u_1}{v}\mb{G}_x+\f1v\mb {P}_1(\xi_1\mb{G}_x)
+\f1v\mb{P}_1(\Phi_x\partial_{\xi_1}F_2) -2Q(\mb{G},\mb{G})].\label{Pi-1}
\end{equation}
By the new micro-macro decomposition to $F_2$,
$$
F_2=n_2 v\mb{M}+\mb{P}_c F_2,
$$
the macroscopic part $n_2$ satisfy the following equation
\begin{equation}\label{n21}
n_{2t}-\f{u_1}{v}n_{2x}+\f1v\Big(\int \xi_1F_2d\xi\Big)_x=0,
\end{equation}
or
\begin{equation}\label{n22}
n_{2t}+\f{u_{1x}}{v}n_{2}+\f1v\Big(\int \xi_1\mb{P}_c F_2d\xi\Big)_x=0,
\end{equation}
or equivalently,
\begin{equation}\label{n2}
\begin{array}{ll}
\di n_{2t}+\f{u_{1x}}{v}n_2+ \f1v \Big(\f{\k_1(\t)}{R\t v}\Phi_x\Big)_x- \f1v\Big(\f{\k_1(\t)}{v}(n_2 v)_x\Big)_x\\
\di
 =-\f1v\Big(n_2\int \xi_1 \mb{N}_{\mb{M}}^{-1}\big[\mb{P}_c(\xi_{1}\mb{M}_x) \big]d\xi\Big)_x
 -\f1v\Big(\int \xi_1 \mb{N}_{\mb{M}}^{-1}\big[\f1v\mb{P}_c(\xi_1 (\mb{P}_cF_2)_x)\big]d\xi\Big)_x
 - \f1v\Big(\int \xi_1 \mb{N}_{\mb{M}}^{-1}\big[\f{\Phi_x}{v}\mb{G}_{\xi_1}\big]d\xi\Big)_x \\
\di
 -\f1v\Big(\int \xi_1 \mb{N}_{\mb{M}}^{-1}\Big[(\mb{P}_c F_2)_t-\f{u_1}{v}(\mb{P}_c F_2)_x+\big((\mb{M}v)_t-\f{u_1}{v}(\mb{M}v)_x\big)~n_2
 -2Q(F_2,\mb{G})\Big]d\xi\Big)_x.
\end{array}
\end{equation}
And the microscopic part $\mb{P}_c F_2$ satisfies the equation
\begin{equation}\label{F2-pc-L}
(\mb{P}_c F_2)_t-\mb{N}_{\mb{M}}(\mb{P}_c F_2)=-\f{\xi_1}{v}F_{2x}+\f{u_1}{v}F_{2x}-(n_2\mb{M}v)_t
-\f1v \mb{P}_c (\Phi_x\partial_{\xi_1}F_1)+2Q(F_2,\mb{G}).
\end{equation}

%
%
\section{Constructions of Ansatz and Main Result}
\label{result}
\setcounter{equation}{0}

In this section, we first construct the viscous contact wave, rarefaction waves and superposition waves to the bipolar VPB system and then state our main result.

\subsection{Construction of superposition wave}
In the present paper, we consider the case of the superposition  of a viscous
contact wave in the second characteristic field and two rarefaction waves in the first and third characteristic fields, respectively, that is, for given constant state $(v_-, u_-, \t_-)$ with $v_->0$, $\t_->0$ and $u_-\in\mb{R}$, the right end state $(v_+, u_+, \t_+)$ satisfying
\begin{equation}\label{zheng}
\di (v_+,u_+,\t_+)\in R_1-CD_2-R_3(v_-,u_-,\t_-)\subset\Omega(v_-,u_-,\t_-),
\end{equation}
where $\Omega(v_-,u_-,\t_-)$ is some neighborhood of the state $(v_-, u_-, \t_-)$ in the phase plane $\{(v,u,\t)|v>0,u>0\}$  and
\begin{equation*}
\begin{array}{ll}
\di R_1-CD_2-R_3(v_-,u_-,\t_-)\triangleq\Bigg\{(v,u,\t)\in\Omega(v_-,u_-,\t_-)\Bigg|s\neq s_-,\\
\di u\geq u_--\int_{v_-}^{e^{\f35(s_--s)}v}\lambda_-(\eta,s_-)d\eta,
u\geq u_--\int^v_{e^{\f35(s_--s)}v_-}\lambda_+(\eta,s)d\eta\Bigg\}
\end{array}
\end{equation*}
and the entropy $s$ and the eigenvalues $\lambda_{\pm}(v,s)$ are given by
$$
s=\f23\ln v+\ln(\f{4\pi}{3}\t)+1,\quad s_{\pm}=\f23\ln v_{\pm}+\ln(\f{4\pi}{3}\t_{\pm})+1,~~
\lambda_{\pm}(v,s)=\pm\sqrt{\f{5}{6\pi} v^{-\f83}e^{s}}.
$$
By the standard argument (e.g. \cite{Smoller}), if the wave strength $\delta$ is suitably small, there exists a unique pair of states $(v_-^*, u^*, \t_-^*)$
and $(v_+^*, u^*, \t_+^*)$ in $\Omega(v_-,u_-,\t_-)$ satisfying
\begin{equation}\label{RH}
\frac{2\t_-^*}{3v_-^*}=\frac{2\t_+^*}{3v_+^*}\triangleq p^*,
\end{equation}
i. e., the states $(v_\pm^*, u^*, \theta_\pm^*)$ are connected by 2-contact discontinuity solution
\begin{equation}\label{cdi}
(\tilde V,\tilde U,\tilde\Theta)(x,t)=\left\{
\begin{array}{ll}
\di (v_-^*,u^*,\t_-^*),&\di x<u^*t,t>0,\\
\di (v_+^*,u^*,\t_+^*),&\di x>u^*t,t>0,
\end{array}
\right.
\end{equation}
to the Riemann problem  \eqref{euler}-\eqref{R-in} provide that \eqref{RH} holds, and the states $(v_-^*,u^*,\t_-^*)$ and $(v_+^*,u^*,\t_+^*)$ belong to the 1-rarefaction wave
curve $R_-(v_-,u_-,\t_-)$ and the 3-rarefaction wave curve $R_+(v_+,u_+,\t_+)$, respectively, with the rarefaction curves defined by
$$
R_{\pm}(v_{\pm},u_{\pm},\t_{\pm})=\left\{(v,u,\t)\Bigg|s=s_{\pm},
u=u_{\pm}-\int_{v_{\pm}}^v\lambda_{\pm}(\eta,s_{\pm})d\eta,v>v_{\pm}
\right\}.
$$
Now we first construct the viscous contact wave, which is the viscous version of 2-contact discontinuity solution \eqref{cdi}. Without loss of generality, we assume $u^*=0$ in what follows.
Motivated by \cite{Huang-Yang},  for the bipolar VPB system, the viscous contact wave is constructed as the unique self-similar solution $\Theta(\f{x}{\sqrt{1+t}})$ of the following nonlinear diffusion equation
\begin{equation}\label{1.8}
\di \Theta_t=(a(\Theta)\Theta_x)_x,\quad
\Theta(\pm\infty,t)=\theta_{\pm}^*,\quad a=\frac{9p^*\k(s)}{10s}>0.
\end{equation}
We define
\begin{equation}\label{contact}
\di V^{cd}=\frac{2}{3p^*}\Theta,\quad  U^{cd}_1=\frac{2a(\Theta)}{3p^*}\Theta_x+u^*, \quad U^{cd}_i=0,~i=2,3,
\quad  \Theta^{cd}=\Theta.
\end{equation}
There exists some positive constant $\bar\delta$, such that for $\d^{cd}=|\t_+^* -\t_-^*|\leqq\bar \d$, $\Theta$ satisfies
\begin{equation}\label{1.9}
\di (1+t)|\Theta_{xx}|+(1+t)^{\frac{1}{2}}|\Theta_x|+|\Theta-\theta_{\pm}^*|
\leq c_1\delta^{cd} e^{-\frac{c_2x^2}{1+t}} \quad \mathrm{as}~ |x|\rightarrow\infty,
\end{equation}
where $c_1$ and $c_2$ are positive constants depending only on $\t_{\pm}^*$ and $\bar \d$. It can verified by a straightforward
computation that $(V^{cd},  U^{cd},  \Theta^{cd})$ satisfies
\begin{equation}\label{contact equation}
\left\{
\begin{array}{ll}
\di V^{cd}_t- U^{cd}_{1x}=0,\\
\di U^{cd}_{1t}+ P^{cd}_x=\f43\left(\f{\mu(\Theta^{cd})}{V^{cd}}U^{cd}_{1x}\right)_x+R_1,
\quad \di U^{cd}_{it}=0,\quad i=2,3,\\
\di \Theta^{cd}_t + P^{cd}U^{cd}_{1x}=\left(\f{\k(\Theta^{cd})}{V^{cd}}\Theta^{cd}_x\right)_x+\frac{4\mu(\Theta^{cd})}{3V^{cd}}(U^{cd}_{1x})^2
+\sum_{i=2}^3\frac{\mu(\Theta^{cd})}{V^{cd}}(U^{cd}_{ix})^2+R_2,
\end{array}
\right.
\end{equation}
where $ P^{cd}=p_+^*=p_-^*=p^*$ and
\begin{equation}\label{r1r2}
\begin{array}{l}
\di  R_1=U^{cd}_{1t}-\frac43\left(\frac{\mu(\Theta^{cd})}{V^{cd}}U^{cd}_{1x}\right)_x=O(1)\d^{cd}(1+t)^{-3/2}e^{-\f{c_2x^2}{1+t}},\\
\di  R_2=-\left(\f43\frac{\mu(\Theta^{cd})}{V^{cd}} (U^{cd}_{1x})^2+\sum_{i=2}^3\frac{\mu(\Theta^{cd})}{V^{cd}}(U^{cd}_{ix})^2\right)
=O(1)\d^{cd}(1+t)^{-2}e^{-\f{c_2x^2}{1+t}}.\\
\end{array}
\end{equation}

Now we turn to the construction of approximate rarefaction wave. The 1-rarefaction wave $(v_-^r,u_-^r,\t_-^r)(\frac{x}{t})$ (respectively 3-rarefaction wave $(v_+^r,u_+^r,\t_+^r)(\frac{x}{t})$)
connecting $(v_-,u_-,\t_-)$ and $(v_-^*,0,\t_-^*)$ (respectively connecting $(v_+^*,0,\t_+^*)$ and $(v_+,u_+,\t_+)$) is the
weak solution of the Riemann problem of the Euler system \eqref{euler} with the following initial Riemann data
\begin{equation}\label{Riemann2}
\di(v_{\pm},u_{\pm},\t_{\pm})(x,0)=\left\{
\begin{array}{ll}
(v_{\pm}^*,0,\t_{\pm}^*), &\quad\pm x<0,\\
(v_{\pm},u_{\pm},\t_{\pm}), &\quad\pm x>0.
\end{array}
\right.
\end{equation}
Since the rarefaction wave $(v_{\pm}^r,u_{\pm}^r,\t_{\pm}^r)$ are weak solutions, it is
convenient to construct approximate rarefaction wave which is smooth. Motivated by Matsumura-Nishihara \cite{MN-86}, the smooth solutions of Euler
system \eqref{euler}, $(V_{\pm}^r,U_{\pm}^r,\Theta_{\pm}^r)$, which approximate $(v_{\pm}^r,u_{\pm}^r,\t_{\pm}^r)$,
are given by
\begin{equation}\label{appro-rare}
\di\left\{
\begin{array}{ll}
\di\lambda_{\pm}(V_{\pm}^r(x,t),s_{\pm})=w_{\pm}(x,t),\\
\di U_{1\pm}^r=u_{1\pm}-\int_{v_{\pm}}^{V_{\pm}^r(x,t)}\lambda_{\pm}(\eta,s_{\pm})d\eta,\quad U^{r}_{i\pm}=0,~i=2,3,\\
\Theta_{\pm}^r=\t_{\pm}(v_{\pm})^{\f23}(V_{\pm}^r)^{-\f23},
\end{array}
\right.
\end{equation}
 where $w_-$ (respectively $w_+$) is the solution of the initial problem for the typical Burgers equation:
\begin{equation}\label{burgers}
\di\left\{
\begin{array}{ll}
\di w_t+ww_x=0,\quad(x,t)\in\mathbb{R}\times(0,\infty),\\
\di w(x,0)=\frac{w_r+w_l}{2}+\frac{w_r-w_l}{2}\tanh x,
\end{array}
\right.
\end{equation}
with $w_l=\lambda_-(v_-,s_-)$, $w_r=\lambda_-(v_-^*,s_-)$ (respectively $w_l=\lambda_+(v_+^*,s_+)$, $w_r=\lambda_+(v_+,s_+)$).
It can be verified by a straightforward computation that $(V^r_\pm, U^r_\pm, \Theta^r_{\pm})$ satisfies
\begin{equation}\label{rare-equation}
\di\left\{
\begin{array}{ll}
\di(V^r_\pm)_t-(U^r_{1\pm})_x=0,\\
\di (U^r_{1\pm})_t+(P^r_\pm)_x=0, \quad (U^{r}_{i\pm})_t=0,~i=2,3,\\
\di (\Theta^r_{\pm})_t+P^r_\pm(U^r_{1\pm})_x=0,
\end{array}
\right.
\end{equation}
where $P^r_\pm=\frac{2\Theta^r_{\pm}}{3V^r_{\pm}}$.
\begin{lemma}\label{bur-pro}
For given $w_l\in\mb{R}$ and $\bar w>0$, let $w_r\in\{0<\tilde w\triangleq w-w_l<\bar w\}$.
Then the problem \eqref{burgers} has a unique smooth global solution in time satisfying the following properties.
\begin{itemize}
\item[(i)] $w_l<w(x,t)<w_r$, $w_x>0$ $(x\in\mb{R},t>0)$.
\item[(ii)] For $p\in[1,\infty]$, there exists some positive constant $C=C(p,w_l,\bar w)$ such that
for $\tilde w\geqq0$ and $t\geqq0$,
$$
\|w_x(t)\|_{L^p}\leq C\min\{\tilde w,\tilde w^{1/p}t^{-1+1/p}\},\quad
\|w_{xx}(t)\|_{L^p}\leq C\min\{\tilde w,t^{-1}\}.
$$
\item[(iii)] If $w_l>0$, for any $(x,t)\in(-\infty,0]\times[0,\infty)$,
$$
|w(x,t)-w_l|\leq\tilde we^{-2(|x|+w_lt)},\quad
|w_x(x,t)|\leq2\tilde we^{-2(|x|+w_lt)}.
$$
\item[(iv)] If $w_r<0$, for any $(x,t)\in[0,\infty)\times[0,\infty)$,
$$
|w(x,t)-w_r|\leq\tilde we^{-2(x+|w_r|t)},\quad
|w_x(x,t)|\leq2\tilde we^{-2(x+|w_r|t)}.
$$
\item[(v)] For the Riemann solution $w^r(x/t)$ of the scalar equation \eqref{burgers} with the Riemann
initial data
\begin{equation*}
w(x,0)=\left\{
\begin{array}{ll}
w_l,&\quad x<0,\\
w_r,&\quad x>0,
\end{array}
\right.
\end{equation*}
we have
$$
\lim_{t\rightarrow+\infty}\sup_{x\in\mb{R}}|w(x,t)-w^r(x/t)|=0.
$$
\end{itemize}

\end{lemma}

In order to cope with the wave interactions of viscous contact wave and two rarefaction waves, we divide the half space $\mb{R}\times\mb{R}^+$ into
three parts, that is $\mb{R}\times\mb{R}^+=\Omega_-\cup\Omega_c\cup\Omega_+$ with
$$
\Omega_{\pm}=\big\{(x,t)\big|\pm2x>\pm\lambda_{\pm}(v_{\pm}^*,s_{\pm})t\big\},
$$
and
$$
\Omega_{c}=\big\{(x,t)\big|\lambda_-(v_-^*,s_-)t\leq2x\leq\lambda_{+}(v_{+}^*,s_{+})t\big\}.
$$
Then \eqref{1.9} and Lemma \ref{bur-pro}  lead to

\begin{lemma}\label{rare-pro}
For any given left end state $(v_-,u_-,\t_-)$, we assume that \eqref{zheng} holds.
Then the smooth rarefaction waves $(V_{\pm}^r,U_{\pm}^r,\Theta_{\pm}^r)$
constructed in \eqref{appro-rare} and the viscous contact wave $(V^{cd},U^{cd},\Theta^{cd})$ constructed in \eqref{contact} satisfying
the following:
\begin{itemize}
\item[(i)] $(U_{1\pm}^r)_x\geq0$,~$(x\in\mb{R},t>0)$.
\item[(ii)] For $p\in[1,\infty]$, there exists a positive constant $C=C(v_-,u_-,\t_-,\d)$
such that for $\d=|\t_+-\t_-|$,
$$
\|\big((V_{\pm}^r)_x,(U_{1\pm}^r)_x,(\Theta_{\pm}^r)_x\big)(t)\|_{L^p}
\leq C\min\Big\{\d,\d^{1/p}t^{-1+1/p}\Big\}
$$
and
$$
\|\big((V_{\pm}^r)_{xx},(U_{1\pm}^r)_{xx},(\Theta_{\pm}^r)_{xx}\big)(t)\|_{L^p}
\leq C\min\Big\{\d,t^{-1}\Big\}.
$$
\item[(iii)] There exists a positive constant $C=C(v_-,u_-,\t_-,\d)$
such that for
$$
c_0=\frac{1}{10}\min\Big\{|\lambda_-(v_-^*,s_-)|,\lambda_+(v_+^*,s_+),c_1\lambda_-^2(v_-^*,s_-),
c_1\lambda_+^2(v_+^*,s_+),1\Big\},
$$
we have in $\Omega_c$
$$
(U_{1\pm}^r)_x+|(V_{\pm}^r)_x|+|V_{\pm}^r-v_{\pm}^*|+|(\Theta_{\pm}^r)_x|+|\Theta_{\pm}^r-\t_{\pm}^*|
\leq C\d e^{-c_0(|x|+t)},
$$
and in $\Omega_{\mp}$
$$
\left\{
\begin{array}{ll}
|V^{cd}-v_{\mp}^*|+|V^{cd}_x|+|\Theta^{cd}-\t_{\mp}^*|+|U^{cd}_{1x}|+|\Theta^{cd}_x|
\leq C\d e^{-c_0(|x|+t)},\\
(U_{1\pm}^r)_x+|(V_{\pm}^r)_x|+|V_{\pm}^r-v_{\pm}^*|+|(\Theta_{\pm}^r)_x|+|\Theta_{\pm}^r-\t_{\pm}^*|
\leq C\d e^{-c_0(|x|+t)}.
\end{array}
\right.
$$
\item[(iv)] For the rarefaction waves $(v_{\pm}^r,u_{\pm}^r,\t_{\pm}^r)(x/t)$, it holds that
$$
\lim_{t\rightarrow+\infty}\sup_{x\in\mb{R}}
\big|(V_{\pm}^r,U_{\pm}^r,\Theta_{\pm}^r,)(x,t)-(v_{\pm}^r,u_{\pm}^r,\t_{\pm}^r,)(x/t)\big|=0.
$$
\end{itemize}
\end{lemma}

Finally, we define the superposition wave pattern by
\begin{equation}\label{ansatz}
\left(
\begin{array}{l}
\di \bar v\\
\di \bar u_1 \\
\di \bar \t
\end{array}
\right)(x,t) = \left(
\begin{array}{l}
\di V^{cd}+V_-^r+V_+^r\\
\di U^{cd}_1+U_{1-}^r+U_{1+}^r\\
\di \Theta^{cd}+\Theta_-^r+\Theta_+^r
\end{array}\right)(x,t)-\left(
\begin{array}{l}
\di v_-^*+v_+^*\\
\di \qquad 0,\\
\di \t_-^*+\t_+^*
\end{array}\right), \quad \bar u_i=0, i=2,3,
\end{equation}
which is the linear combination of  a viscous contact wave and two rarefaction waves.
Then the composite wave pattern $(\bar v,\bar u, \bar\t)$ satisfies
\begin{equation}\label{ansatz-system}
\di\left\{
\begin{array}{ll}
\di\bar v_t-\bar u_{1x}=0,\\
\di \bar u_{1t}+\bar p_x=(\bar p-p_+-p_-)_x+\frac43\Big(\frac{\mu(\Theta^{cd})}{V^{cd}}U^{cd}_{1x}\Big)_x+R_1,\\
\di \bar u_{it}=0,\\
\di \bar\t_t+\bar p \bar u_{1x}=(\bar p- p^*)U^{cd}_{1x}+(\bar p-p^r_+)(U^r_{1+})_x+(\bar p-p^r_-)(U^r_{1-})_x\\
\di\qquad +\Big(\frac{\k(\Theta^{cd})}{V^{cd}}\Theta^{cd}_x\Big)_x+\frac{4\mu(\Theta^{cd})}{3V^{cd}}(U^{cd}_{1x})^2
+\sum_{i=2}^3\frac{\mu(\Theta^{cd})}{V^{cd}}(U^{cd}_{ix})^2+R_2,
\end{array}
\right.
\end{equation}
where $\bar p=\frac{2\bar\t}{3\bar v}$, and $R_1$ and $R_2$ are defined as in \eqref{r1r2}.

\subsection{Main Result.}
Denote
\begin{equation}\label{Et}
\begin{array}{ll}
\di \mathcal{E}(t)=
\di \sup_{\tau\in[0,t]}\bigg\{\|(v-\bar v,u-\bar u,\t-\bar \t)\|^2_{H^1(\mb{R})}+\|(\Phi_{x},n_{2},n_{2x})\|^2
+\sum_{0\leq|\beta|\leq2}\int_{\mb{R}}\int\f{|\partial^\beta(\mb{G},\mb{P}_cF_{2})|^2}{\mb{M}_\star }d\xi
dx\\
\di\qquad\quad +\sum_{|\alpha'|=1,0\leq|\beta'|\leq1}\int\int\f{|\partial^{\alpha'}\partial^{\beta'}(\mb{G},\mb{P}_cF_{2})|^2}{\mb{M}_\star }d\xi
dx
 +\sum_{|\alpha|=2}\int\int\f{|\partial^\alpha(F_{1},F_{2})|^2}{\mb{M}_\star }d\xi dx\bigg\},
\end{array}
\end{equation}
where and in the sequel $\partial^\alpha=\partial_{x,t}^\alpha$, $\partial^\beta=\partial_{\xi}^\beta$. Then, we can state our main result as follows.

\begin{theorem}\label{thm}
There exist positive constants $\delta_0$ and $\varepsilon_0$ and a global Maxellian $\mb{M}_\star =\mb{M}_{[v_\star,u_\star,\t_\star]}$ with $v_\star>0,\t_\star>0$, such that if
the wave strength $\delta=|(v_+-v_-,u_+-u_-,\t_+-\t_-)| \leq \delta_0$ and the initial data satisfies
\begin{equation}
\mathcal{E}(0)\leq \varepsilon_0,
\end{equation}
then the Cauchy problem of the bipolar VPB system \eqref{vpb-l}--\eqref{vpb-l-data} admits a unique global classical solution $(F_1,F_2,\Phi)$ which satisfies the uniform estimates
$$
\mathcal{E}(t)\leq C(\mathcal{E}(0)+\delta_0^{\frac12})
$$
for the  positive constant $C$ independent of time, and tends time-asymptotically towards the composite wave pattern $(\bar v,\bar u,\bar \t)$
$$
\begin{array}{ll}
\di \|\big(F_1(t,x,\xi)-\mb{M}_{[\bar v,\bar u,\bar \t]}(t,x,\xi), F_2(t,x,\xi)\big)\|_{L^\infty_x
L_{\xi}^2(\frac1{\sqrt{\mb{M_\star}}})}+ \|(\Phi_x,n_2)(t,x)\|_{L^\infty_x}
\rightarrow
0, ~~~~{\rm as}~~t\to \infty,.
\end{array}
$$
Consequently, it holds that
$$
\|\big(F_A-\mb{M}_{[\bar v,\bar u,\bar \t]}, F_B-\mb{M}_{[\bar v,\bar u,\bar \t]}\big)\|_{L^\infty_x
L_v^2(\frac1{\sqrt\mb{M_*}})}+ \|(\Phi_x ,n_2)\|_{L^\infty_x}
\rightarrow
0, ~~~~{\rm as}~~t\to \infty.
$$
Here and in the sequel $f(\xi)\in L_{\xi}^2(\f{1}{\sqrt{\mb{M}_\star }})$ means that
$\frac{f(\xi)}{\sqrt{\mb{M_\star}}}\in L_{\xi}^2(\mb{R}^3)$.
\end{theorem}

\begin{remark}
Theorem \ref{thm} implies that the composite wave pattern by the linear superposition of 2-viscous contact wave and two rarefaction waves in the first and third characteristic fields is nonlinearly stable to the 1D bipolar VPB system \eqref{VPB} under the combined effects of the binary collisions, mutual interactions, and the electrostatic potential forces.
\end{remark}

Set the perturbation around the superposition wave pattern $(\bar v,\bar u, \bar\t)(x,t)$ by
\begin{equation}\label{perturb}
(\phi,\psi,\zeta)(x,t)=(v-\bar v, u-\bar u,\t-\bar\t)(x,t)
\end{equation}
where $(v, u, \t)(x,t)$ is the fluid variables to the solution $F_1(x,t,\xi)$ of the VPB equation \eqref{vpb-l}.
By \eqref{F1-f} and \eqref{ansatz-system}, we can write the system for the perturbation $(\phi,\psi,\zeta)$ in \eqref{perturb} as follows.
\begin{equation}\label{sys-h}
\left\{
\begin{array}{ll}
\di \phi_t-\psi_{1x}=0,\\[2mm]
\di \psi_{1t}+(p-\bar p)_x-\Phi_x n_2
=\f43\left(\frac{\mu(\t)}{v}u_{1x}-\frac{\mu(\bar\t)}{\bar v} \bar u_{1x}\right)_x-\int \xi_1^2\Pi_{1x} d\xi-\bar R_1,\\[2mm]
\di \psi_{it}=\left(\frac{\mu(\t)}{v}u_{ix}-\frac{\mu(\bar\t)}{\bar v} \bar u_{ix}\right)_x-\int \xi_1\xi_i\Pi_{1x} d\xi, ~i=2,3,\\[2mm]
\di \zeta_{t}+(p u_{1x}-\bar p\bar u_{1x})-\Phi_x\int\xi_1\mb{P}_cF_2d\xi
=\left(\frac{\k(\t)}{v}\t_{x}-\frac{\k(\bar\t)}{\bar v} \bar \t_{x}\right)_x
+\f43\left(\frac{\mu(\t)}{v}u^2_{1x}-\frac{\mu(\bar\t)}{\bar v} \bar u^2_{1x}\right)\\[2mm]
\di\quad+\sum_{i=2}^3\left(\frac{\mu(\t)}{v}u^2_{ix}-\frac{\mu(\bar\t)}{\bar v} \bar u^2_{ix}\right)
-\int \f12 \xi_1|\xi|^2\Pi_{1x}d\xi+\sum_{i=1}^3 u_i\int \xi_1\xi_i\Pi_{1x}d\xi-\bar R_2,
\end{array} \right.
\end{equation}
where
$$
\begin{array}{l}
\di \bar p=\f{2\bar \t}{3\bar v},\quad p^*=\frac{2\t_-^*}{3v_-^*}=\frac{2\t_+^*}{3v_+^*} ,\quad  P^r_\pm=\f{2\Theta^r_{\pm}}{3 V^{r}_{\pm}},\\
\di \bar R_1=(\bar p-P^r_--P^r_+)_x+U^{cd}_{1t}-\frac43\left(\frac{\mu(\bar\t)}{\bar v}\bar u_{1x}\right)_x,\\
\di \bar R_2=(\bar p-p^*)U^{cd}_{1x}+(\bar p-P^r_-)(U^{r}_{1-})_x+(\bar p-P^r_+)(U^{r}_{1+})_x-\f43\frac{\mu(\bar\t)}{\bar v}|\bar u_{1x}|^2
-\left(\k(\bar\t)\f{\bar\t_x}{\bar v}-\k(\Theta^{cd})\f{\Theta^{cd}_x}{V^{cd}}\right)_x.
\end{array}
$$

Denote
\begin{equation}\label{bG}
\bar{\mb{G}}=\f{3}{2v\t }\mb{L}_{\mb{M}}^{-1}\Big[\mb{P}_1\big(\xi_1(\f{|\xi-u|^2}{2\t}\bar\t_x+\xi_1\bar
u_{1x})\big)\mb{M}\Big],
\end{equation}
and let
\begin{equation}\label{tG}
\widetilde{\mb{G}}=\mb{G}-\bar{\mb{G}}.
\end{equation}
Note that for technical reasons, we introduce $\bar{\mb{G}}$ in \eqref{bG} to circumvent the difficulty caused by
the fact that $\|\bar\t_x\|^2\sim(1+t)^{-\f12}$ is not integrable with respect to $t$.
Then it holds that
\begin{equation}\label{GE}
\begin{array}{ll}
\di \widetilde{\mb{G}}_t-\mb{L}_{\mb{M}}\widetilde{\mb{G}}=
-\f{3}{2v\t}\mb{P}_1\Big[\xi_1(\f{|\xi-u|^2}{2\t}\zeta_x+\xi\cdot\psi_{x})\Big]\mb{M}+\f{u_1}{v}\mb{G}_x
-\f1v \mb{P}_1(\xi_1\mb{G}_x)-\f1v \mb{P}_1(\Phi_x\partial_{\xi_1}F_2)+2Q(\mb{G},\mb{G})-\bar{\mb{G}}_t.
\end{array}
\end{equation}

We do the following a priori assumption:
\begin{equation}\label{assumption-ap}
\begin{array}{ll}
\di \mathcal{N}(T)\di=\sup_{0\leq t \leq
T}\Bigg\{
\|(\phi,\psi,\zeta)(t,\cdot)\|^2_{H^1}+\|(\Phi_x,n_2,n_{2x})\|^2+\sum_{0\leq|\beta|\leq 2}\int\int\f{|\partial^\beta(\wt{\mb{G}},\mb{P}_cF_2)|^2}{\mb{M}_\star }d\xi dx\\
\di \qquad\quad +\sum_{|\alpha^\prime|=1,0\leq|\beta^\prime|\leq1}\int\int\f{|\partial^{\alpha^\prime}\partial^{\beta^\prime}(\mb{G},\mb{P}_cF_2)|^2}{\mb{M}_\star }d\xi dx
+\sum_{|\alpha|=2}\int\int\f{|\partial^\alpha(F_1,F_2)|^2}{\mb{M}_\star }d\xi dx\Bigg\}\leq \chi^2,
\end{array}
\end{equation}
where and in the sequel $\chi$ is a small positive constant depending on the initial data
and wave strengths.

By \eqref{ns-l}, \eqref{r1r2} and \eqref{ansatz-system}, one has
\begin{equation}\label{sys-h-o}
\left\{
\begin{array}{ll}
\di \phi_t-\psi_{1x}=0,\\[2mm]
\di \psi_{1t}+(p-\bar p)_x-\Phi_x n_2=-\int\xi_1^2\mb{G}_xd\xi-  U^{cd}_{1t}-(\bar p-P^r_+-P^r_-)_x,\\[2mm]
\di \psi_{it}=-\int \xi_1\xi_i\mb{G}_{x} d\xi, ~i=2,3,\\[2mm]
\di \zeta_{t}+(p u_{1x}-\bar p\bar u_{1x})-\Phi_x\int\xi_1\mb{P}_cF_2d\xi
=-\int \f12 \xi_1|\xi|^2\mb{G}_{x}d\xi+\sum_{i=1}^3 u_i\int \xi_1\xi_i\mb{G}_{x}d\xi\\
\di \quad -\Big((\bar p-p^*)U^{cd}_{1x}+(\bar p-P^r_-)(U^{r}_{1-})_x+(\bar p-P^r_+)(U^{r}_{1+})_x\Big)
-\left(\frac{\k(\Theta^{cd})}{V^{cd}}\Theta^{cd}_x\right)_x.
\end{array} \right.
\end{equation}

In fact, by the a priori assumption \eqref{assumption-ap}, one also has from the system \eqref{sys-h-o} that
\begin{equation}\label{4.18}
\|(\phi,\psi,\zeta,n_2)\|_{L^\i_x}^2\leq C\chi^2,
\end{equation}
and
\begin{equation}
\|(\phi_{t},\psi_{t},\zeta_{t})\|^2\leq C(\chi+\d)^2,
\label{(4.19)}
\end{equation}
hence, one has
\begin{equation}
\|(v_{t},u_{t},\t_{t})\|^2\di
\leq C\|(\phi_{t},\psi_{t},\zeta_{t})\|^2+ C \|(\bar v_{t},\bar u_{t},\bar\t_{t})\|^2\leq C(\chi+\d)^2.
\end{equation}
For $|\alpha|=2$, it follows from \eqref{macro}  and \eqref{assumption-ap} that
\begin{equation}
\|\partial^\a\left(\rho,\rho u,\rho(\t+\f{|u|^2}{2}),n_2\right)\|^2 \leq
C\int\int\f{|\partial^\a (F_1,F_2)|^2}{\mb{M}_\star }d\xi dx\leq C\chi^2,
\label{(4.8)}
\end{equation}
and
\begin{equation}
\begin{array}{ll}
\di\|\partial^\a(v,u,\t)\|^2&\di\leq C\|\partial^\a\left(\rho,\rho
u,\rho(\t+\f{|u|^2}{2})\right)\|^2
+C\sum_{|\a^\prime|=1}\int|\partial^{\a^\prime}\left(\rho,\rho
u,\rho(\t+\f{|u|^2}{2})\right)|^4dx\\
 &\di\leq C(\chi+\delta)^2.
\end{array}
\label{3.25}
\end{equation}
Therefore, for $|\a|=2$, we have
\begin{equation}
\|\partial^\a(\phi,\psi,\zeta,n_2)\|^2\leq C(\chi+\d)^2. \label{3.26}
\end{equation}
By \eqref{n22} and a priori assumption \eqref{assumption-ap}, it holds that
\begin{equation}\label{n2t-ap}
\|n_{2t}\|^2\leq C\|u_{1x}\|^2+C\int\int\f{|(\mb{P}_c F_2)_x|^2}{\mb{M}_\star }d\xi dx\leq C(\chi+\d)^2.
\end{equation}
By \eqref{assumption-ap}, \eqref{3.26} and \eqref{n2t-ap}, for $|\alpha^\prime|=1$, it holds that
\begin{equation}\label{w-3.30}
\begin{array}{ll}
\|\partial^{\alpha^\prime}(\phi,\psi,\zeta,n_2)\|^2_{L^{\i}}\leq
C(\chi+\d)^2.
\end{array}
\end{equation}
By \eqref{vpb-l}$_3$ and \eqref{n22}, one has
\begin{equation}\label{wt-3.31}
\left(\f{\Phi_x}{2v}\right)_t+\int \xi_1 \mb{P}_cF_2d\xi=0.
\end{equation}
Then by $(\ref{vpb-l})_3$ and \eqref{wt-3.31},  one has
\begin{equation}\label{phi-xt-ap}
\begin{array}{l}
\di\|\Phi_{xx}\|^2\leq C\|n_2\|^2+C\int |\Phi_x|^2|v_x|^2dx\leq C(\chi+\d)^2,\\
\di\|\Phi_{xt}\|^2\leq C\int |\Phi_x|^2|v_t|^2dx+C\int\int\frac{|\mb{P}_cF_2|^2}{\mb{M}_\star }d\xi dx\leq C(\chi+\d)^2.
\end{array}
\end{equation}
By \eqref{3.26}, \eqref{w-3.30} and \eqref{phi-xt-ap}, it holds that
\begin{equation}\label{phixtt-ap}
\begin{array}{l}
\di\|(\Phi_{xxx},\Phi_{xxt},\Phi_{xtt})\|^2\leq C\int|(v_x,v_t)|^2|(\Phi_{xx},\Phi_{xt})|^2dx
+C\int|\Phi_x|^2(|(v_{xx},v_{xt},v_{tt})|^2+|(v_x,v_t)|^4)dx\\
\di +C\int(|n_{2x}|^2+|n_2|^2|v_x|^2)dx +C\int|(v_x,v_t)|^2\left(\int\frac{|\mb{P}_cF_2|^2}{\mb{M}_\star }d\xi\right) dx
+C\sum_{|\alpha'|=1}\int\int\frac{|\pa^{\alpha'}(\mb{P}_cF_2)|^2}{\mb{M}_\star }d\xi dx\\
\di\leq C(\chi+\d)^2.
\end{array}
\end{equation}
By \eqref{n2t-ap}, \eqref{phi-xt-ap}, \eqref{phixtt-ap} and \eqref{3.26}, it holds that
\begin{equation}\label{phi-xt-in}
\|(\Phi_x,\Phi_{xx},\Phi_{xt})\|_{L^\infty}^2\leq  C(\chi+\d)^2.
\end{equation}
Moreover, it holds that
\begin{equation}
\begin{array}{ll}
 \di\|\int\f{|(\wt{\mb{G}},\mb{P}_c F_2)|^2}{\mb{M}_\star }d\xi \|_{L^\infty_x} \leq
C\left(\int\int\f{|(\wt{\mb{G}},\mb{P}_c F_2)|^2}{\mb{M}_\star }d\xi
dx\right)^{\f{1}{2}}\cdot\left(\int\int\f{|(\wt{\mb{G}},\mb{P}_c F_2)_x|^2}{\mb{M}_\star }d\xi dx\right)^{\f{1}{2}}
\leq C\chi.
\end{array}
\label{G-infty}
\end{equation}
Furthermore,  for $|\a^\prime|=1$, $|\beta^\prime|=1$,  it holds that
\begin{equation}
\begin{array}{ll}
 \di\|\int\f{|\partial^{\a^\prime}
(\mb{G}, \mb{P}_c F_2)|^2}{\mb{M}_\star }d \xi\|_{L_x^{\i}}\leq
C\left(\int\int\f{|\partial^{\a^\prime} (\mb{G}, \mb{P}_c F_2)|^2}{\mb{M}_\star }d\xi
dx\right)^{\f{1}{2}}\cdot\left(\int\int\f{|\partial^{\a^\prime}
(\mb{G}, \mb{P}_c F_2)_x|^2}{\mb{M}_\star }d\xi dx\right)^{\f{1}{2}}\leq C\chi
\end{array}
\label{(4.17+)}
\end{equation}
and
\begin{equation}
\begin{array}{ll}
 \di\|\int\f{|\partial^{\beta^\prime}
(\widetilde{\mb{G}}, \mb{P}_c F_2)|^2}{\mb{M}_\star }d \xi\|_{L_x^{\i}}\leq
C\left(\int\int\f{|\partial^{\beta^\prime} (\widetilde{\mb{G}}, \mb{P}_c F_2)|^2}{\mb{M}_\star }d\xi
dx\right)^{\f{1}{2}}\cdot\left(\int\int\f{|\partial^{\beta^\prime}
(\widetilde{\mb{G}}, \mb{P}_c F_2)_x|^2}{\mb{M}_\star }d\xi dx\right)^{\f{1}{2}}\leq C\chi.
\end{array}
\label{(4.17-)}
\end{equation}
Finally, by noticing the facts that $F_1=\mb{M}+\mb{G}$ and
$F_2=n_2 v\mb{M}+\mb{P}_c F_2$  and  \eqref{3.25} with $|\a|=2$, it holds that
\begin{equation}
\begin{array}{ll}
\di\int\int\f{|\partial^{\alpha}(\mb{G},\mb{P}_cF_2)|^2}{\mb{M}_\star }d\xi dx\di
\le C\int\int\f{|\partial^\alpha (F_1,F_2)|^2}{\mb{M}_\star }d\xi dx
+C\int\int\f{|\partial^\alpha\mb{M}|^2+|\partial^\a\big(n_2 v\mb{M}\big)|^2}{\mb{M}_\star }d\xi dx
\leq C(\chi+\d)^2,
\end{array}
\end{equation}
where in the last inequality we have used a similar argument as used for
\eqref{(4.17+)}.

Now we list some lemmas for later use. The following lemmas are based on the celebrated H-theorem. The first lemma is from \cite{Liu-Yu}.

\begin{lemma}\label{Lemma 4.1} There exists a positive
constant $C$ such that
$$
\int\f{\nu(|\xi|)^{-1}Q(f,g)^2}{\wt{\mb{M}}}d\xi\leq
C\left\{\int\f{\nu(|\xi|)f^2}{\wt{\mb{M}}}d\xi\cdot\int\f{g^2}{\wt{\mb{M}}}d\xi+
\int\f{f^2}{\wt{\mb{M}}}d\xi\cdot\int\f{\nu(|\xi|)g^2}{\wt{\mb{M}}}d\xi\right\},
$$
where $\wt{\mb{M}}$ can be any Maxwellian so that the above
integrals are well-defined.
\end{lemma}

Based on Lemma \ref{Lemma 4.1}, the following three lemmas are taken
from \cite{Liu-Yang-Yu-Zhao}. Their proofs are straightforward by
using Cauchy inequality.

\begin{lemma}\label{Lemma 4.2} If $\t/2<\t_\star<\t$, then there exist two
positive constants $\wt\sigma=\wt\sigma(v,u,\t; v_ \star, u_ \star,\t_ \star)$ and
$\eta_0=\eta_0(v,u,\t; v_\star,u_\star,\t_\star)$ such that if $|v-v_\star|+|u-u_\star|+|\t-\t_\star|<\eta_0$, we have for
$g(\xi)\in  \mathfrak{N}^\bot$,
$$
-\int\f{g\mb{L}_\mb{M}g}{\mb{M}_\star}d\xi\geq
\wt\sigma\int\f{\nu(|\xi|)g^2}{\mb{M}_\star}d\xi,\qquad -\int\f{g\mb{N}_\mb{M}g}{\mb{M}_\star}d\xi\geq
\wt\sigma\int\f{\nu(|\xi|)g^2}{\mb{M}_\star}d\xi.
$$
\end{lemma}

 \begin{lemma}\label{Lemma 4.3} Under the assumptions in Lemma \ref{Lemma 4.2}, we
have  for each $g( \xi)\in  \mathfrak{N}^\bot$,
$$
 \int\f{\nu(|\xi|)}{\mb{M}_\star}|\mb{L}_\mb{M}^{-1}g|^2d\xi
\leq \wt\sigma^{-2}\int\f{\nu(|\xi|)^{-1}g^2}{\mb{M}_\star}d\xi,~~{\rm
and}~~
\int\f{\nu(|\xi|)}{\mb{M}_\star}|\mb{N}_\mb{M}^{-1}g|^2d\xi\leq
\wt\sigma^{-2}\int\f{\nu(|\xi|)^{-1}g^2}{\mb{M}_\star}d\xi.
$$
\end{lemma}

\begin{lemma}\label{hlm}
For $0<t\leq+\infty$, suppose that $h(x,t)$ satisfies
$$
h\in L^{\infty}(0,t;L^2(\mb{R})),\quad h_x\in L^{2}(0,t;L^2(\mb{R})),\quad h_t\in L^{2}(0,t;H^{-1}(\mb{R})).
$$
Then
\begin{equation}
\int_0^t\int h^2\hat{w}^2dxds\leq 4\pi\|h(0)\|^2+4\pi\alpha^{-1}\int_0^t\|h_x\|^2ds
+8\alpha\int_0^t<h_t,hg^2>_{H^{-1}\times H^1}ds
\end{equation}
for $\alpha>0$, and
$$
\hat{w}(x,t)=(1+t)^{-\frac{1}{2}}\exp\left(-\frac{\alpha x^2}{1+t}\right),\quad
g(x,t)=\int_{-\infty}^x \hat{w}(y,t)dy.
$$
\end{lemma}

\begin{remark}
In Lemmas \ref{Lemma 4.2}-\ref{Lemma 4.3}, $\eta_0$ may not be sufficiently small positive constant. However, in the proof of Theorem \ref{thm} in the following sections, the smallness of $\eta_0$ is crucially used to close the a priori assumption \eqref{assumption-ap}. From the proof of Lemmas \ref{Lemma 4.2}-\ref{Lemma 4.3} in \cite{Liu-Yang-Yu-Zhao}, it can be seen that Lemmas \ref{Lemma 4.2}-\ref{Lemma 4.3} still hold true for more smaller positive constant $\eta_0.$
The proof of Lemma \ref{hlm} can be found in Huang-Li-Matsumura \cite{Huang-Li-matsumura} for the stability of viscous contact wave for the compressible Navier-Stokes equations.
\end{remark}

\section{Lower Order Estimates}
\setcounter{equation}{0}
The proof of Theorem \ref{thm} is shown by the continuum argument for the local solution to the system \eqref{VPB} or equivalently the system \eqref{F1-f}-\eqref{F2-pc-L}, while the local-in-time solution can be constructed similarly as in \cite{Guo1, Yang-Zhao}. Therefore, to prove Theorem \ref{thm}, it is sufficient to close the a priori assumption \eqref{assumption-ap} and verify the time-asymptotic behaviors of the solution. For the a priori estimates, we start from the lower order estimates.
\begin{proposition} \label{Prop3.1}
For each $(v_\pm,u_\pm,\t_\pm)$, there exists a positive constant $C$ such that, if
$\d\leq \d_0$ for a suitably small positive constant $\d_0$, then it holds that
\begin{equation}\label{LE-A}
\begin{array}{ll}
\di \sup_{0\leq t\leq +\i}\|(\phi,\psi,\zeta,\phi_x,\Phi_x ,n_2)(\cdot,t)\|^2+\int\int\f{|(\widetilde{\mb{G}},\mb{P}_c F_2)|^2}{\mb{M}_\star }(x,\xi,t)d\xi dx
+\d\int_0^t\||\hat w|(\phi,\psi,\zeta)\|^2ds\\
\di +\int_0^t\Big[\sum_{|\alpha'|=1}\|\pa^{\alpha'}(\phi,\psi,\zeta,n_2)\|^2+\|(\Phi_x,\Phi_{xt},n_2)\|^2\Big]ds
+\int_0^t\int\int\f{\nu(|\xi|)|(\widetilde{\mb{G}},\mb{P}_c F_2)|^2}{\mb{M}_\star }d\xi dxds\\
\di \leq C\|(\phi_0,\psi_0,\zeta_0,\phi_{0x},\Phi_{0x},n_{20})\|^2+C\int\int\f{|(\widetilde{\mb{G}},\mb{P}_c F_2)|^2}{\mb{M}_\star }(x,\xi,0) d\xi dx+C\d\\
\di +C\sum_{1\leq |\alpha|\leq 2}\int_0^t\int\int \f{\nu(|\xi|)|\pa^{\alpha}(\mb{G},\mb{P}_c F_2)|^2}{\mb{M}_\star } d\xi dxds
+C(\chi+\d)\int_0^t\int\int \f{\nu(|\xi|)|\pa_{\xi_1}(\mb{P}_c F_2)_x|^2}{\mb{M}_\star } d\xi dxds\\
\di+C(\chi+\d)\sum_{|\beta'|=1}\int_0^t\int\int \f{\nu(|\xi|)|\pa^{\beta'}(\widetilde{\mb{G}},\mb{P}_c F_2)|^2}{\mb{M}_\star } d\xi dxds.
\end{array}
\end{equation}

\end{proposition}

\textbf{Proof:} The proof of the lower order estimates in Proposition \ref{Prop3.1} includes the following steps.

\

\underline{Step 1. Estimation on $\|(\phi,\psi,\zeta)(\cdot, t)\|^2$.}

\

First, multiplying $\eqref{sys-h}_2$ by $\psi_1$ leads to
\begin{equation}\label{6.1}
\begin{array}{ll}
\di \left(\frac{\psi_1^2}{2}\right)_t+\left((p-\bar p)\psi_1-\f43\left(\f{\mu(\t)}{v}u_{1x}-\f{\mu(\bar\t)}{\bar v}\bar u_{1x}\right)\psi_1
+\psi_1\int \xi_1^2\Pi_1d\xi\right)_x+\f43\f{\mu(\t)}{v}\psi_{1x}^2-\f23\frac{\zeta}{v}\psi_{1x}\\
\di\quad-\f23\bar\t\left(\frac{1}{v}-\frac{1}{\bar v}\right)\phi_t
+\f43\left(\f{\mu(\t)}{v}-\f{\mu(\bar\t)}{\bar v}\right)\bar u_{1x}\psi_{1x}
-\Phi_x n_2\psi_1-\psi_{1x}\int \xi_1^2\Pi_1d\xi=-\bar R_1\psi_1,
\end{array}
\end{equation}
multiplying $\eqref{sys-h}_3$ by $\psi_i~ (i=2,3)$,   gives
\begin{equation}\label{6.2}
\begin{array}{ll}
\di \left(\frac{\psi_i^2}{2}\right)_t+\left(-\Big(\f{\mu(\t)}{v}u_{ix}\Big)\psi_i+\psi_i\int\xi_1\xi_i\Pi_1d\xi\right)_x
+\f{\mu(\t)}{v}\psi_{ix}^2-\psi_{ix}\int\xi_1\xi_i\Pi_1d\xi=0,
\end{array}
\end{equation}
multiplying $\eqref{sys-h}_4$ by $\z\t^{-1}$, we have,
\begin{equation}\label{6.3}
\begin{array}{ll}
\di \f{\z\z_t}{\t}+\left(-\Big(\f{\k(\t)}{v}\t_x-\f{\k(\bar\t)}{\bar v}\bar\t_x\Big)\f{\z}{\t}+\f{\z}{\t}\int\f12\xi_1|\xi|^2\Pi_1d\xi
-\sum_{i=1}^3u_i\f{\z}{\t}\int\xi_1\xi_i\Pi_1d\xi\right)_x+\f{2\z}{3v}\psi_{1x}\\
\di +(p-\bar p)\bar u_{1x}\f{\z}{\t}+\f{\k(\t)\bar\t}{\t^2v}\z_x^2-\f{\k(\t)}{\t^2v}\z\z_x\bar\t_x
+\Big(\f{\k(\t)}{v}-\f{\k(\bar\t)}{\bar v}\Big)\f{\bar\t}{\t^2}\bar\t_x\z_x-\Big(\f{\k(\t)}{v}-\f{\k(\bar\t)}{\bar v}\Big)\f{\z\bar\t_x^2}{\t^2}\\
\di-\f43\f{\mu(\t)}{\t v}\z\psi_{1x}^2-\sum_{i=2}^3\f{\mu(\t)}{\t v}\z\psi_{ix}^2-\f{8}{3}\f{\mu(\t)}{\t v}\z\psi_{1x}\bar u_{1x}
-\f43\Big(\f{\mu(\t)}{v}-\f{\mu(\bar \t)}{\bar v}\Big)\f{\z}{\t}\bar u_{1x}^2\\
\di-\Big(\f{\z}{\t}\Big)_x\int\f12\xi_1|\xi|^2\Pi_1d\xi+\sum_{i=1}^3\Big(u_i\f{\z}{\t}\Big)_x\int\xi_1\xi_i\Pi_1d\xi
-\Phi_x\f{\z}{\t}\int \xi_1\mb{P}_cF_2d\xi=-\bar R_2\f{\z}{\t}.
\end{array}
\end{equation}
Note that
\begin{equation}\label{6.4}
\di-\f23\bar\t\left(\frac{1}{v}-\frac{1}{\bar v}\right)\phi_t=\left(\f23\bar\t\Phi\left(\frac{v}{\bar v}\right)\right)_t
-\f23\bar\t_t\Phi\left(\frac{v}{\bar v}\right)+\frac{\bar p\phi^2}{v\bar v}\bar v_t,
\end{equation}
\begin{equation}
\di \f{\z\z_t}{\t}=\Big(\bar\t\Phi\Big(\f{\t}{\bar\t}\Big)\Big)_t+\bar\t_t\Phi\Big(\f{\bar\t}{\t}\Big),
\end{equation}
\begin{equation}
\begin{array}{ll}
\di-\f23\bar\t_t\di=\f23P^r_-(U_{1-}^r)_x+\f23P^r_+(U_{1+}^r)_x-p^* U_{1x}^{cd}\\[2mm]
\di\quad\qquad=\f23\bar p(U_{1-}^r)_x+\f23\bar p(U_{1+}^r)_x+\f23(P^r_- -\bar p)(U_{1-}^r)_x
+\f23(P^r_+-\bar p)(U_{1+}^r)_x-p^* U_{1x}^{cd}
\end{array}
\end{equation}
and
\begin{equation}
-\f23\bar\t_t\Phi\left(\frac{v}{\bar v}\right)+\frac{\bar p\phi^2}{v\bar v}\bar v_t+\bar\t_t\Phi\left(\frac{\bar\t}{\t}\right)
+(p-\bar p)\bar u_{1x}\frac{\z}{\t}
=\bar p\left(\Phi\left(\frac{\t \bar v}{\bar\t v}\right)+\f53\Phi\left(\frac{v}{\bar v}\right)\right)\big((U_{1-}^r)_x+(U_{1+}^r)_x\big)+\widetilde Q,
\end{equation}
where
\begin{equation}
\begin{array}{ll}\label{6.5}
\di \widetilde Q=U_{1x}^{cd}\left(\frac{\bar p\phi^2}{v\bar v}-p^m\Phi\left(\frac{v}{\bar v}\right)
+\frac{3}{2}p^m\Phi\left(\frac{\bar\t}{\t}\right)+\frac{\z}{\t}(p-\bar p)\right)\\[3mm]
\di\qquad+\f23(p_--\bar p)(U_{1-}^r)_x\left(\Phi\left(\frac{v}{\bar v}\right)-\frac{3}{2}\Phi\left(\frac{\bar\t}{\t}\right)\right)
+\f23(p_+-\bar p)(U_{1+}^r)_x\left(\Phi\left(\frac{v}{\bar v}\right)-\frac{3}{2}\Phi\left(\frac{\bar\t}{\t}\right)\right).
\end{array}
\end{equation}
It follows from \eqref{6.1}-\eqref{6.5} that,
\begin{equation}\label{E1}
\begin{array}{l}
\di\left(\f23\bar\t\Phi\left(\frac{v}{\bar v}\right)+\bar\t\Phi\left(\frac{\theta}{\bar\t}\right)
+\f12\sum_{i=1}^3\psi_i^2\right)_t
+\f43\f{\mu(\t)\bar\t}{\t v}\psi_{1x}^2+\sum_{i=2}^3\f{\mu(\t)\bar\t}{\t v}\psi_{ix}^2
+\f{\k(\t)\bar\t}{\t^2 v}\zeta_{x}^2\\
\di+\bar p\left(\Phi\left(\frac{\t \bar v}{\bar\t v}\right)+\f53\Phi\left(\frac{v}{\bar v}\right)\right)\big((U^r_{1-})_x+(U^r_{1+})_x\big)+ Q_1+Q_2-\Phi_x n_2\psi_1-\Phi_x\f{\zeta}{\t}\int\xi_1\mb{P}_cF_2d\xi\\
\di +(\cdots)_x
=-\psi_1\bar R_1-\frac{\zeta}{\theta}\bar R_2,
\end{array}
\end{equation}
where $(\cdots)_x$ means some conservative terms which vanish after integration on $\mathbf{R}$ and
$$
\Psi(z)=z-\ln z-1, \quad z>0,
$$
\begin{equation}\label{q1}
\begin{array}{ll}
\di  Q_1=\widetilde Q-\f{\k(\t)}{\t^2v}\zeta\zeta_x\bar\t_x+\left(\f{\k(\t)}{v}-\f{\k(\bar\t)}{\bar v}\right)\f{\bar\t}{\t^2}\bar\t_x\zeta_x
-\left(\f{\k(\t)}{v}-\f{\k(\bar\t)}{\bar v}\right)\f{\zeta\bar\t_x^2}{\t^2}\\
\di\quad\quad+\f43\left(\f{\mu(\t)}{v}-\f{\mu(\bar\t)}{\bar v}\right)\bar u_{1x}\psi_{1x} -\f83\frac{\mu(\t)}{\t v}\zeta\psi_{1x}\bar u_{1x}-\f43\left(\f{\mu(\t)}{v}-\f{\mu(\bar\t)}{\bar v}\right)\f{\zeta}{\t}\bar u_{1x}^2
\end{array}
\end{equation}
and
\begin{equation}\label{q2}
\begin{array}{ll}
\di Q_2=-\left(\f{\bar\t\zeta_x}{\t^2}-\f{\zeta\bar\t_x}{\t^2}\right)\left(\int\f12\xi_1|\xi|^2\Pi_1d\xi
-\sum_{i=1}^3u_i\int\xi_1\xi_i\Pi_1d\xi\right)\\[2mm]
\di\quad\quad\quad-\sum_{i=1}^3\psi_{ix}\int\xi_1\xi_i\Pi_1d\xi+\f{\zeta}{\t}\sum_{i=1}^3(\psi_{ix}+\bar u_{ix})\int\xi_1\xi_i\Pi_1d\xi.
\end{array}
\end{equation}
Recalling (iii) in Lemma \ref{rare-pro}, we can compute
\begin{equation}
\begin{array}{ll}
\di|(P^r_--\bar p)(U_{1-}^r)_x|\\
\di\leq C\big(|\Theta^{cd}-\t_-^*|+|\Theta^{r}_+-\t_+^*|+|V^{cd}-v_-^*|+|V^{r}_+-v_+^*|\big)|(U_{1-}^r)_x|\\
\di\leq C\big(|\Theta^{cd}-\t_-^*|+|\Theta^{r}_+-\t_+^*|+|V^{cd}-v_-^*|+|V^{r}_+-v_+^*|\big)\big|_{\Omega_-}
+C|(U_{1-}^r)_x|\big|_{\Omega_c\cap\Omega_+}\\
\di\leq C\d e^{-c_0(|x|+t)},
\end{array}
\end{equation}
which leads to
\begin{equation}
\di|\widetilde Q|\leq C|U_{1x}^{cd}|(\phi^2+\zeta^2)+C\d e^{-c_0(|x|+t)}(\phi^2+\zeta^2)
\end{equation}
and
\begin{equation}
\di| Q_1|\leq |\widetilde Q|+\frac{\mu(\t)\bar \t}{3\t v}\psi_{1x}^2+\frac{\k(\t)\bar\t}{3\t^2 v}\z_x^2
+C(\phi^2+\zeta^2)(|\bar u_{1x}|+|\bar \t_x|^2).
\end{equation}
Note that
\begin{equation}
\begin{array}{ll}
(\phi^2+\zeta^2)(|\bar u_{1x}|+|\bar \t_x|^2)\leq C(\phi^2+\zeta^2)((\Theta_x^{cd})^2+|(U^r_{1-})_x|^2+|(U^r_{1+})_x|^2)\\
\quad \leq C\d(1+t)^{-1}(\phi^2+\zeta^2)e^{-\frac{c_1x^2}{1+t}}+C\d \bar p\left(\Phi\left(\frac{\t \bar v}{\bar\t v}\right)+\f53\Phi\left(\frac{v}{\bar v}\right)\right)((U^r_{1-})_x+(U^r_{1+})_x).
\end{array}
\end{equation}
Following the same calculations as in \cite{Huang-Li-matsumura}, it holds that
\begin{equation}
\di \|(\bar R_1,\bar R_2)\|_{L^1}\leq C\d^{\f18}(1+t)^{-\f78}.
\end{equation}
Then we have
\begin{equation}\label{int-F}
\begin{array}{ll}
\di\int_0^t\int\left(\bar R_1\psi_1+\bar R_2\frac{\z}{\t}\right)dxds\leq \int_0^t\|(\bar R_1,\bar R_2)\|_{L^1}\|(\psi_1,\zeta)\|_{L^{\infty}}ds\\
\di\quad\leq \d^{\f18}\int_0^t(1+s)^{-\f78}\|(\psi_1,\zeta)\|^{\f12}\|(\psi_{1x},\zeta_x)\|^{\f12}ds\\
\di\quad\leq \int_0^t\int\left(\frac{\mu(\t)\bar \t}{3\t v}\psi_{1x}^2+\frac{\k(\t)\bar\t}{3\t^2 v}\z_x^2\right)dxds
+C\d^{\f16}\int_0^t(1+s)^{-\f76}\Big(1+\Big\|\Big(\psi_1,\sqrt{\Phi\Big(\frac{\t}{\bar\t}\Big)}\Big)\Big\|^2\Big)ds.
\end{array}
\end{equation}
By Cauchy's inequality and \eqref{4.18}, we have
\begin{equation}
\begin{array}{ll}
\di \left|\int_0^t\int \Phi_x\Big(n_2\psi_1+\f{\zeta}{\t}\int \xi_1\mb{P}_cF_2\Big) dxds\right|\\
\di\leq C\d\int_0^t\|\Phi_x\|\Big(\|\psi_1\|_{L^{\infty}}\|n_2\|+\|\zeta\|_{L^{\i}}
\Big(\int\int\frac{\nu(|\xi|)}{\mb{M}_\star }|\mb{P}_cF_2|^2d\xi dx\Big)^{\f12}\Big)ds\\
\di \leq C (\chi+\d)\int_0^t\|(\Phi_x,n_2)\|^2ds+C(\chi+\d)\int_0^t\int\int\frac{\nu(|\xi|)}{\mb{M}_\star }|\mb{P}_cF_2|^2d\xi dxds
\end{array}
\end{equation}
and
\begin{equation}
\begin{array}{ll}
\di \left|\int_0^t\int \left(\f{\bar\t\zeta_x}{\t^2}-\f{\zeta\bar \t_x}{\t^2}\right)\left(\int\f12\xi_1|\xi|^2\Pi_1d\xi
-\sum_{i=1}^3u_i\int\xi_1\xi_i\Pi_1d\xi\right) dxds\right|\\
\di \leq \eta \int_0^t\|(\zeta_x,\zeta \bar\t_{x})\|^2ds+C_\eta \int_0^t\int\Big[|\int
\f12\xi_1|\xi|^2\Pi_1 d\xi|^2+\sum_{i=1}^3|\int \xi_1\xi_i\Pi_1 d\xi|^2\Big]dxds,
\end{array}
\end{equation}
with some small positive constant $\eta>0$ to be determined and the positive constant $C_\eta$ depending on $\eta$.
Note that by \eqref{Pi-1}, it holds that
\begin{equation}\label{I3+}
\begin{array}{ll}
\di |\f12\int \xi_1|\xi|^2\Pi_1 d\xi|=|\int \f12 \xi_1|\xi|^2
\mb{L}_\mb{M}^{-1}[\mb{G}_t-\f{u_1}{v}\mb{G}_x+\f1v\mb {P}_1(\xi_1\mb{G}_x)
+\f1v\mb{P}_1(\Phi_x\partial_{\xi_1}F_2)-2Q(\mb{G}, \mb{G})] d\xi| :=\sum_{i=1}^5I_{i}.
\end{array}
\end{equation}
Choose the global Maxellian $\mb{M}_\star =\mb{M}_{[v_\star,u_\star,\t_\star]}$ such that
\begin{equation}\label{ma1}
v_\star>0,\qquad \f12\t(t,x)<\t_\star<\t(t,x),
\end{equation}
and
\begin{equation}\label{ma2}
|v(x,t)-v_\star|+|u(x,t)-u_\star|+|\t(x,t)-\t_\star|<\eta_0
\end{equation}
with $\eta_0$ being the small positive constant in Lemma \ref{Lemma 4.2}. In fact, if the total wave strength $\delta=|(v_+-v_-,u_+-u_-,\t_+-\t_-)|$ is suitably small, then it is holds that
$$
\f12\sup_{(x,t)}\bar\t(x,t)=\f12\t_+<\inf_{(x,t)}\bar v(x,t)=\t_-.
$$
Therefore, we can choose the global Maxwellian $\mb{M}_\star =\mb{M}_{[v_\star,u_\star,\t_\star]}$ satisfying \eqref{ma1} and \eqref{ma2}
provided that the solution $(v,u,\t)(x,t)$ is near the ansatz $(\bar v,\bar u,\bar\t)(x,t)$ as in a priori assumption \eqref{assumption-ap}.
Then with such chosen $\mb{M}_\star $, it holds that
\begin{equation}\label{I31}
\begin{array}{ll}
\di I_{1}+I_2 \leq \Big(\int\f{\nu(|\xi|)|\mb{L}_\mb{M}^{-1}[\mb{G}_t-\f{u_1}{v}\mb{G}_x]|^2}{\mb{M}_\star }
d\xi\Big)^{\f12} \Big(\int\mb{M}_\star \nu(|\xi|)^{-1}(\f12\xi_1|\xi|^2)^2
d\xi\Big)^{\f12}\\
 \di\qquad \leq C \Big(\int\f{\nu(|\xi|)|(\mb{G}_x,\mb{G}_t)|^2}{\mb{M}_\star }
d\xi\Big)^{\f12},
\end{array}
\end{equation}
and
\begin{equation}\label{I32}
\begin{array}{ll}
\di I_{3}=|\int \f12 \xi_1|\xi|^2 \mb{L}_\mb{M}^{-1}[\f1v\mb
{P}_1(\xi_1\mb{G}_x)] d\xi|\leq C
\Big(\int\f{\nu(|\xi|)|\mb{L}_\mb{M}^{-1}[\f1v\mb
{P}_1(\xi_1\mb{G}_x)]|^2}{\mb{M}_{[v_\star,u_\star,2\t_\star]}}
d\xi\Big)^{\f12} \\
 \di\quad ~~~ \leq C \Big(\int\f{\nu(|\xi|)^{-1}|\mb
{P}_1(\xi_1\mb{G}_x)|^2}{\mb{M}_{[v_\star,u_\star,2\t_\star]}} d\xi\Big)^{\f12}\leq
C \Big(\int\f{\nu(|\xi|)|\mb{G}_x|^2}{\mb{M}_\star } d\xi\Big)^{\f12}.
\end{array}
\end{equation}
Furthermore, one has
\begin{equation}\label{I33}
\begin{array}{ll}
\di I_{4}=|\int \f12 \xi_1|\xi|^2 \mb{L}_\mb{M}^{-1}[\f1v\mb
{P}_1(\Phi_x\partial_{\xi_1}F_2)] d\xi|
  \leq C|\Phi_x| \Big(\int\f{\nu(|\xi|)^{-1}|n_2v\mb
{M}_{\xi_1}+\partial_{\xi_1}(\mb{P}_c F_2)|^2}{\mb{M}_\star }
d\xi\Big)^{\f12}\\
\di\quad \leq C|\Phi_x||n_2|+C|\Phi_x|
\Big(\int\f{\nu(|\xi|)|\partial_{\xi_1}(\mb{P}_c F_2)|^2}{\mb{M}_\star }
d\xi\Big)^{\f12},
\end{array}
\end{equation}
and
\begin{equation}\label{I34}
\begin{array}{ll}
\di I_{5}=|\int \f12 \xi_1|\xi|^2 \mb{L}_\mb{M}^{-1}[2Q(\mb{G},
\mb{G})] d\xi|\leq C \Big(\int\f{\nu(|\xi|)|\mb{L}_\mb{M}^{-1}[Q(\mb{G},
\mb{G})]|^2}{\mb{M}_\star }
d\xi\Big)^{\f12} \\
 \di\quad ~~ \leq C\Big(\int\f{\nu(|\xi|)^{-1}|Q(\mb{G},
\mb{G})|^2}{\mb{M}_\star } d\xi\Big)^{\f12}\leq
C\Big(\int\f{\nu(|\xi|)|\mb{G}|^2}{\mb{M}_\star }
d\xi\Big)^{\f12}\Big(\int\f{|\mb{G}|^2}{\mb{M}_\star }
d\xi\Big)^{\f12} \\
\di \leq C\Big(\int\f{\nu(|\xi|)|\widetilde{\mb{G}}|^2}{\mb{M}_\star }
d\xi\Big)^{\f12}\Big(\int\f{|\widetilde{\mb{G}}|^2}{\mb{M}_\star }
d\xi\Big)^{\f12}+C|(\bar\t_x,\bar
u_{1x})|\Big(\int\f{\nu(|\xi|)|\widetilde{\mb{G}}|^2}{\mb{M}_\star }
d\xi\Big)^{\f12}+C|(\bar\t_x,\bar u_{1x})|^2.
\end{array}
\end{equation}
Substituting \eqref{I31}-\eqref{I34} into \eqref{I3+}, one can arrive at
\begin{equation}\label{I3-e}
\begin{array}{ll}
\di \left|\int_0^t\int \Big(\f{\bar\t\zeta_x}{\t^2}-\f{\zeta\t_x}{\t^2}\Big)\Big(\int\f12\xi_1|\xi|^2\Pi_1d\xi
-\sum_{i=1}^3u_i\int\xi_1\xi_i\Pi_1d\xi\Big) dxds\right|\\
\di \leq \eta \int_0^t\|(\zeta_x,\zeta \bar\t_{x})\|^2ds+C_\eta \int_0^t\int\int\f{\nu(|\xi|)|(\mb{G}_x,\mb{G}_t)|^2}{\mb{M}_\star }
d\xi dxds+ C_\eta(\chi+\d)\int_0^t\|\Phi_x\|^2ds\\
\di \quad+C_\eta(\chi+\d)\int_0^t\int\int\f{\nu(|\xi|)|(\widetilde{\mb{G}},\pa_{\xi_1}(\mb{P}_c F_2))|^2}{\mb{M}_\star }d\xi dxds+C\d.
\end{array}
\end{equation}
Similarly, one has
\begin{equation}\label{I5}
\begin{array}{ll}
\di \left|\int_0^t\int \Big(-\sum_{i=1}^3\psi_{ix}\int\xi_1\xi_i\Pi_1d\xi+\f{\zeta}{\t}\sum_{i=1}^3(\psi_{ix}+\bar u_{ix})\int\xi_1\xi_i\Pi_1d\xi\Big) dxds\right|\\
\di \leq \eta \int_0^t\|(\psi_x,\zeta \bar u_{1x})\|^2ds+C_\eta \int_0^t\int\int\f{\nu(|\xi|)|(\mb{G}_x,\mb{G}_t)|^2}{\mb{M}_\star }
d\xi dxds+ C_\eta(\chi+\d)\int_0^t\|\Phi_x\|^2ds\\[3mm]
\di \quad+C_\eta(\chi+\d)\int_0^t\int\int\f{\nu(|\xi|)|(\widetilde{\mb{G}},\pa_{\xi_1}(\mb{P}_c F_2))|^2}{\mb{M}_\star }d\xi dxds+C\d.
\end{array}
\end{equation}
Integrating the equation \eqref{E1} over $x$ and $t$, choosing $\eta$ suitable small and applying Gronwall's inequality,
which together with the above estimates give that
\begin{equation}\label{le}
\begin{array}{ll}
\di\|(\phi,\psi,\zeta)(\cdot,t)\|^2+\int_0^t\|(\psi_x,\zeta_x)\|^2ds+\int_0^t\int \big((U^r_{1-})_x+(U^r_{1+})_x\big)(\phi^2+\zeta^2) dxds\\
\di\leq C\|(\phi,\psi,\zeta)(\cdot,0)\|^2+C\d
+C\d\int_0^t(1+s)^{-1}\int(\phi^2+\zeta^2)e^{-\frac{c_1x^2}{1+s}}dxds+ C(\chi+\d) \int_0^t\|(\Phi_x, n_2)\|^2ds\\
\di+C(\chi+\delta)\int_0^t\int\int\f{\nu(|\xi|)|(\widetilde{\mb{G}},\mb{P}_c F_2,\pa_{\xi_1}(\mb{P}_c F_2))|^2}{\mb{M}_\star }
d\xi dxds+C \int_0^t\int\int\f{\nu(|\xi|)|(\mb{G}_x,\mb{G}_t)|^2}{\mb{M}_\star }d\xi dxds.
\end{array}
\end{equation}

\

\underline{Step 2. Estimation on $\|\phi_x(\cdot, t)\|^2$.}

\

First, we state the following lemma which is helpful to get the lower order estimate, and the proof of the Lemma can be  done similarly as in \cite{Huang-Li-matsumura},
the only difference here is that we need to additionally care about the microscopic terms and the electric terms.
\begin{lemma}\label{le-1}
For $\alpha\in(0,\frac{c_1}{4}]$ and $\hat w$ defined in Lemma \ref{hlm}, there exists some positive
constant $C$ depending on $\alpha$, such that the following estimate holds
\begin{equation}\label{important}
\begin{array}{l}
\di\int_0^t\int(\phi^2+\psi^2+\zeta^2)\hat w^2 dxds\leq C+C\|(\phi,\psi,\zeta)(\cdot,t)\|^2+C\int_0^t\|(\phi_x,\psi_x,\z_x)\|^2ds\\
\di\quad+C\int_0^t\int \big((U_{1-}^r)_x+(U_{1+}^r)_x\big)(\phi^2+\zeta^2)dxds+C\int_0^t\int\int\f{\nu(|\xi|)|(\mb{G}_x,\mb{G}_t)|^2}{\mb{M}_\star }d\xi dxds\\
\di\quad+C(\chi+\d)\int_0^t\|(\Phi_x,n_2)\|^2ds+C(\chi+\delta)\int_0^t\int\int\f{\nu(|\xi|)|(\widetilde{\mb{G}},\mb{P}_c F_2,\pa_{\xi_1}(\mb{P}_c F_2))|^2}{\mb{M}_\star }d\xi dxds.
\end{array}
\end{equation}
\end{lemma}
Next, we rewrite \eqref{sys-h}$_2$ as
$$
\begin{array}{l}
\displaystyle
\di\f43\mu(\t)\Big(\f{\phi_x}{v}\Big)_t+\f{2\t}{3v}\f{\phi_x}{v}=\psi_{1t}+\f{2\z_x}{3v}+\frac23\bar\t_x\Big(\frac1v-\frac1{\bar v}\Big)
-\frac23\bar v_x\Big(\frac{\t}{v^2}-\frac{\bar\t}{\bar v^2}\Big)-\Phi_xn_2\\[2mm]
\di\quad-\f43\mu(\t)\Big(\f{\bar v_x}{v}\Big)_t-\f43\mu(\t)_x\f{u_{1x}}{v}+(\bar p-P^r_+-P^r_-)_x+U^{cd}_{1t}+\int\xi_1^2\Pi_{1x} d\xi.
\end{array}
$$
Multiplying the above equation by $\frac{\phi_x}{v}$ and noting that
$$
\begin{array}{ll}
\di
\psi_{1t}\f{\phi_x}{v}=\Big(\psi_1\f{\phi_x}{v}\Big)_t-\psi_1\Big(\f{\phi_x}{v}\Big)_t=\Big(\psi_1\f{\phi_x}{v}\Big)_t
-\f{\psi_1\psi_{1xx}}{v}+\f{\psi_1\phi_x\psi_{1x}+\psi_1\phi_x\bar u_{1x}}{v^2}\\
\di \qquad\quad
=\Big(\psi_1\f{\phi_x}{v}\Big)_t-\Big(\f{\psi_1\psi_{1x}}{v}\Big)_x+\f{\psi_{1x}^2}{v}-\f{\psi_1\psi_{1x}\bar v_x}{v^2}
+\f{\psi_1\phi_x\bar u_{1x}}{v^2},
\end{array}
$$
we can obtain
\begin{equation}\label{ff}
\begin{array}{l}
\displaystyle
\Big(\f23\mu(\t)\f{\phi_x^2}{v^2}-\psi_1\f{\phi_x}{v}\Big)_t+\f{2\t}{3v}\f{\phi_x^2}{v^2}=
-\Big(\f{\psi_1\psi_{1x}}{v}\Big)_x+\f{\psi_{1x}^2}{v}-\f{\psi_1\psi_{1x}\bar v_x}{v^2}
+\f{\psi_1\phi_x\bar u_{1x}}{v^2}\\[2mm]
\di\quad +\f{2\z_x\phi_x}{3v^2}+\Big(\frac23\bar\t_x\Big(\frac1v-\frac1{\bar v}\Big)
-\frac23\bar v_x\Big(\frac{\t}{v^2}-\frac{\bar\t}{\bar v^2}\Big)\Big)\frac{\phi_x}{v}
-\Phi_x n_2\f{\phi_x}{v}-\f43\mu(\t)\Big(\f{\bar v_x}{v}\Big)_t\f{\phi_x}{v}\\[2mm]
\di\quad-\f43\mu(\t)_x\f{u_{1x}\phi_x}{v^2}+\f23\mu(\t)_t\f{\phi_x^2}{v^2}
+(\bar p-P^r_+-P^r_-)_x\frac{\phi_x}{v}+\f{U^{cd}_{1t}\phi_x}{v}+\f{\phi_x}{v}\int\xi_1^2\Pi_{1x} d\xi.
\end{array}
\end{equation}
One can carry out the similar steps as in obtaining \eqref{I3-e} to get
\begin{equation}\label{pix}
\begin{array}{ll}
\di \int_0^t\int|\int \xi_1^2\Pi_{1x} d \xi|^2dxds\leq
C\sum_{|\alpha|=2}\int_0^t\int\int\frac{\nu(|
\xi|)}{\mb{M}_\star }|\partial^\alpha \mathbf{G}|^2d \xi
dxds+C\d\\
\di+C(\chi+\delta)\int_0^t\Big[\|(\Phi_x,n_2,\phi_x,\psi_x,\zeta_x)\|^2+\sum_{|\alpha^\prime|=
1}\int\int\frac{\nu(|
\xi|)}{\mb{M}_\star }|(\widetilde{\mathbf{G}},\partial^{\alpha^\prime}\mb{G})|^2
d \xi dx\Big]ds\\
\di+C(\chi+\delta)\int_0^t\int\int\frac{\nu(|\xi|)}{\mb{M}_\star }|(\pa_{\xi_1}(\mb{P}_c F_2),\pa_{\xi_1}(\mb{P}_c F_2)_x)|^2d\xi dxds.
\end{array}
\end{equation}
On the other hand, By Lemma \ref{rare-pro}, we have
$$
|(\bar p-P^r_+-P^r_-)_x|\leq C\d e^{-c_0(|x|+t)},
$$
therefore, the other terms on the right-hand side of \eqref{ff} can be controlled by
\begin{equation}\label{uw}
\begin{array}{ll}
\di\Big|\int_0^t\int\Big(-\Big(\f{\psi_1\psi_{1x}}{v}\Big)_x+\cdots+\f{\bar u_{1t}\phi_x}{v}\Big) dxds\Big|
\leq C\d+ C(\chi+\d+\eta)\int_0^t\|\phi_x\|^2ds\\
\di+C_\eta\int_0^t\|(\psi_x,\z_x)\|^2ds+C_\eta\int_0^t\||(\bar v_x,\bar u_{1x}, \bar\t_x)|(\phi,\psi,\z)\|^2ds+C\chi\int_0^t\|\Phi_x\|^2ds.
\end{array}
\end{equation}
Integrating the equation \eqref{ff} with respect to $x,t$, then using Cauchy's inequality and \eqref{pix}-\eqref{uw}
and choosing $\chi,\delta,\eta$ suitably small, one can obtain
\begin{equation}\label{phi-x-E}
\begin{array}{l}
\displaystyle
\|\phi_x(\cdot,t)\|^2+\int_0^t\|\phi_x\|^2ds \leq
C\Big[\|\psi_1(\cdot,t)\|^2+\|(\phi_{0x},\psi_{10})\|^2+\delta\Big]+C\int_0^t\|(\psi_{1x},\zeta_x)\|^2ds
\\
\quad\displaystyle+C\int_0^t\||(\bar v_x,\bar u_{1x}, \bar\t_x)|(\phi,\psi_1,\zeta)\|^2ds
+C\sum_{|\alpha|=2}\int_0^t\int\int\frac{\nu(| \xi|)}{\mb{M}_\star }|\partial^\alpha\mathbf{G}|^2d \xi dxds\\
\di \quad+C(\chi+\d)\int_0^t\|(\Phi_x,n_2)\|^2ds
+C(\chi+\delta)\sum_{|\alpha^\prime|=
1}\int_0^t\int\int\frac{\nu(| \xi|)}{\mb{M}_\star }|(\widetilde{\mathbf{G}},\partial^{\alpha^\prime}\mb{G})|^2 d \xi dxds\\
\di\quad+C(\chi+\delta)\int_0^t\int\int\frac{\nu(|\xi|)}{\mb{M}_\star }|(\pa_{\xi_1}(\mb{P}_c F_2),\pa_{\xi_1}(\mb{P}_c F_2)_x)|^2d\xi dxds.
\end{array}
\end{equation}
Then we estimate $\|(\phi,\psi,\zeta)_t\|^2$. For this, we use the
system \eqref{sys-h-o}. Multiplying the equation $\eqref{sys-h-o}_2$ by $\psi_{1t}$ yields that
$$
\psi_{1t}^2=\Big(-\f{2\z_x}{3v}+\f{2\t}{3v}\f{\phi_x}{v}-\frac23\bar\t_x\Big(\frac{1}{v}-\frac{1}{\bar v}\Big)
+\frac23\bar v_x\Big(\frac{\t}{v^2}-\frac{\bar\t}{\bar v^2}\Big)
+\Phi_x n_2-\int \xi_1^2\mb{G}_xd\xi - U^{cd}_{1t}-(\bar p-p_+-p_-)_x\Big)\psi_{1t}.
$$
Integrating the above equality with respect to $x, t$ and using
Cauchy inequality and smallness of $\chi$ and $\delta$, one has
\begin{equation}\label{psi-t-E}
\begin{array}{l}
\displaystyle \int_0^t\| \psi_{1t}\|^2 ds \leq C\int_0^t\|(\phi_x,\zeta_x)\|^2ds
+C\int_0^t\||(\bar v_x,\bar\t_x)|(\phi,\zeta)\|^2ds+C(\chi+\delta)\int_0^t\|\Phi_x\|^2ds\\
\di \qquad\qquad\qquad +C\int_0^t\int\int\frac{\nu(| \xi|)}{\mb{M}_\star }|\mathbf{G}_x|^2
d \xi dxds+C\d.
\end{array}
\end{equation}
Similar estimates hold for $\phi_t,\psi_{2t},\psi_{3t}$ and $\zeta_t$.  Therefore, one can arrive at
\begin{equation}\label{t-E}
\begin{array}{l}
\displaystyle \int_0^t\|(\phi_{t},\psi_t,\zeta_t)\|^2 ds \leq C\d+ C\int_0^t\|(\phi_x,\psi_x,\zeta_x)\|^2ds
+C\int_0^t\||(\bar v_x,\bar u_{1x}, \bar\t_x)|(\phi,\zeta)\|^2ds\\
\di \quad +C(\chi+\delta)\int_0^t\|\Phi_x\|^2ds+C\int_0^t\int\int\frac{\nu(| \xi|)}{\mb{M}_\star }|\mathbf{G}_x|^2
d \xi dxds +C\chi\int_0^t\int\int\frac{\nu(| \xi|)}{\mb{M}_\star }|\mathbf{P}_cF_2|^2
d \xi dxds .
\end{array}
\end{equation}
By \eqref{le}, \eqref{important}, \eqref{phi-x-E} and \eqref{t-E}, it holds that
\begin{equation}\label{le1}
\begin{array}{ll}
\di
\|(\phi,\psi,\zeta,\phi_x)(\cdot,t)\|^2+\sum_{|\alpha|=1}\int_0^t\|\partial^\alpha(\phi,\psi,\zeta)\|^2ds\\
\di\quad+\d\int_0^t\int(\phi^2+\psi^2+\z^2)\hat w^2dxds+\int_0^t\int \big((U^r_{1-})_x+(U^r_{1+})_x\big)(\phi^2+\zeta^2)dxds\\
\di \leq C\|(\phi,\psi,\zeta,\phi_x)(\cdot,0)\|^2+C\d+C\sum_{1\leq|\alpha|\leq2}
\int_0^t\int\int\f{\nu(|\xi|)|\partial^\alpha\mb{G}|^2}{\mb{M}_\star }
d\xi
dxds\\
\di
\quad+C(\chi+\delta)\int_0^t\int\int\f{\nu(|\xi|)|(\widetilde{\mb{G}},\mb{P}_c F_2)|^2}{\mb{M}_\star }
d\xi dxds+ C(\chi+\delta) \int_0^t\|(\Phi_x, n_2)\|^2ds\\
\di\quad +C(\chi+\delta)\int_0^t\int\int\frac{\nu(|\xi|)}{\mb{M}_\star }|(\pa_{\xi_1}(\mb{P}_c F_2),\pa_{\xi_1}(\mb{P}_c F_2)_x)|^2d\xi dxds.
\end{array}
\end{equation}

\

\underline{Step 3. Estimation on the non-fluid part.}

\

Next we do the microscopic estimates for the Vlasov-Poisson-Boltzmann system.
Multiplying the equation \eqref{GE} and the equation \eqref{F2-pc-1} by
$\f{\widetilde{\mb{G}}}{\mb{M}_\star }$ and $\f{\mb{P}_c
F_2}{\mb{M}_\star }$ , respectively, one has
\begin{equation}\label{Gle0}
\begin{array}{ll}
\di \Big(\f{|\widetilde{\mb{G}}|^2}{2\mb{M}_\star }\Big)_t-\f{\widetilde{\mb{G}}}{\mb{M}_\star }\mb{L}_{\mb{M}}\widetilde{\mb{G}}=\Big\{-\f{3}{2\t v}\mb{P}_1\Big[\xi_1\big(\f{|\xi-u|^2}{2\t}\zeta_x+\xi\cdot\psi_{x}\big)\mb{M}\Big]\\[3mm]
\di
+\f{u_1}{v}\mb{G}_x-\f1v\mb{P}_1(\xi_1\mb{G}_x)-\f1v\mb{P}_1(\Phi_x\partial_{\xi_1}F_2)+2Q(\mb{G},\mb{G})-\bar{
\mb{G}}_t\Big\}\f{\widetilde{\mb{G}}}{\mb{M}_\star },
\end{array}
\end{equation}
and
\begin{equation}\label{Pc-F2-le}
\begin{array}{ll}
\di \Big(\f{|\mb{P}_c F_2|^2}{2\mb{M}_\star }\Big)_t-\f{\mb{P}_c
F_2}{\mb{M}_\star }\mb{N}_{\mb{M}}(\mb{P_c}F_2)=\Big[-\f{\xi_1}{v}F_{2x}+\f{u_1}{v}F_{2x}-(n_2v\mb{M})_t
\\\di\qquad\qquad\qquad\qquad\qquad -\f1v\mb{P}_c(\Phi_x\partial_{\xi_1}F_1)+2Q(F_2,\mb{G})\Big]\f{\mb{P}_c F_2}{\mb{M}_\star }.
\end{array}
\end{equation}
It can be computed that
\begin{equation}\label{P1-1}
\begin{array}{ll}
\di \mb{P}_1(\Phi_x\partial_{\xi_1}F_2)=\Phi_x\mb{P}_1\big[\partial_{\xi_1}\big(n_2v\mb{M}+\mb{P}_cF_2\big)\big]
=\Phi_x\mb{P}_1\big[\partial_{\xi_1}\big(\mb{P}_cF_2\big)\big]\\
\di\qquad\qquad\qquad =\Phi_x\big(\mb{P}_cF_2\big)_{\xi_1}-\Phi_x\sum_{j=0}^4\int\big(\mb{P}_cF_2\big)_{\xi_1}\f{\chi_j}{\mb{M}}d\xi\chi_j\\
\di\qquad\qquad\qquad =\Phi_x\big(\mb{P}_cF_2\big)_{\xi_1}+\Phi_x\sum_{j=0}^4\int(\mb{P}_cF_2) \big(\f{\chi_j}{\mb{M}}\big)_{\xi_1}d\xi\chi_j,
\end{array}
\end{equation}
and
\begin{equation}\label{Pc-1}
\di \mb{P}_c(\Phi_x\partial_{\xi_1}F_1)=\Phi_x\partial_{\xi_1}F_1-\Phi_x\mb{P}_d\big(\partial_{\xi_1}F_1\big)=\Phi_x\partial_{\xi_1}F_1=\Phi_x \mb{M}_{\xi_1}+\Phi_x \widetilde{\mb{G}}_{\xi_1}+\Phi_x \bar{\mb{G}}_{\xi_1}.
\end{equation}
Substituting \eqref{P1-1} and \eqref{Pc-1} into \eqref{Gle0} and \eqref{Pc-F2-le}, respectively, then summing the resulting equations together and noting that
$$
\Phi_x\big(\mb{P}_cF_2\big)_{\xi_1}\f{\widetilde{\mb{G}}}{\mb{M}_\star }+\Phi_x \widetilde{\mb{G}}_{\xi_1}\f{\mb{P}_c F_2}{\mb{M}_\star }=\Big(\Phi_x\f{\mb{P}_cF_2\cdot\widetilde{\mb{G}}}{\mb{M}_\star }\Big)_{\xi_1}
+\Phi_x\f{\mb{P}_cF_2\cdot\widetilde{\mb{G}}}{\mb{M}_\star ^2}\big(\mb{M}_\star \big)_{\xi_1}.
$$
It holds that
\begin{equation*}
\begin{array}{ll}
\di \Big(\f{|\widetilde{\mb{G}}|^2+|\mb{P}_c F_2|^2}{2\mb{M}_\star }\Big)_t-\f{\widetilde{\mb{G}}}{\mb{M}_\star }\mb{L}_{\mb{M}}\widetilde{\mb{G}}
-\f{\mb{P}_c F_2}{\mb{M}_\star }\mb{N}_{\mb{M}}(\mb{P_c}F_2)+(\cdots)_{\xi_1}
=\Big\{\f{u_1}{v}\mb{G}_x-\f1v \mb{P}_1(\xi_1\mb{G}_x)-\bar{\mb{G}}_t\\[3mm]
\di
+2Q(\mb{G},\mb{G})-\f{3}{2\t v}\mb{P}_1\Big[\xi_1\big(\f{|\xi-u|^2}{2\t}\zeta_x+\xi\cdot\psi_{x}\big)\mb{M}\Big]
-\f{\Phi_x}{v}\sum_{j=0}^4\int(\mb{P}_cF_2) \big(\f{\chi_j}{\mb{M}}\big)_{\xi_1}d\xi\chi_j\Big\}\f{\widetilde{\mb{G}}}{\mb{M}_\star }\\
\di+\Big[-\f{\xi_1}{v}F_{2x}+\f{u_1}{v}F_{2x}-(n_2v\mb{M})_t+2Q(F_2,\mb{G})-\f{\Phi_x}{v} \mb{M}_{\xi_1}
-\f{\Phi_x}{v} \bar{\mb{G}}_{\xi_1}\Big]\f{\mb{P}_c F_2}{\mb{M}_\star }
-\f{\Phi_x}{v}\f{\mb{P}_cF_2\cdot\widetilde{\mb{G}}}{\mb{M}_\star ^2}\big(\mb{M}_\star \big)_{\xi_1}.
\end{array}
\end{equation*}
Integrating the above equation with respect to $\xi, x, t$ yields that
\begin{equation}\label{GE0}
\begin{array}{ll}
\di \int\int \f{|\widetilde{\mb{G}}|^2+|\mb{P}_c F_2|^2}{2\mb{M_\star}}(x,\xi,t)
d\xi dx-\int\int \f{|\widetilde{\mb{G}}|^2+|\mb{P}_c F_2|^2}{2\mb{M}_\star }(x,\xi,0)
d\xi dx\\
\di  -\int_0^t\int\int \Big[\f{\widetilde{\mb{G}}}{\mb{M}_\star }\mb{L}_{\mb{M}}\widetilde{\mb{G}}+\f{\mb{P}_c
F_2}{\mb{M}_\star }\mb{N}_{\mb{M}}(\mb{P_c}F_2)\Big]
d\xi dxds
=\int_0^t\int\int\Bigg\{\Big\{\f{u_1}{v}\mb{G}_x-\f1v\mb{P}_1(\xi_1\mb{G}_x)\\[3mm]
\di-\bar{\mb{G}}_t+2Q(\mb{G},\mb{G})-\f{3}{2\t v}\mb{P}_1\Big[\xi_1\big(\f{|\xi-u|^2}{2\t}\zeta_x+\xi\cdot\psi_{x}\big)\mb{M}\Big]
-\f{\Phi_x}{v}\sum_{j=0}^4\int(\mb{P}_cF_2) \big(\f{\chi_j}{\mb{M}}\big)_{\xi_1}d\xi\chi_j\Big\}\f{\widetilde{\mb{G}}}{\mb{M}_\star }\\
\di+\Big[-\f{\xi_1}{v}F_{2x}+\f{u_1}{v}F_{2x}-(n_2v\mb{M})_t+2Q(F_2,\mb{G})-\f{\Phi_x}{v}\mb{M}_{\xi_1}-\f{\Phi_x}{v}\bar{\mb{G}}_{\xi_1}\Big]\f{\mb{P}_c F_2}{\mb{M}_\star }\\
\di-\f{\Phi_x}{v}\f{\mb{P}_cF_2\cdot\widetilde{\mb{G}}}{\mb{M}_\star ^2}\big(\mb{M}_\star \big)_{\xi_1}\Bigg\}d\xi dxds :=\sum_{i=1}^{13}J_i.
\end{array}
\end{equation}
Now we calculate the right hand side of \eqref{GE0} terms by terms. First, it holds
that
\begin{equation}\label{J1}
\begin{array}{ll}
\di |J_1|+|J_2|\leq \f{\tilde\sigma}{16}\int_0^t\int\int\f{\nu(|\xi|)|\widetilde{\mb{G}}|^2}{\mb{M}_\star }d\xi dxds
+C\int_0^t\int\int\f{\nu(|\xi|)|\mb{G}_x|^2}{\mb{M}_\star }d\xi dxds,
\end{array}
\end{equation}
and
\begin{equation}\label{J3}
\begin{array}{ll}
\di |J_3| \leq\f{\tilde\sigma}{16}\int_0^t\int\int\f{\nu(|\xi|)|\widetilde{\mb{G}}|^2}{\mb{M}_\star }d\xi dxds
+C\int_0^t\int\int \f{\nu(|\xi|)|\bar{\mb{G}}_t|^2}{\mb{M}_\star }d\xi dxds\\
\di \quad \leq
\f{\tilde\sigma}{16}\int_0^t\int\int\f{\nu(|\xi|)|\widetilde{\mb{G}}|^2}{\mb{M}_\star }d\xi dxds
+C\int_0^t\int\Big[|(\bar\t_{xt},\bar u_{1xt})|^2+|( v_t, u_t, \t_t)|^2|(\bar\t_x,\bar
u_{1x})|^2\Big]dxds\\
\di\quad\leq
\f{\tilde\sigma}{16}\int_0^t\int\int\f{\nu(|\xi|)|\widetilde{\mb{G}}|^2}{\mb{M}_\star }d\xi dxds
+C\d\int_0^t\|(\phi_t,\psi_t,\zeta_t)\|^2ds+C\delta.
\end{array}
\end{equation}
Then it holds that
\begin{equation}\label{J4}
\begin{array}{ll}
\di |J_4| \leq
\f{\tilde\sigma}{16}\int_0^t\int\int\f{\nu(|\xi |)|\widetilde{\mb{G}}|^2}{\mb{M}_\star }d\xi dxds+C\int_0^t\int\int
\f{\nu(|\xi |)^{-1}|Q(\mb{G},\mb{G})|^2}{\mb{M}_\star }d\xi dxds\\
\di \quad \leq
\f{\tilde\sigma}{16}\int_0^t\int\int\f{\nu(|\xi |)|\widetilde{\mb{G}}|^2}{\mb{M}_\star }d\xi dxds+C\int_0^t\int\Big(\int\f{\nu(|\xi |)|\mb{G}|^2}{\mb{M}_\star }
d\xi \int\f{|\mb{G}|^2}{\mb{M}_\star } d\xi \Big)dxds\\
\di\quad \leq
\f{\tilde\sigma}{16}\int_0^t\int\int\f{\nu(|\xi |)|\widetilde{\mb{G}}|^2}{\mb{M}_\star }d\xi dxds+C\int_0^t\int\Big(\int\f{\nu(|\xi |)|\widetilde{\mb{G}}|^2}{\mb{M}_\star }
d\xi \int\f{|\widetilde{\mb{G}}|^2}{\mb{M}_\star }
d\xi \Big)dxds\\
\di\qquad +C\int_0^t\int|(\bar\t_x,\bar
u_{1x})|^2\Big(\int\f{\nu(|\xi |)|\widetilde{\mb{G}}|^2}{\mb{M}_\star }d\xi\Big)
 dxds+C\int_0^t\int|(\bar\t_x,\bar u_{1x})|^4dxds\\
\di\quad\leq
(\f{\tilde\sigma}{16}+C(\chi+\d))\int_0^t\int\int\f{\nu(|\xi |)|\widetilde{\mb{G}}|^2}{\mb{M}_\star }d\xi dxds
+C\delta,
\end{array}
\end{equation}
\begin{equation}\label{J5}
\begin{array}{ll}
\di |J_5| \leq
\f{\tilde\sigma}{16}\int_0^t\int\int\f{\nu(|\xi|)|\widetilde{\mb{G}}|^2}{\mb{M}_\star }d\xi dxds
+C\int_0^t\|(\psi_x,\z_x)\|^2ds,
\end{array}
\end{equation}
and
\begin{equation}\label{J6}
\begin{array}{ll}
\di
|J_6| \leq C\int_0^t\int |\Phi_x|\Big(\int\f{\nu(|\xi |)|\mb{P}_cF_2|^2}{\mb{M}_\star }d\xi \Big)^{\f12}
\Big(\int\f{\nu(|\xi |)|\widetilde{\mb{G}}|^2}{\mb{M}_\star }d\xi \Big)^{\f12} dxds\\
\di\qquad \leq C\chi
\int_0^t\int\int\f{\nu(|\xi |)|\mb{P}_cF_2|^2}{\mb{M}_\star }d\xi dxds+C\chi\int_0^t\int\int\f{\nu(|\xi |)|\widetilde{\mb{G}}|^2}{\mb{M}_\star }d\xi dxds.
\end{array}
\end{equation}
Note that
$$
\begin{array}{ll}
\di -\f{\xi_1}{v}F_{2x}\f{\mb{P}_c F_2}{\mb{M}_\star }=-\f{\xi_1}{v}\big(n_2v\mb{M}+\mb{P}_cF_2\big)_x\f{\mb{P}_c F_2}{\mb{M}_\star }\\[2mm]
\di=-n_{2x}\f{\mb{M}}{\mb{M}_\star }\xi_1\mb{P}_c F_2-\f{n_2}{v}(v\mb{M})_x\f{\xi_1\mb{P}_c F_2}{\mb{M}_\star }
-\f{v_x}{v}\f{\xi_1|\mb{P}_c F_2|^2}{2\mb{M}_\star }-\Big(\f{\xi_1}{v}\f{|\mb{P}_c F_2|^2}{2\mb{M}_\star }\Big)_x\\[2mm]
\di=-\Big(\f{\mb{M}}{\mb{M}_\star }-1\Big)n_{2x}\xi_1\mb{P}_c F_2-\f{n_2}{v}(v\mb{M})_x\f{\xi_1\mb{P}_c F_2}{\mb{M}_\star }
-\f{v_x}{v}\f{\xi_1|\mb{P}_c F_2|^2}{2\mb{M}_\star }-n_{2x}\xi_1\mb{P}_c F_2+(\cdots)_x.
\end{array}
$$
Then integration by parts yields that
\begin{equation}\label{J7}
\begin{array}{ll}
\di J_7=-\int_0^t\int\int \f{\xi_1}{v} F_{2x}\f{\mb{P}_c F_2}{\mb{M}_\star }d\xi dxds \\
\di\quad =-\int_0^t\int\int \Big[\Big(\f{\mb{M}}{\mb{M}_\star }-1\Big)n_{2x}\xi_1\mb{P}_c F_2
+\f{n_2}{v}(v\mb{M})_x\f{\xi_1\mb{P}_c F_2}{\mb{M}_\star }
+\f{v_x}{v}\f{\xi_1|\mb{P}_c F_2|^2}{2\mb{M}_\star }\Big]d\xi dxds\\
\di\quad\quad  +\int_0^t\int\int n_{2}\Big(\xi_1\mb{P}_c F_2 d\xi\Big)_x dxds.
\end{array}
\end{equation}
It follows from \eqref{n22} that
\begin{equation}\label{J71}
\begin{array}{ll}
\di \int_0^t\int\int n_{2}\Big(\xi_1\mb{P}_c F_2 d\xi\Big)_x dxds
=-\int_0^t\int \Big(v n_2 n_{2t}+u_{1x}n^2_2\Big)dxds\\
\di=-\int\f{vn_2^2}{2}(x,t)dx+\int\f{vn_{2}^2}{2}(x,0)dx
-\int_0^t\int \f{u_{1x}}{2}n_2^2dxds.
\end{array}
\end{equation}
Substituting \eqref{J71} into \eqref{J7} and estimating the other
terms yield that
\begin{equation}\label{J7-E}
\begin{array}{ll}
\di J_7+\int\f{vn_2^2}{2}(x,t)dx \leq C\|n_{20}\|^2+C(\chi+\d)\int_0^t\|n_2\|^2ds
+C\eta_0\int_0^t\|n_{2x}\|^2ds\\
\di\qquad+C(\chi+\delta+\eta_0)\int_0^t\int\int\f{\nu(|\xi |)|\mb{P}_c F_2|^2}{\mb{M}_\star }d\xi dxds.
\end{array}
\end{equation}
Similarly, it holds that
\begin{equation}\label{J8}
\begin{array}{ll}
\di |J_8|\leq C(\chi+\d)\int_0^t\|n_2\|^2ds+C\eta_0\int_0^t\|n_{2x}\|^2ds
+C(\chi+\delta+\eta_0)\int_0^t\int\int\f{\nu(|\xi |)|\mb{P}_c F_2|^2}{\mb{M}_\star }d\xi dxds.
\end{array}
\end{equation}
\begin{equation}\label{J9}
\begin{array}{ll}
\di |J_9|=|\int_0^t\int\int \Big[n_{2t}v\f{\mb{M}}{\mb{M}_\star }\mb{P}_c F_2+n_2(v\mb{M})_t\f{\mb{P}_c F_2}{\mb{M}_\star }\Big] d\xi dxds|\\
\di=|\int_0^t\int\int \Big[n_{2t}v\Big(\f{\mb{M}}{\mb{M}_\star }-1\Big)\mb{P}_c F_2+n_2(v\mb{M})_t\f{\mb{P}_c F_2}{\mb{M}_\star }\Big] d\xi dxds|\\
\di\leq C(\chi+\d)\int_0^t\|n_2\|^2ds+C\eta_0\int_0^t\|n_{2t}\|^2ds
+C(\chi+\delta+\eta_0)\int_0^t\int\int\f{\nu(|\xi |)|\mb{P}_c F_2|^2}{\mb{M}_\star }d\xi dxds.
\end{array}
\end{equation}
Then one has
\begin{equation}\label{J10}
\begin{array}{ll}
\di |J_{10}| \leq \f{\tilde\sigma}{16}\int_0^t\int\int\frac{\nu(|\xi |)|\mb{P}_c
F_2|^2}{\mb{M}_\star }d\xi dxds+C\int_0^t
\int\int\frac{\nu(|\xi |)^{-1}|Q(F_2,\mb{G})|^2}{\mb{M}_\star }d\xi dxds\\
\di \leq \f{\tilde\sigma}{16}\int_0^t\int\int\frac{\nu(|\xi |)|\mb{P}_c
F_2|^2}{\mb{M}_\star }d\xi dxds+C\int_0^t
\int\int\frac{\nu(|\xi |)|F_2|^2}{\mb{M}_\star }d\xi \int
\frac{|\mb{G}|^2}{\mb{M}_\star }d\xi  dxds\\
\di\qquad+C\int_0^t \int\int\frac{|F_2|^2}{\mb{M}_\star }d\xi \int
\frac{\nu(|\xi |)|\mb{G}|^2}{\mb{M}_\star }d\xi  dxds\\
\di\leq (\f{\tilde\sigma}{16}+C(\chi+\d))\int_0^t\int\int\frac{\nu(|\xi |)|\mb{P}_c
F_2|^2}{\mb{M}_\star }d\xi dxds+C(\chi+\delta)\int_0^t
\Big[\|n_2\|^2+\int\int\frac{\nu(|\xi |)|\widetilde{\mb{G}}|^2}{\mb{M}_\star }d\xi dx\Big]ds.
\end{array}
\end{equation}
By \eqref{wt-3.31}, one has
\begin{equation}\label{J11}
\begin{array}{ll}
\di J_{11}=\int_0^t\int\int \f{\Phi_x}{v}\f{\xi_1-u_1}{R\t}\mb{M}\f{\mb{P}_c
F_2}{\mb{M}_\star } d\xi dxds\\
\di \quad=\int_0^t\int\int
\Phi_x\f{\xi_1-u_1}{R\t v}\Big(\f{\mb{M}}{\mb{M}_\star }-1\Big)\mb{P}_c F_2
d\xi dxds+\int_0^t\int\int \Phi_x\f{\xi_1-u_1}{R\t v}\mb{P}_c F_2 d\xi dxds
\\
\di \quad=\int_0^t\int\int
\Phi_x\f{\xi_1-u_1}{R\t v}\Big(\f{\mb{M}}{\mb{M}_\star }-1\Big)\mb{P}_c F_2
d\xi dxds+\int_0^t\int\f{\Phi_x}{R\t v}\Big(\int \xi_1\mb{P}_c F_2
d\xi\Big)dxds\\
\di\quad =\int_0^t\int\int
\Phi_x\f{\xi_1-u_1}{R\t v}\Big(\f{\mb{M}}{\mb{M}_\star }-1\Big)\mb{P}_c F_2
d\xi dxds-\int_0^t\int\f{\Phi_x}{R\t v}\Big(\f{\Phi_{x}}{2v}\Big)_tdxds\\
\di\quad \leq -\int\f{\Phi_x^2}{4R\t v^2}(x,t)dx+C\|\Phi_{0x}\|^2
+C(\chi+\delta+\eta_0)\int_0^t\|\Phi_x\|^2ds\\
\di\qquad+C\eta_0\int_0^t\int\int \f{\nu(|\xi|)|\mb{P}_c F_2|^2}{\mb{M}_\star }d\xi dxds,
\end{array}
\end{equation}
where in the last inequality one has used the fact that
$$
\begin{array}{ll}
\di
-\int_0^t\int\f{\Phi_x}{R\t v}\Big(\f{\Phi_x}{2v}\Big)_tdxds=-\int\f{\Phi_x^2}{4R\t v^2}(x,t)dx+\int\f{\Phi_{x}^2}{4R\t v^2}(x,0)dx
-\int_0^t\int\f{3\t_t\Phi_x^2}{8\t^2 v^2}dxds\\
\di\qquad\qquad  \leq -\int\f{\Phi_x^2}{4R\t v^2}(x,t)dx
+C\|\Phi_{0x}\|^2+C(\chi+\delta)\int_0^t\|\Phi_x\|^2ds.
\end{array}
$$
Moreover, one has
\begin{equation}\label{J12}
\begin{array}{ll}
\di |J_{12}|\leq C\d \int_0^t\int\int\frac{\nu(|\xi|)|\mb{P}_c
F_2|^2}{\mb{M}_\star }d\xi dxds+C\delta\int_0^t\|\Phi_x\|^2ds,
\end{array}
\end{equation}
and
\begin{equation}\label{J13}
\begin{array}{ll}
\di |J_{13}| \leq C\chi \int_0^t\int\int\frac{\nu(|\xi|)|\mb{P}_c
F_2|^2}{\mb{M}_\star }d\xi dxds+ C\chi\int_0^t\int\int\frac{\nu(|\xi|)|\widetilde{\mb{G}}|^2}{\mb{M}_\star }d\xi dxds.
\end{array}
\end{equation}

Substituting \eqref{J1}-\eqref{J13} into \eqref{GE0}, using the lower order estimates in \eqref{le} and thenchoosing
$\chi,\delta$ sufficiently small give that
\begin{equation}\label{FE}
\begin{array}{ll}
\di \int\int \f{|(\widetilde{\mb{G}},\mb{P}_c
F_2)|^2}{\mb{M}_\star }(x,\xi,t)d\xi dx+\|(\Phi_x,n_2)(\cdot,t)\|^2+\int_0^t\int\int\frac{\nu(|\xi|)|\big(\widetilde{\mb{G}},
\mb{P}_c
F_2\big)|^2}{\mb{M}_\star }d\xi dxds\\
\di \leq C\int\int \f{|(\widetilde{\mb{G}}, \mb{P}_c
F_2)|^2}{\mb{M}_\star }(x,\xi,0)d\xi dx+C\|(\Phi_{0x},n_{20})\|^2+C\|(\phi_0,\psi_0,\zeta_0)\|^2+
C\delta\\
\di +C(\chi+\delta+\eta_0)\int_0^t\|(\Phi_x,n_{2x},n_{2t})\|^2
ds+C(\chi+\delta) \int_0^t
\Big[\|n_2\|^2+\sum_{|\alpha'|=1}\|\partial^{\alpha'}(\phi,\psi,\zeta)\|^2\Big]ds\\
\di +C(\chi+\d)\int_0^t\int\int\f{\nu(|\xi|)|\pa_{\xi_1}(\mb{P}_c F_2)|^2}{\mb{M}_\star } d\xi dxds
+C\sum_{|\alpha^\prime|=1}\int_0^t\int\int\frac{\nu(|\xi|)|\partial^{\alpha^\prime} \mb{G}|^2}{\mb{M}_\star }d\xi dxds.
\end{array}
\end{equation}

\

\underline{Step 4. Estimation on the Poisson term.}

\

Now we estimate the Poisson term in the electric fields, which is one of the key ingredients of the present paper. Multiplying the equation
\eqref{n2} by $-v\Phi$  and noting that
\begin{equation}
\begin{array}{l}
\di-v\Phi n_{2t}-u_{1x}\Phi n_2=-\f{\Phi}{2}\Big(\f{\Phi_x}{v}\Big)_{xt}
+\f{v_t\Phi}{2v}\Big(\f{\Phi_x}{v}\Big)_x-\f{u_{1x}\Phi}{2v}\Big(\f{\Phi_x}{v}\Big)_x\\[2mm]
\di=(\cdots)_x+\f{\Phi_x}{2}\Big(\f{\Phi_x}{v}\Big)_t=(\cdots)_x+\f{\Phi_x\Phi_{xt}}{2v}-\f{v_t\Phi_x^2}{2v^2}=(\cdots)_x+\Big(\f{\Phi_x^2}{4v}\Big)_t-\f{v_t\Phi_x^2}{4v^2}
\end{array}
\end{equation}
and
\begin{equation}
\begin{array}{l}
\di\Phi\Big(\f{\k_1(\t)}{v}(n_2v)_x\Big)_x=(\cdots)_x-\Phi_x\f{\k_1(\t)}{v}(n_2 v)_x=(\cdots)_x+\Big(\k_1(\t)\f{\Phi_x}{v}\Big)_xn_2 v\\[2mm]
\di\qquad\qquad\qquad\qquad=(\cdots)_x+\k_1(\t)_x\Phi_x n_2+\k_1(\t)\Big(\f{\Phi_x}{v}\Big)_xn_2 v  \\[2mm]
\di\qquad\qquad\qquad\qquad=(\cdots)_x+\k_1(\t)_x\Phi_x n_2+2\k_1(\t)n_2^2v^2,
\end{array}
\end{equation}
which yields that
\begin{equation}
\begin{array}{ll}
\di
\Big(\f{\Phi_x^2}{4v}\Big)_t+\f{3\k_1(\t)}{2\t v}\Phi_x^2  +2\k_1(\t)n_2^2v^2+(\cdots)_x=\f{v_t\Phi_x^2}{4v^2}-\k_1(\t)_x\Phi_x n_2\\
\di -\Phi_x n_2\int \xi_1 \mb{N}_{\mb{M}}^{-1}\big[\mb{P}_c(\xi_1\mb{M}_x)\big]d\xi
-\f{\Phi_x}{v}\int \xi_1 \mb{N}_{\mb{M}}^{-1}\big[\mb{P}_c(\xi_1(\mb{P}_cF_2)_x)\big]d\xi
-\f{\Phi_x^2}{v}\int \xi_1 \mb{N}_{\mb{M}}^{-1}\big[\mb{G}_{\xi_1}\big]d\xi\\
\di -\Phi_x\int \xi_1 \mb{N}_{\mb{M}}^{-1}\Big[(\mb{P}_c F_2)_t-\f{u_1}{v}(\mb{P}_c F_2)_x
+\big((\mb{M}v)_t-\f{u_1}{v}(\mb{M}v)_x\big)n_2-2Q(F_2,\mb{G})\Big]d\xi.
\end{array}
\end{equation}
Integrating the above equation with respect to $x,t$ gives
\begin{equation}\label{Phi-1}
\begin{array}{ll}
\di \|\Phi_x(\cdot,t)\|^2+\int_0^t\|(\Phi_x,n_2)\|^2ds\leq
C\|\Phi_{0x}\|^2+C(\chi+\delta)\int_0^t \|(\Phi_{x},n_2)\|^2ds\\
\di   +C\int_0^t\int|\Phi_x n_2\int \xi_1 \mb{N}_{\mb{M}}^{-1}\big[\mb{P}_c(\xi_1\mb{M}_x) \big]d\xi|dxds
+C\int_0^t\int|\f{\Phi_x}{v}\int \xi_1\mb{N}_{\mb{M}}^{-1}\big[\mb{P}_c(\xi_1 (\mb{P}_cF_2)_x)\big]d\xi|dxds\\
\di
+C\int_0^t\int| \f{\Phi_x^2}{v}\int \xi_1 \mb{N}_{\mb{M}}^{-1}\big[\mb{G}_{\xi_1}\big]d\xi|dxds
+C\int_0^t\int|\Phi_x\int \xi_1 \mb{N}_{\mb{M}}^{-1}\big[ Q(F_2,\mb{G})\big]d\xi|dxds\\
 \di +C\int_0^t\int|\Phi_x\int \xi_1 \mb{N}_{\mb{M}}^{-1}\Big[(\mb{P}_c F_2)_t-\f{u_1}{v}(\mb{P}_c F_2)_x
+\big((\mb{M}v)_t-\f{u_1}{v}(\mb{M}v)_x\big)~n_2\Big]d\xi|dxds\\
\di :=C\|\Phi_{0x}\|^2+C(\chi+\delta)\int_0^t \|(\Phi_{x},n_2)\|^2ds+\sum_{i=1}^5K_i.
\end{array}
\end{equation}
Now we estimate $K_i~(i=1,\cdots,5)$ in \eqref{Phi-1} one by one. First, it holds that
\begin{equation}\label{K1}
\begin{array}{ll}
\di K_1\leq C \int_0^t\int |n_2||\Phi_x|\big(\int\f{\nu(|\xi|)|\mb{N}_{\mb{M}}^{-1}[\mb{P}_c(\xi_1\mb{M}_x)]|^2}{\mb{M}_\star } d\xi\big)^{\f12}\big(\int \nu(|\xi|)^{-1}\xi_1^2\mb{M}_\star d\xi\big)^{\f12}dxds\\
\di \qquad \leq C\int_0^t\int  |n_2||\Phi_x| \big(\int\f{\nu(|\xi|)^{-1}|\mb{P}_c(\xi_1\mb{M}_x)|^2}{\mb{M}_\star } d\xi\big)^{\f12}dxds\\
\di \qquad \leq C\int_0^t\int  |n_2||\Phi_x| \big(\int\f{\nu(|\xi|)^{-1}|\xi_1\mb{M}_x-\mb{M}v(\f{u_1}{v})_x|^2}{\mb{M}_\star } d\xi\big)^{\f12}dxds\\
\di \qquad \leq C\int_0^t\int  |n_2||\Phi_x| |(v_x,u_x,\t_x)|dxds
\leq C(\chi+\delta)\int_0^t\|(n_2,\Phi_x)\|^2dxds,
\end{array}
\end{equation}
and
\begin{equation}\label{K2}
\begin{array}{ll}
\di K_2\leq C \int_0^t\int |\Phi_x|\big(\int\f{\nu(|\xi|)|\mb{N}_{\mb{M}}^{-1}\big[\mb{P}_c(\xi_1 (\mb{P}_cF_2)_x)\big]|^2}{\mb{M}_{[v_\star,u_\star,2\t_\star]}} d\xi\big)^{\f12}\big(\int \nu(|\xi|)^{-1}\xi_1^2\mb{M}_{[v_\star,u_\star,2\t_\star]}d\xi\big)^{\f12}dxds\\
\di \qquad \leq C\int_0^t\int  |\Phi_x| \big(\int\f{\nu(|\xi|)^{-1}|\mb{P}_c(\xi_1 (\mb{P}_cF_2)_x)|^2}{\mb{M}_{[v_\star,u_\star,2\t_\star]}} d\xi\big)^{\f12}dxds\\
\di \qquad \leq C\int_0^t\int  |\Phi_x| \big(\int\f{\nu(|\xi|)^{-1}|(\mb{P}_cF_2)_x|^2}{\mb{M}_\star } d\xi\big)^{\f12}dxds \\
\di \qquad \leq \f18\int_0^t\|\Phi_x\|^2ds+C~\int_0^t\int
\int\f{\nu(|\xi|)|(\mb{P}_cF_2)_x|^2}{\mb{M}_\star } d\xi dxds.
\end{array}
\end{equation}
Then one has
\begin{equation}\label{K3}
\begin{array}{ll}
\di K_3\leq C \int_0^t\int |\Phi_x|^2\big(\int\f{\nu(|\xi|)^{-1}|\mb{G}_{\xi_1}|^2}{\mb{M}_\star } d\xi\big)^{\f12}dxds\\
\di\qquad \leq C \int_0^t\int |\Phi_x|^2\Big[\big(\int\f{\nu(|\xi|)|\widetilde{\mb{G}}_{\xi_1}|^2}{\mb{M}_\star } d\xi\big)^{\f12}
+|(\bar\t_x,\bar u_{1x})|\Big]dxds\\
\di \qquad \leq
C(\chi+\delta)\int_0^t\Big[\|\Phi_x\|^2+\int\int\f{\nu(|\xi|)|\widetilde{\mb{G}}_{\xi_1}|^2}{\mb{M}_\star }
d\xi dx\Big]ds,
\end{array}
\end{equation}
\begin{equation}\label{K4}
\begin{array}{ll}
\di K_4\leq C \int_0^t\int |\Phi_x|\big(\int\f{\nu(|\xi|)^{-1}|Q(F_2,\mb{G})|^2}{\mb{M}_\star } d\xi\big)^{\f12}dxds\\
\di \qquad \leq C \int_0^t\int |\Phi_x|\Big(\int\f{\nu(|\xi|)|F_2|^2}{\mb{M}} d\xi\int\f{\nu(|\xi|)|\mb{G}|^2}{\mb{M}} d\xi\Big)^{\f12}dxds\\
\di\qquad  \leq C \int_0^t\int |\Phi_x|\Big(\int\f{|F_2|^2}{\mb{M}_\star } d\xi\int\f{|\mb{G}|^2}{\mb{M}_\star } d\xi\Big)^{\f12}dxds\\
\di\qquad \leq C \int_0^t\int |\Phi_x|\Big(|n_2|^2+\int\f{|\mb{P}_c F_2|^2}{\mb{M}_\star } d\xi\Big)^{\f12}
\Big(|(\bar\t_x,\bar u_{1x})|^2+\int\f{|\widetilde{\mb{G}}|^2}{\mb{M}_\star } d\xi\Big)^{\f12}dxds\\
\di\qquad \leq C(\chi+\delta)\int_0^t\Big[\|(\Phi_x,n_2)\|^2
+\int\int\f{\nu(|\xi|)|(\widetilde{\mb{G}},\mb{P}_cF_2)|^2}{\mb{M}_\star } d\xi dx\Big] ds
\end{array}
\end{equation}
and
\begin{equation}\label{K5}
\begin{array}{ll}
\di K_5\leq  \f18\int_0^t\|\Phi_x\|^2ds+C\sum_{|\alpha'|=1}\int_0^t\int |\Phi_x||n_2||\pa^{\alpha'}(v,u,\t)|dxds\\
\di\qquad+C\sum_{|\alpha'|=1}\int_0^t\int\int\f{\nu(|\xi|)|\pa^{\alpha'}(\mb{P}_c F_2)|^2}{\mb{M}_\star } d\xi dxds\\
\di\qquad \leq \f18\int_0^t\|\Phi_x\|^2ds+C(\chi+\d)\int_0^t\|(\Phi_x,n_2)\|^2ds
+C\sum_{|\alpha'|=1}\int_0^t\int\int\f{\nu(|\xi|)|\pa^{\alpha'}(\mb{P}_c F_2)|^2}{\mb{M}_\star } d\xi dxds.
\end{array}
\end{equation}
Substituting \eqref{K1}-\eqref{K5} into \eqref{Phi-1} and choosing
$\chi,\delta$ sufficiently small yield that
\begin{equation}\label{Phi-E}
\begin{array}{ll}
\di \|\Phi_x(\cdot,t)\|^2+\int_0^t\|(\Phi_x,n_2)\|^2ds\leq
C\|\Phi_{0x}\|^2+C\sum_{|\alpha^\prime|=1}\int_0^t\int
\int\f{\nu(|\xi|)|\partial^{\alpha^\prime}(\mb{P}_cF_2)|^2}{\mb{M}_\star }
d\xi dxds\\
\di \qquad\qquad +C(\chi+\delta)\int_0^t\int
\int\f{\nu(|\xi|)|(\widetilde{\mb{G}},\mb{P}_cF_2,\widetilde{\mb{G}}_{\xi_1})|^2}{\mb{M}_\star } d\xi dxds.
\end{array}
\end{equation}
Multiplying the equation \eqref{n2} by $vn_2$ and recalling \eqref{vpb-l}$_3$ yield that
\begin{equation}\label{n2-1}
\begin{array}{ll}
\di
\Big(\f{vn_2^2}{2}\Big)_t+\f{3\k_1(\t)}{\t}vn_2^2+\k_1(\t)n_{2x}^2=(\cdots)_x-\f{u_{1x}}{2}n^2_2-\Big(\f{3k_1(\t)}{2\t}\Big)_x\f{\Phi_x n_2}{v}
-\f{\k_1(\t)v_x}{v}n_2n_{2x}\\
\di
=n_{2x}n_2\int \xi_1 \mb{N}_{\mb{M}}^{-1}\big[\mb{P}_c(\xi_1\mb{M}_x)\big]d\xi+\f{n_{2x}}{v}\int \xi_1 \mb{N}_{\mb{M}}^{-1}\big[\mb{P}_c(\xi_1
(\mb{P}_cF_2)_x)\big]d\xi
+\f{n_{2x}\Phi_x}{v}\int \xi_1 \mb{N}_{\mb{M}}^{-1}\big[\mb{G}_{\xi_1}\big]d\xi\\
\di
+n_{2x}\int \xi_1 \mb{N}_{\mb{M}}^{-1}\Big[(\mb{P}_c F_2)_t-\f{u_1}{v}(\mb{P}_c F_2)_x+\big((\mb{M}v)_t-\f{u_1}{v}(\mb{M}v)_x\big)n_2-2Q(F_2,\mb{G})\Big]d\xi.
\end{array}
\end{equation}
Integrating the above equation with respect to $x,t$ and using similar estimates as for $K_i(i=1,\cdots, 5)$, one has
\begin{equation}\label{n2-E}
\begin{array}{ll}
\di \|n_{2}(\cdot,t)\|^2+\int_0^t\|(n_2,n_{2x})\|^2ds\leq C\|n_{20}\|^2+ C(\chi+\delta)\int_0^t \|\Phi_x\|^2ds\\
\di
+C\sum_{|\alpha^\prime|=1}\int_0^t\int
\int\f{\nu(|\xi|)|\partial^{\alpha^\prime}(\mb{P}_cF_2)|^2}{\mb{M}_\star }
d\xi dxds +C(\chi+\delta)\int_0^t\int
\int\f{\nu(|\xi|)|(\widetilde{\mb{G}},\mb{P}_cF_2,\widetilde{\mb{G}}_{\xi_1})|^2}{\mb{M}_\star } d\xi dxds,
\end{array}
\end{equation}
which together with \eqref{Phi-E} and \eqref{FE}, we have
\begin{equation}
\begin{array}{ll}
\di \|(\Phi_x,n_2)(\cdot,t)\|^2+\int\int\f{|(\widetilde{\mb{G}},\mb{P}_cF_2)|^2}{\mb{M}_\star }(x,\xi,t)d\xi dx
+\int_0^t\|(\Phi_x,n_2,n_{2x})\|^2ds\\
\di\quad+\int_0^t\int\int\f{\nu(|\xi|)|(\widetilde{\mb{G}},\mb{P}_cF_2)|^2}{\mb{M}_\star }d\xi dxds
\leq C\|(\phi_0,\psi_0,\z_0,\Phi_{0x},n_{20})\|^2+C\d\\
\di\quad+\int\int\f{|(\widetilde{\mb{G}},\mb{P}_cF_2)|^2}{\mb{M}_\star }(x,\xi,0)d\xi dx
+C(\chi+\d+\eta_0)\int_0^t\|n_{2t}\|^2ds\\
\di\quad+C(\chi+\d)\sum_{|\alpha'|=1}\int_0^t\|\pa^{\alpha'}(\phi,\psi,\z)\|^2ds
+C(\chi+\d)\int_0^t\int\int\f{\nu(|\xi|)|(\widetilde{\mb{G}},\mb{P}_cF_2)_{\xi_1}|^2}{\mb{M}_\star }d\xi dxds\\
\di\quad+C\sum_{|\alpha'|=1}\int_0^t\int\int\f{\nu(|\xi|)|\pa^{\alpha'}(\mb{G},\mb{P}_cF_2)|^2}{\mb{M}_\star }d\xi dxds.
\end{array}
\end{equation}
By the equation \eqref{n22}, one has
\begin{equation}\label{n2t-E}
\begin{array}{ll}
\di \int_0^t\|n_{2t}\|^2d \leq C(\chi+\d)\int_0^t\|n_{2}\|^2ds+C\int_0^t\int\int\f{\nu(|\xi|)|(\mb{P}_cF_2)_{x}|^2}{\mb{M}_\star }d\xi dxds.
\end{array}
\end{equation}
By \eqref{wt-3.31}, it holds that
\begin{equation}\label{phixt}
 \int_0^t\|\Phi_{xt}\|^2ds \leq
C(\chi+\d)\int_0^t\|\Phi_x\|^2ds+C\int_0^t\int\int\f{\nu(|\xi|)|\mb{P}_c F_2|^2}{\mb{M}_\star }d\xi dxds.
\end{equation}
In summary, collecting all the above lower order estimates and choosing suitably small $\chi$, $\delta$ and $\eta_0$, we arrive that the estimate \eqref{LE-A} in Proposition \ref{Prop3.1}.

\section{Higher order estimates}
\setcounter{equation}{0}
\begin{proposition} \label{Prop3.2}
For each $(v_\pm,u_\pm,\t_\pm)$, there exists a positive uniform-in-time constant $C$ such that, if
$\d\leq \d_0$ for a suitably small positive constant $\d_0$, then it holds that
\begin{equation}\label{high}
\begin{array}{ll}
\di \|(\phi_x,\psi_x,\zeta_x,n_{2x})(\cdot,t)\|^2
+\sum_{1\leq|\beta|\leq 2}\int\int\f{|\pa^{\beta}(\widetilde{\mb{G}},\mb{P}_c F_2)|^2}{\mb{M}_\star }(x,\xi,t)d\xi dx\\
\di+\sum_{|\alpha'|=1,0\leq |\beta'|\leq 1}\int\int\f{|\pa^{\alpha'}\pa^{\beta'}(\mb{G},\mb{P}_c F_2)|^2}{\mb{M}_\star }(x,\xi,t)d\xi dx
+\sum_{|\alpha|= 2}\int\int\f{|\pa^{\alpha}(F_1,F_2)|^2}{\mb{M}_\star }(x,\xi,t)d\xi dx\\
\di +\sum_{|\alpha|= 2}\int_0^t\|\pa^{\alpha}(\phi,\psi,\z,n_2)\|^2ds
+ \sum_{1\leq|\beta|\leq 2}\int_0^t\int\int\f{|\pa^{\beta}(\widetilde{\mb{G}},\mb{P}_c F_2)|^2}{\mb{M}_\star }d\xi dxds\\
\di+\sum_{|\alpha'|=1,0\leq |\beta'|\leq 1}\int_0^t\int\int\f{|\pa^{\alpha'}\pa^{\beta'}(\mb{G},\mb{P}_c F_2)|^2}{\mb{M}_\star }d\xi dxds
+\sum_{|\alpha|= 2}\int_0^t\int\int\f{|\pa^{\alpha}(\mb{G},\mb{P}_c F_2)|^2}{\mb{M}_\star }d\xi dxds\\
\di \leq C\mathcal{N}^2(0)+C \d^{\frac12}+C(\chi+\d+\eta_0)\int_0^t\Big[\sum_{|\alpha'|=1}\|\pa^{\alpha'}(\phi,\psi,\z,n_2)\|^2
+\|(\Phi_x,\Phi_{xt},n_2)\|^2\Big]ds\\
\di+C(\chi+\d)\int_0^t\int\int \f{\nu(|\xi|)|(\widetilde{\mb{G}},\mb{P}_c F_2)|^2}{\mb{M}_\star } d\xi dxds.
\end{array}
\end{equation}

\end{proposition}

\textbf{Proof:} The proof of the higher order estimates in Proposition \ref{Prop3.2} consists the following steps.

\

\underline{Step 1. Estimation on $\|(\phi_x,\psi_x,\z_x)(\cdot, t)\|^2$.}

\

Applying $\partial_x$ to the system \eqref{sys-h}, yields that
\begin{equation}\label{sys-h-h}
\left\{
\begin{array}{ll}
\di \phi_{xt}-\psi_{1xx}=0,\\[2mm]
\di \psi_{1xt}-\f{p}{v}\phi_{xx}+\f{2\z_{xx}}{3v}
=\f43\Big(\frac{\mu(\t)}{v}u_{1x}-\frac{\mu(\bar\t)}{\bar v} \bar u_{1x}\Big)_{xx}-\int \xi_1^2\Pi_{1xx} d\xi-\bar R_{1x},\\[2mm]
\di\quad
-\Big[\f23\z_x\Big(\f1v\Big)_x-\f23\Big(\f{\t}{v^2}\Big)_x\phi_x+\frac23\bar\t_{xx}\Big(\frac{1}{v}-\frac{1}{\bar v}\Big)
+\frac23\bar\t_x\Big(\frac{1}{v}-\frac{1}{\bar v}\Big)_x\\[2mm]
\di\quad-\Big(\frac{p}{v}-\frac{\bar p}{\bar v}\Big)\bar v_{xx}
-\Big(\frac{p}{v}-\frac{\bar p}{\bar v}\Big)_x\bar v_{x}\Big]+(\Phi_x n_2)_x,\\[2mm]
\di \psi_{ixt}=\Big(\frac{\mu(\t)}{v}u_{ix}\Big)_{xx}-\int \xi_1\xi_i\Pi_{1xx} d\xi, ~i=2,3,\\[2mm]
\di \zeta_{xt}+p\psi_{1xx}=-(p_x\psi_{1x}+(p-\bar p)_x\bar u_{1x}+(p-\bar p)\bar u_{1xx})+\Big(\Phi_x\int\xi_1\mb{P}_c F_2d\xi\Big)_x\\[2mm]
\di\quad+\Big(\frac{\k(\t)}{v}\t_{x}-\frac{\k(\bar\t)}{\bar v} \bar \t_{x}\Big)_{xx}
+\f43\Big(\frac{\mu(\t)}{v}u^2_{1x}-\frac{\mu(\bar\t)}{\bar v} \bar u^2_{1x}\Big)_x
+\sum_{i=2}^3\Big(\f{\mu(\t)}{v}u_{ix}^2\Big)_x\\[2mm]
\di\quad-\int \f12 \xi_1|\xi|^2\Pi_{1xx}d\xi+\sum_{i=1}^3 \Big(u_i\int \xi_1\xi_i\Pi_{1x}d\xi\Big)_x-\bar R_{2x}.
\end{array} \right.
\end{equation}
Multiplying $\eqref{sys-h-h}_1$ by $p\phi_x$, $\eqref{sys-h-h}_2$ by
$v\psi_{1x}$, $\eqref{sys-h-h}_3$ by $\psi_{ix}$, $\eqref{sys-h-h}_4$
by $\f{v}{\t}\z_x$, integrating the summation of the resulting equality over $x,t$, and recalling \eqref{pix} and estimating the terms on the right hand side as in the previous section, one has
\begin{equation}\label{5.3}
\begin{array}{ll}
\di \|(\phi_x,\psi_x,\zeta_x)(\cdot,t)\|^2+\int_0^t\|(\psi_{xx},\zeta_{xx})\|^2ds
\leq C\|(\phi_{0x},\psi_{0x},\zeta_{0x})\|^2+C\d^{\frac12} \\
\di +C(\chi+\d)\int_0^t\|(\phi_x,\psi_x,\z_x,\Phi_x,n_2)\|^2ds +C\sum_{|\alpha|=2}\int_0^t\int\int\frac{\nu(|
\xi|)}{\mb{M}_\star }|\partial^\alpha \mathbf{G}|^2d \xi dxds\\
\di+C(\chi+\delta)\sum_{|\alpha^\prime|= 1}\int_0^t\int\int\frac{\nu(|
\xi|)}{\mb{M}_\star }|\partial^{\alpha^\prime}\mb{G}|^2d \xi dxds
+C(\chi+\delta)\int_0^t\int\int\f{\nu(|\xi|)}{\mb{M}_\star }|\pa_{\xi_1}(\mb{P}_c F_2)_x|^2 d\xi dxds\\
\di +C(\chi+\d)\int_0^t\int\int\f{\nu(|\xi|)}{\mb{M}_\star }|(\widetilde{\mb{G}},\mb{P}_c F_2, \pa_{\xi_1}\mb{P}_c F_2)|^2 d\xi dxds.
\end{array}
\end{equation}
To get the estimate of $\phi_{xx}$, applying $\pa_x$ to the equation $\eqref{sys-h-o}_2$:
\begin{equation}\label{jj12}
\begin{array}{l}
\di \psi_{1xt}+(p-\bar p)_{xx}-(\Phi_{xx}n_2+\Phi_x n_{2x})=-\int \xi^2_1\mb{G}_{xx}d\xi- U^{cd}_{1xt}-(\bar p-p_+-p_-)_{xx},
\end{array}
\end{equation}
where
\begin{equation}
\begin{array}{ll}
\di (p-\bar p)_{xx}=\Big(\f{2\z_x}{3v}-\f{2\t\phi_x}{3v^2}+\frac23\bar\t_x\big(\frac{1}{v}
-\frac{1}{\bar v}\big)-\bar v_x\big(\frac{p}{v}-\frac{\bar p}{\bar v}\big)\Big)_x\\[2mm]
\di ~\qquad\qquad =\f{2\z_{xx}}{3v}-\f{2\t\phi_{xx}}{3v^2}+\f23\z_x\big(\f1v\big)_x-\f23\big(\f{\t}{v^2}\big)_x\phi_x
+\frac23\bar\t_{xx}\big(\frac{1}{v}-\frac{1}{\bar v}\big) \\[2mm]
\di ~~\quad\qquad\qquad +\frac23\bar\t_x\big(\frac{1}{v}-\frac{1}{\bar v}\big)_x
-\bar v_{xx}\big(\frac{p}{v}-\frac{\bar p}{\bar v}\big)
-\bar v_{x}\big(\frac{p}{v}-\frac{\bar p}{\bar v}\big)_x.
\end{array}
\end{equation}
Multiplying the above equation by $-\phi_{xx}$ and noting
that
$$
\begin{array}{l}
-\psi_{1xt}\phi_{xx}=-(\psi_{1x}\phi_{xx})_t+\psi_{1x}\phi_{xxt}=-(\psi_{1x}\phi_{xx})_t+\psi_{1x}\psi_{1xxx}\\
\di~~~\quad\quad\quad\quad =-(\psi_{1x}\phi_{xx})_t+(\psi_{1x}\psi_{1xx})_x-\psi_{1xx}^2,
\end{array}
$$
we have:
\begin{equation}\label{cd1}
\begin{array}{l}
\displaystyle \int\psi_{1x}\phi_{xx}(x,0)dx-\int\phi_{xx}\psi_{1x}(x,t)dx+\int_0^t\int\f{2\t}{3v^2}\phi_{xx}^2dxds\leq C\int_0^t\|(\psi_{1xx},\zeta_{xx})\|^2ds\\
\di+C(\chi+\delta)\int_0^t\|(\phi_x,\zeta_x,\Phi_x, n_2,n_{2x})\|^2 ds
+C\sum_{|\alpha|=2}\int_0^t\int\int\frac{\nu(| \xi|)}{\mb{M}_\star }|\partial^\alpha \mathbf{G}|^2d \xi dxds+C\d^{\frac12}.
\end{array}
\end{equation}
We also use the original fluid-type equation \eqref{sys-h-o} to estimate
$\|(\phi_{xt},\psi_{xt},\zeta_{xt})\|$ and $\|(\phi_{tt},\psi_{tt},\zeta_{tt})\|$.
For example, we only estimate $\|\psi_{1xt}\|^2$ here. Multiplying \eqref{jj12} by $\psi_{1xt}$ and integrating over $x$ and $t$, we
have
\begin{equation}\label{cd2}
\begin{array}{l}
\displaystyle \int_0^t\|\psi_{1xt}\|^2 ds\leq C\int_0^t\|(\phi_{xx},\zeta_{xx})\|^2 ds
+C(\chi+\d)\int_0^t\|(\phi_x,\z_x,\Phi_x,n_2,n_{2x})\|^2ds \\
\di \qquad\qquad\qquad+C\sum_{|\alpha|=2}\int_0^t\int\int\frac{\nu(| \xi|)}{\mb{M}_\star }|\partial^\alpha
\mathbf{G}|^2d \xi dxds + C\d^{\frac12}.
\end{array}
\end{equation}
Applying $\partial_t $ the equation $\eqref{sys-h-o}_2$ gives that
\begin{equation*}\label{jj}
\begin{array}{l}
\di \psi_{1tt}+\f{2\z_{xt}}{3v}-\f{2\t\phi_{xt}}{3v^2}+\f23\z_x\big(\f1v\big)_t-\f23\big(\f{\t}{v^2}\big)_t\phi_x
+\frac23\bar\t_{xt}\big(\frac{1}{v}-\frac{1}{\bar v}\big)+\frac23\bar\t_{x}\big(\frac{1}{v}-\frac{1}{\bar v}\big)_t\\
\di\quad -\bar v_{xt}\big(\frac{p}{v}-\frac{\bar p}{\bar v}\big)-\bar v_{x}\big(\frac{p}{v}-\frac{\bar p}{\bar v}\big)_t
-(\Phi_{xt}n_2+\Phi_x n_{2t})=-\int\xi_1^2\mb{G}_{xt}d\xi- U^{cd}_{1tt}-(\bar p-p_+-p_-)_{xt},
\end{array}
\end{equation*}
which implies that
\begin{equation}\label{cd3}
\begin{array}{l}
\displaystyle \int_0^t\|\psi_{1tt}\|^2 ds\leq C\int_0^t\|(\phi_{xt},\zeta_{xt})\|^2 ds
+C(\chi+\delta)\sum_{|\alpha'|=1}\int_0^t\|(\pa^{\alpha'}(\phi,\z),n_2,n_{2t})\|^2ds\\
\di \qquad\qquad\qquad+C\sum_{|\alpha|=2}\int_0^t\int\int\frac{\nu(| \xi|)}{\mb{M}_\star }|\partial^\alpha
\mathbf{G}|^2d \xi dxds+C\d^{\frac12}.
\end{array}
\end{equation}
The combinations of \eqref{5.3}, \eqref{cd1}-\eqref{cd3} gives that
\begin{equation}\label{55e}
\begin{array}{l}
\displaystyle
\|(\phi_x,\psi_x,\zeta_x)(\cdot,t)\|^2+\sum_{|\alpha|=2}\int_0^t\|\partial^\alpha(\phi,\psi,\zeta)\|^2ds\leq
 C\|(\phi_{0x},\psi_{0x},\zeta_{0x},\phi_{0xx})\|^2+C\d^{\frac12}\\
 \di +C\int|\phi_{xx}\psi_{1x}|(x,t)dx+C(\chi+\delta)\int_0^t\Big[\|(\Phi_{x},n_2)\|^2
 +\sum_{|\alpha^\prime|=1}\|\partial^{\alpha^\prime}(\phi,\psi,\zeta,n_2)\|^2\Big]ds\\
 \di +C(\chi+\d)\sum_{|\alpha'|=1}\int_0^t\int\int\f{\nu(|\xi|)}{\mb{M}_\star }|(\widetilde{\mb{G}},\mb{P}_c F_2,\partial^{\alpha'}(\mb{G},\mb{P}_c F_2))|^2 d\xi dxds\\
\di +C(\chi+\d)\int_0^t\int\int\f{\nu(|\xi|)}{\mb{M}_\star }|(\pa_{\xi_1}(\mb{P}_c F_2),\pa_{\xi_1}(\mb{P}_c F_2)_x)|^2 d\xi dxds
+C\sum_{|\alpha|=2}\int_0^t\int\int\frac{\nu(|\xi|)}{\mb{M}_\star }|\partial^\alpha
\mathbf{G}|^2d \xi dxds.
\end{array}
\end{equation}
Now we estimate $\|n_{2xx}\|^2.$ Applying $\partial_x$ to the equation \eqref{n2} leads to
\begin{equation}\label{n2x}
\begin{array}{ll}
\di n_{2xt}+(\f{u_1x}{v}n_2)_{x}+ \Big(\f1v\big(\f{3\k_1(\t)}{2\t v}\Phi_x\big)_x\Big)_{x}- \Big(\f1v\big(\f{\k_1(\t)}{v}(n_2v)_x\big)_x\Big)_{x}\\
\di
 =-\Big(\f1v\big(n_2\int \xi_1 \mb{N}_{\mb{M}}^{-1}\big[\mb{P}_c(\xi_1\mb{M}_x)\big]d\xi)_x\Big)_{x}
 -\Big(\f1v\big(\int \xi_1 \mb{N}_{\mb{M}}^{-1}\big[\f1v \mb{P}_c(\xi_1 (\mb{P}_cF_2)_x)\big]d\xi\big)_x\Big)_{x}\\
\di\quad
 - \Big(\f1v\big(\int \xi_1 \mb{N}_{\mb{M}}^{-1}\big[\f{\Phi_x}{v}\mb{G}_{\xi_1}\big]d\xi\big)_x\Big)_{x}
 -\Big(\f1v\big(\int \xi_1 \mb{N}_{\mb{M}}^{-1}\Big[(\mb{P}_c F_2)_t-\f{u_1}{v}(\mb{P}_c F_2)_x\\
\di\quad +\big((\mb{M}v)_t-\f{u_1}{v}(\mb{M}v)_x\big)_t~n_2
 -2Q(F_2,\mb{G})\Big]d\xi\big)_x\Big)_{x}.
\end{array}
\end{equation}
Multiplying the equation \eqref{n2x} by $n_{2x}$ gives that
\begin{equation}\label{n2x-1}
\begin{array}{ll}
\di \big(\f{n_{2x}^2}{2}\big)_t+\f{3\k_1(\t)}{\t}n_{2x}^2+\f{\k_1(\t)}{v}n_{2xx}^2=(\cdots)_x
-\big(\f{u_{1x}}{v}\big)_xn_2 n_{2x}-\f{u_{1x}}{v}n_{2x}^2\\
\di +\f{3}{2v^2}\big(\f{\k_1(\t)}{v}\big)_x\Phi_x n_{2xx}
-\big(\f{3\k_1(\t)}{\t}\big)_x n_2 n_{2x}-\f{\k_1(\t)}{v^2}(2n_{2x}v_x+n_2v_{xx})n_{2xx}\\
\di -\big(\f{\k_1(\t)}{v}\big)_x\big(n_{2x}+n_2\f{v_x}{v}\big)n_{2xx}
+\f{n_{2xx}}{v}\Big(n_2\int \xi_1 \mb{N}_{\mb{M}}^{-1}\big[\mb{P}_c(\xi_1\mb{M}_x) \big]d\xi \Big)_x\\
\di +\f{n_{2xx}}{v}\Big(\int \xi_1 \mb{N}_{\mb{M}}^{-1}\big[\f1v\mb{P}_c(\xi_1 (\mb{P}_cF_2)_x)\big]d\xi\Big)_x
+\f{n_{2xx}}{v}\Big(\int \xi_1 \mb{N}_{\mb{M}}^{-1}\big[\f{\Phi_x}{v}\mb{G}_{\xi_1}\big]d\xi\Big)_x\\
\di+\f{n_{2xx}}{v}\Big(\int \xi_1 \mb{N}_{\mb{M}}^{-1}\Big[(\mb{P}_c F_2)_t-\f{u_1}{v}(\mb{P}_c F_2)_x
+\big((\mb{M}v)_t-\f{u_1}{v}(\mb{M}v)_x\big)~n_2-2Q(F_2,\mb{G})\Big]d\xi\Big)_x.
\end{array}
\end{equation}
Integrating the above equation with respect to $x,t$ yields that
\begin{equation}\label{n2x-2}
\begin{array}{ll}
\di \|n_{2x}(\cdot,t)\|^2+\int_0^t\|(n_{2x},n_{2xx})\|^2ds\leq C\|n_{20x}\|^2
+C(\chi+\d)\int_0^t\Big[\|n_2\|^2+\sum_{1\leq |\alpha|\leq 2}\|\pa^{\alpha}(\phi,\psi,\z)\|^2\Big]ds\\
\di \quad +C(\chi+\d)\int_0^t\Big[\sum_{|\alpha'|=1}\int\int\f{\nu(|\xi|)}{\mb{M}_\star }|\pa^{\alpha'}(\mb{G},\mb{P}_c F_2)|^2
+\sum_{0\leq|\beta'|\leq1}\int\int\f{\nu(|\xi|)}{\mb{M}_\star }|\pa^{\beta'}(\widetilde{\mb{G}},\mb{P}_c F_2)|^2\Big]ds\\
\di\quad +C(\chi+\d)\int_0^t\int\int\f{\nu(|\xi|)}{\mb{M}_\star }|\pa_{\xi_1}\mb{G}_x|^2d\xi dxds
+C\sum_{|\alpha|=2}\int_0^t\int\int\f{\nu(|\xi|)}{\mb{M}_\star }|\pa^{\alpha}(\mb{P}_c F_2)|^2d\xi dxds+C\d^{\frac12}.
\end{array}
\end{equation}

By the equation \eqref{n22}, it holds that
\begin{equation}\label{n2xt}
\begin{array}{ll}
\di \int_0^t\|n_{2xt}\|^2ds\leq C(\chi+\d)\int_0^t\Big[\|(\psi_{1xx},n_2,n_{2x})\|^2
+\int\int\f{\nu(|\xi|)}{\mb{M}_\star }|(\mb{P}_cF_2)_x|^2d\xi dx\Big]ds\\
\di\qquad\qquad\qquad \leq C\int_0^t\int\int\f{\nu(|\xi|)|(\mb{P}_cF_2)_{xx}|^2}{\mb{M}_\star }d\xi dxds,
\end{array}
\end{equation}
and
\begin{equation}\label{n2tt}
\begin{array}{ll}
\di \int_0^t\|n_{2tt}\|^2ds\leq C(\chi+\d)\int_0^t\Big[\|(\psi_{1xt},n_2,n_{2t})\|^2
+\int\int\f{\nu(|\xi|)}{\mb{M}_\star }|(\mb{P}_cF_2)_x|^2d\xi dx\Big]ds\\
\di\qquad\qquad\qquad \leq C\int_0^t\int\int\f{\nu(|\xi|)|(\mb{P}_cF_2)_{xt}|^2}{\mb{M}_\star }d\xi dxds.
\end{array}
\end{equation}
By \eqref{n2x-2}, \eqref{n2xt} and \eqref{n2tt}, one has
\begin{equation}\label{n2-he}
\begin{array}{ll}
\di \|n_{2x}(\cdot,t)\|^2+\int_0^t\Big[\|n_{2x}\|^2+\sum_{|\alpha|=2}\|\partial^\alpha n_2\|^2\Big] ds\leq C\|n_{20x}\|^2+C\d^{\frac12}\\
\di\quad+C(\chi+\delta)\int_0^t\Big[\|(n_2,n_{2t})\|^2+\sum_{1\leq|\alpha|\leq2}\|\partial^\alpha (\phi,\psi,\zeta)\|^2
+\int\int\f{\nu(|\xi|)}{\mb{M}_\star }|\pa_{\xi_1}\mb{G}_x|^2d\xi dx\Big]ds\\
\di\quad+C(\chi+\d)\sum_{|\alpha'|=1}\int_0^t\int\int \f{\nu(|\xi|)}{\mb{M}_\star }|\partial^{\alpha'}(\mb{G},\mb{P}_c F_2)|^2d\xi dxds\\
\di\quad+C(\chi+\d)\sum_{0\leq|\beta'|\leq1}\int_0^t\int\int \f{\nu(|\xi|)}{\mb{M}_\star }|\partial^{\beta'}(\widetilde{\mb{G}},\mb{P}_c F_2)|^2d\xi dxds\\
\di\quad +C\sum_{|\alpha|=2}\int_0^t\int\int \f{\nu(|\xi|)}{\mb{M}_\star }|\partial^\alpha(\mb{P}_c F_2)|^2d\xi dxds.
\end{array}
\end{equation}

\

\underline{Step 2. Estimation on the non-fluid part.}

\

To close the above estimate, we need to estimate the derivatives
on the non-fluid component $|\pa^{\beta}(\widetilde{\mb{G}},\mb{P}_cF_2)|(1\leq|\beta|\leq2)$ and
$|\partial^{\alpha'}\partial^{\b'}(\mathbf{G}, \mb{P}_cF_2)|(|\alpha'|=1,|\beta'|=1)$.
First, applying $\partial_x$ to \eqref{F1-G} and \eqref{F2-pc-L}, respectively, we have
\begin{equation}\label{Gx}
\begin{array}{ll}
\di \mathbf{G}_{xt}-\mb{L}_\mathbf{M}\mathbf{G}_x=\big(\f{u_1}{v}\mb{G}_x\big)_x-\big(\f1v\mb{P}_1(\xi_1\mb{M}_x)\big)_x
-\big(\f1v\mb{P}_1(\xi_1\mb{G}_x)\big)_x-\big(\f{\Phi_x}{v}\mb{P}_1(\pa_{\xi_1}F_2)\big)_x\\[2mm]
\di ~~\qquad\qquad\qquad+2Q(\mathbf{M}_x,\mathbf{G})+2Q(\mathbf{G},\mathbf{M}_x)+2Q(\mathbf{G}_x,
\mathbf{G})+2Q(\mathbf{G},
\mathbf{G}_x),
\end{array}
\end{equation}
and
\begin{equation}\label{Pc-F2x}
\begin{array}{ll}
\di
(\mb{P}_c F_2)_{xt}-\mb{N}_{\mb{M}}(\mb{P}_c F_2)_x=-\big(\f{\xi_1}{v}F_{2x}\big)_x+\big(\f{u_1}{v}F_{2x}\big)_x
-(n_2\mb{M}v)_{xt}-\big(\f{\Phi_x}{v}\mb{P}_c(\pa_{\xi_1}F_1)\big)_x\\[2mm]
\di\quad\qquad\qquad\qquad\qquad\qquad+2Q(\mb{P}_c F_2,\mathbf{M}_x)+2Q(F_{2x},\mb{G})+2Q(F_{2},\mb{G}_x).
\end{array}
\end{equation}
Multiplying \eqref{Gx} and \eqref{Pc-F2x} by
$\frac{\mathbf{G}_x}{\mb{M}_\star }$ and
$\frac{(\mathbf{P}_cF_2)_x}{\mb{M}_\star }$, respectively, and
integrating with respect to $x, \xi$ and $ t$, and noting that
\begin{equation}\label{wt5.1}
\begin{array}{l}
\di-\int_0^t\int\int\Big(\f{\Phi_x}{v}\mb{P}_c(\pa_{\xi_1}F_1)\Big)_x\f{(\mb{P}_cF_2)_x}{\mb{M}_\star }d\xi dxds
=\int_0^t\int\int\f{\Phi_x}{v}\pa_{\xi_1}F_1\f{(\mb{P}_cF_2)_{xx}}{\mb{M}_\star }d\xi dxds\\
\di=-\int_0^t\int\int\f{\Phi_x}{v}\f{\xi_1-u_1}{R\t}\mb{M}\f{(\mb{P}_cF_2)_{xx}}{\mb{M}_\star }d\xi dxds
+\int_0^t\int\int\f{\Phi_x}{v}(\widetilde{\mb{G}}_{\xi_1}+\bar{\mb{G}}_{\xi_1})\f{(\mb{P}_cF_2)_{xx}}{\mb{M}_\star }d\xi dxds\\
\di=-\int_0^t\int\int\f{\Phi_x}{v}\f{\xi_1-u_1}{R\t}\big(\f{\mb{M}}{\mb{M}_\star }-1\big)(\mb{P}_cF_2)_{xx}d\xi dxds
+\int_0^t\int\int\f{\Phi_x}{v}(\widetilde{\mb{G}}_{\xi_1}+\bar{\mb{G}}_{\xi_1})\f{(\mb{P}_cF_2)_{xx}}{\mb{M}_\star }d\xi dxds\\
\di\quad-\int_0^t\int\int\f{\Phi_x}{v}\f{\xi_1-u_1}{R\t}(\mb{P}_cF_2)_{xx}d\xi dxds\\
\di=-\int_0^t\int\int\f{\Phi_x}{v}\f{\xi_1-u_1}{R\t}\big(\f{\mb{M}}{\mb{M}_\star }-1\big)(\mb{P}_cF_2)_{xx}d\xi dxds
+\int_0^t\int\int\f{\Phi_x}{v}(\widetilde{\mb{G}}_{\xi_1}+\bar{\mb{G}}_{\xi_1})\f{(\mb{P}_cF_2)_{xx}}{\mb{M}_\star }d\xi dxds\\
\di\quad+\int_0^t\int\int\f{\Phi_x}{v}\big(\f{\xi_1-u_1}{R\t}\big)_x(\mb{P}_cF_2)_{x}d\xi dxds
+\int_0^t\int\int\big(\f{\Phi_x}{v}\big)_x\f{\xi_1-u_1}{R\t}(\mb{P}_cF_2)_{x}d\xi dxds,
\end{array}
\end{equation}
where by \eqref{n22}, one has
$$
\begin{array}{ll}
\di\quad\int_0^t\int\int\big(\f{\Phi_x}{v}\big)_x\f{\xi_1-u_1}{R\t}(\mb{P}_cF_2)_{x}d\xi dxds\\
\di=\int_0^t\int\big(\f{\Phi_x}{v}\big)_x\f{1}{R\t}\Big(\int\xi_1\mb{P}_c F_2d\xi\Big)_x dxds
=-\int_0^t\int\f{2n_2v}{R\t}(n_{2t}v+u_{1x}n_2) dxds\\[2mm]
\di=\int_0^t\int\Big(-(\f{n_2^2v^2}{R\t})_t+n_2^2\big(\f{v^2}{R\t}\big)_t-2n_2^2\f{u_{1x}v }{R\t}\Big)dxds\\
\di=-\int\f{n_2^2v^2}{R\t}(x,t)dx+\int\f{n_2^2v^2}{R\t}(x,0)dx-\int_0^t\int\f{n_2^2v^2\t_t}{R\t^2}dxds,
\end{array}
$$
which yields
\begin{equation}\label{wt5.2}
\begin{array}{l}
\di-\int_0^t\int\int\Big(\f{\Phi_x}{v}\mb{P}_c(\pa_{\xi_1}F_1)\Big)_x\f{(\mb{P}_cF_2)_x}{\mb{M}_\star }d\xi dxds
=-\int\f{n_2^2v^2}{R\t}(x,t)dx+\int\f{n_2^2v^2}{R\t}(x,0)dx\\
\di-\int_0^t\int\int\f{\Phi_x}{v}\f{\xi_1-u_1}{R\t}\big(\f{\mb{M}}{\mb{M}_\star }-1\big)(\mb{P}_cF_2)_{xx}d\xi dxds
+\int_0^t\int\int\f{\Phi_x}{v}(\widetilde{\mb{G}}_{\xi_1}+\bar{\mb{G}}_{\xi_1})\f{(\mb{P}_cF_2)_{xx}}{\mb{M}_\star }d\xi dxds\\
\di+\int_0^t\int\int\f{\Phi_x}{v}\big(\f{\xi_1-u_1}{R\t}\big)_x(\mb{P}_cF_2)_{x}d\xi dxds
-\int_0^t\int\f{n_2^2v^2\t_t}{R\t^2}dxds\\
\di\leq -\int\f{n_2^2v^2}{R\t}(x,t)dx+\int\f{n_2^2v^2}{R\t}(x,0)dx+C\|n_{20}\|^2+C(\chi+\d+\eta_0)\int_0^t\|(\Phi_x, n_2)\|^2ds\\
\di \quad +C(\chi+\d)\int_0^t\int\int\f{\nu(|\xi|)}{\mb{M}_\star }|(\widetilde{\mb{G}}_{\xi_1},(\mb{P}_cF_2)_x)|^2d\xi dxds\\
\di\quad+C(\chi+\d+\eta_0)\int_0^t\int\int\f{\nu(|\xi|)}{\mb{M}_\star }|(\mb{P}_cF_2)_{xx}|^2d\xi dxds.
\end{array}
\end{equation}
It follows from \eqref{wt5.2} and Lemmas \ref{Lemma 4.1}-\ref{Lemma 4.3} that
\begin{equation}\label{Gx-E}
\begin{array}{l}
\displaystyle
\int\int\frac{|(\mathbf{G}, \mb{P}_cF_2)_x|^2}{\mb{M}_\star }(x,\xi,t) d \xi dx+\|n_{2}(\cdot,t)\|^2
+\int_0^t\int\int\frac{\nu(| \xi|)}{\mb{M}_\star }|(\mathbf{G}, \mb{P}_cF_2)_x|^2d \xi dxds\\
\leq\displaystyle C\int\int\frac{|(\mathbf{G},
\mb{P}_cF_2)_x|^2}{\mb{M}_\star }(x,\xi,0) d \xi
dx+C\|n_{20}\|^2+C\int_0^t\int\int\frac{\nu(| \xi|)}{\mb{M}_\star }|(\mathbf{G}, \mb{P}_cF_2)_{xx}|^2d \xi dxds\\
\di +C\sum_{|\alpha|=2}\int_0^t\|\pa^{\alpha}(\phi,\psi,\z,n_2)\|^2 ds
+C(\chi+\delta)\sum_{|\alpha^\prime|=1}\int_0^t
\|\partial^{\alpha^\prime}(\phi,\psi,\zeta,n_2)\|^2ds+C\d^{\frac12} \\
\di +C(\chi+\delta+\eta_0)\int_0^t \|(\Phi_x,n_2)\|^2ds
+C(\chi+\delta)\sum_{0\leq|\beta^\prime|\leq1}\int_0^t\int\int\frac{\nu(|
\xi|)}{\mb{M}_\star }|\partial^{\beta^\prime}(\widetilde{\mb{G}},\mb{P}_cF_2)|^2d \xi dxds.
\end{array}
\end{equation}
Similar to \eqref{wt5.1}, one has
$$
\begin{array}{ll}
\di-\int_0^t\int\int\Big(\f{\Phi_x}{v}\mb{P}_c(\pa_{\xi_1}F_1)\Big)_t\f{(\mb{P}_cF_2)_t}{\mb{M}_\star }d\xi dxds
=-\int_0^t\int\int\Big(\f{\Phi_x}{v}\pa_{\xi_1}F_1\Big)_t\f{(\mb{P}_cF_2)_{t}}{\mb{M}_\star }d\xi dxds\\
\di=\int_0^t\int\int\f{\Phi_{xt}}{v}\big(\f{\xi_1-u_1}{R\t}\mb{M}-\widetilde{\mb{G}}_{\xi_1}-\bar{\mb{G}}_{\xi_1}\big)
\f{(\mb{P}_cF_2)_{t}}{\mb{M}_\star }d\xi dxds\\
\di\quad-\int_0^t\int\int\f{\Phi_{x}v_t}{v^2}\big(\f{\xi_1-u_1}{R\t}\mb{M}-\widetilde{\mb{G}}_{\xi_1}-\bar{\mb{G}}_{\xi_1}\big)
\f{(\mb{P}_cF_2)_{t}}{\mb{M}_\star }d\xi dxds\\
\di\quad+\int_0^t\int\int\f{\Phi_{x}}{v}\big(\f{\xi_1-u_1}{R\t}\mb{M}-\mb{G}_{\xi_1}\big)_t\f{(\mb{P}_cF_2)_{t}}{\mb{M}_\star }d\xi dxds\\
\di=\int_0^t\int\int\f{\Phi_{xt}}{v}\f{\xi_1-u_1}{R\t}\big(\f{\mb{M}}{\mb{M}_\star }-1\big)(\mb{P}_cF_2)_{t}d\xi dxds
-\int_0^t\int\int\f{\Phi_{xt}}{v}(\widetilde{\mb{G}}_{\xi_1}+\bar{\mb{G}}_{\xi_1})\f{(\mb{P}_cF_2)_{t}}{\mb{M}_\star }d\xi dxds\\
\di\di\quad-\int_0^t\int\int\f{\Phi_{x}v_t}{v^2}\big(\f{\xi_1-u_1}{R\t}\mb{M}-\widetilde{\mb{G}}_{\xi_1}-\bar{\mb{G}}_{\xi_1}\big)
\f{(\mb{P}_cF_2)_{t}}{\mb{M}_\star }d\xi dxds\\
\di\quad+\int_0^t\int\int\f{\Phi_{x}}{v}\big(\f{\xi_1-u_1}{R\t}\mb{M}-\mb{G}_{\xi_1}\big)_t\f{(\mb{P}_cF_2)_{t}}{\mb{M}_\star }d\xi dxds
+\int_0^t\int\int\f{\Phi_{xt}}{v}\f{\xi_1-u_1}{R\t}(\mb{P}_cF_2)_{t}d\xi dxds,
\end{array}
$$
where by \eqref{wt-3.31}, one has
$$
\begin{array}{ll}
\di\quad \int_0^t\int \f{\Phi_{xt}}{v}\f{\xi_1-u_1}{R\t}(\mb{P}_c F_2)_tdxds\\
\di =\int_0^t\int \f{\Phi_{xt}}{R\t v}\Big(\int \xi_1\mb{P}_c F_2d\xi\Big)_tdxds
=-\int_0^t\int\f{\Phi_{xt}}{R\t v}\big(\f{\Phi_x}{2v}\big)_{tt}dxds\\
\di=-\int_0^t\int\f{\Phi_{xt}}{2R\t v}\Big(\f{\Phi_{xtt}}{v}+2\Phi_{xt}\big(\f1v\big)_t+\Phi_x\big(\f1v\big)_{tt}\Big)dxds\\
\di=-\int\f{\Phi_{xt}^2}{4R\t v^2}(x,t)dx+\int\f{\Phi_{xt}^2}{4R\t v^2}(x,0)dx\\
\di\quad+\int_0^t\int \Big(\f{\Phi_{xt}^2}{4R}\big(\f{1}{\t v^2}\big)_t-\f{\Phi_{xt}^2}{R\t v}\big(\f1v\big)_t
-\f{\Phi_{xt}\Phi_x}{2R\t v}\big(\f1v\big)_{tt}\Big) dxds,
\end{array}
$$
which leads to
\begin{equation}\label{wt5.3}
\begin{array}{l}
\di-\int_0^t\int\int\Big(\f{\Phi_x}{v}\mb{P}_c(\pa_{\xi_1}F_1)\Big)_t\f{(\mb{P}_cF_2)_t}{\mb{M}_\star }d\xi dxds\\
\di\leq-\int\f{\Phi_{xt}^2}{4R\t v^2}(x,t)dx+C\|\Phi_{xt}(\cdot,0)\|^2+C(\chi+\d)\int_0^t\|(\Phi_x,\phi_{tt})\|^2ds\\
\di\quad+C(\chi+\d+\eta_0)\int_0^t\|\Phi_{xt}\|^2ds
+C(\chi+\d)\int_0^t\int\int\f{\nu(|\xi|)}{\mb{M}_\star }|(\widetilde{\mb{G}}_{\xi_1},\mb{G}_{\xi_{1}t})|^2d\xi dxds\\
\di\quad+C(\chi+\d+\eta_0)\int_0^t\int\int\f{\nu(|\xi|)}{\mb{M}_\star }|(\mb{P}_c F_2)_t|^2d\xi dxds.
\end{array}
\end{equation}
By \eqref{wt5.3}, we further have
\begin{equation}\label{Gt-E}
\begin{array}{l}
\displaystyle
\int\int\frac{|(\mathbf{G}, \mb{P}_cF_2)_t|^2}{2\mb{M}_\star }(x,\xi,t)d \xi dx+\|\Phi_{xt}(\cdot,t)\|^2
+\int_0^t\int\int\frac{\nu(| \xi|)}{\mb{M}_\star }|(\mathbf{G}, \mb{P}_cF_2)_t|^2d \xi dxds\\
\leq\displaystyle C\int\int\frac{|(\mathbf{G},
\mb{P}_cF_2)_t|^2}{2\mb{M}_\star }(x,\xi,0)d \xi
dx+C\|\Phi_{0xt}\|^2+C\sum_{|\alpha|=2}\int_0^t\|\pa^{\alpha}(\phi,\psi,\z,n_2)\|^2ds+C\d^{\frac12}\\
\di +C(\chi+\delta)\int_0^t\Big[\sum_{|\alpha^\prime|=1}\|\partial^{\alpha^\prime}(\phi,\psi,\zeta,n_2)\|^2+\|(\Phi_x,n_2)\|^2\Big]ds
+C(\chi+\delta+\eta_0)\int_0^t\|\Phi_{xt}\|^2ds\\
\di+C(\chi+\delta)\int_0^t\Big[\sum_{0\leq|\beta^\prime|\leq1}\int\int\frac{\nu(|\xi|)}{\mb{M}_\star }
|\partial^{\beta^\prime}(\widetilde{\mb{G}}, \mb{P}_cF_2)|^2d\xi dx
+\int\int\frac{\nu(|\xi|)}{\mb{M}_\star }|(\mathbf{G},\mb{P}_cF_2)_x|^2 d\xi dx\Big]ds\\
\di+C(\chi+\delta)\int_0^t\int\int\frac{\nu(| \xi|)}{\mb{M}_\star }|\pa_{\xi_1}(\mb{G},\mb{P}_c F_2)_t|^2d\xi dxds
+C\int_0^t\int\int\frac{\nu(| \xi|)}{\mb{M}_\star }|(\mathbf{G},\mb{P}_cF_2)_{xt}|^2d \xi
dxds.
\end{array}
\end{equation}
The combination of \eqref{55e}, \eqref{n2-he}, \eqref{Gx-E} and \eqref{Gt-E}, and choosing $\chi,\delta$ suitably small  yield that
\begin{equation}\label{xt-e}
\begin{array}{ll}
\displaystyle
\sum_{|\alpha'|=1}\int\int\frac{|\partial^{\alpha'}(\mathbf{G}, \mb{P}_cF_2)|^2}{\mb{M}_\star }(\cdot,t)d \xi dx+\|(n_2,\Phi_{xt})(\cdot,t)\|^2+\sum_{|\alpha'|=1}\int_0^t\int\int\frac{\nu(| \xi|)}{\mb{M}_\star }|\partial^{\alpha'}(\mathbf{G}, \mb{P}_cF_2)|^2d \xi dxds\\
\leq\displaystyle C\sum_{|\alpha'|=1}\int\int\frac{|\partial^{\alpha'}(\mathbf{G},
\mb{P}_cF_2)|^2}{\mb{M}_\star }(x,\xi,0)d \xi dx+C(\varepsilon_0+\delta^{\frac12})
 +C\int |\phi_{xx}\psi_{1x}|(x,t)dx\\
\di\quad +C(\chi+\d)\sum_{|\alpha'|=1}\int_0^t\pa^{\alpha'}\|(\phi,\psi,\z,n_2)\|^2ds
+C\sum_{|\alpha|=2}\int_0^t\int\int\frac{\nu(| \xi|)}{\mb{M}_\star }|\partial^\alpha(\mathbf{G},
\mb{P}_cF_2)|^2d \xi dxds\\
\di\quad +C(\chi+\delta)\sum_{0\leq|\beta^\prime|\leq1}\int_0^t\int\int\frac{\nu(|
\xi|)}{\mb{M}_\star }|\partial^{\beta^\prime}(\widetilde{\mb{G}}, \mb{P}_cF_2\big)|^2d \xi dxds+C(\chi+\delta+\eta_0)\int_0^t\|(\Phi_x,\Phi_{xt},n_2)\|^2ds\\
\di\quad +C(\chi+\delta)\sum_{|\alpha'|=1,|\beta'|=1}\int_0^t\int\int\frac{\nu(|
\xi|)}{\mb{M}_\star }|\pa^{\alpha'}\pa^{\beta'}(\mb{G}, \mb{P}_cF_2)|^2d \xi dxds.
\end{array}
\end{equation}
Applying $\partial_{\xi_j}~(j=1,2,3)$ to the equation \eqref{GE}  and noting that \eqref{P1-1}, one has
\begin{equation}\label{Gv}
\begin{array}{ll}
\di
\widetilde{\mb{G}}_{\xi_jt}-\mb{L}_{\mb{M}}\widetilde{\mb{G}}_{\xi_j}=
-\f{3}{2\t v}\pa_{\xi_1}\mb{P}_1\Big[\xi_1\big(\f{|\xi-u|^2}{2\t}\zeta_x+\xi\cdot\psi_{x}\big)\mb{M}\Big]
+\f{u_1}{v}\mb{G}_{\xi_jx}-\f1v\partial_{\xi_j}\mb{P}_1(\xi_1\mb{G}_x)\\[3mm]
\di\quad
-\Phi_x n_2\mb{M}_{\xi_j\xi_1} -\Phi_x n_2\sum_{k=0}^4\int\mb{M} \big(\f{\chi_k}{\mb{M}}\big)_{\xi_1}d\xi(\chi_k)_{\xi_j}
-\f{\Phi_x}{v}\big(\mb{P}_cF_2\big)_{\xi_1\xi_j}-\bar{\mb{G}}_{\xi_jt}\\
\di\quad-\f{\Phi_x}{v}\sum_{k=0}^4\int \mb{P}_cF_2 \big(\f{\chi_k}{\mb{M}}\big)_{\xi_1}d\xi(\chi_k)_{\xi_j}
+2\partial_{\xi_j}Q(\mb{G},\mb{G})+2Q(\mb{M}_{\xi_j},\widetilde{\mb{G}})+2Q(\widetilde{\mb{G}},\mb{M}_{\xi_j}).
\end{array}
\end{equation}
Similarly, applying $\partial_{\xi_j}~(j=1,2,3)$ to the equation \eqref{F2-pc-L}  and noting that \eqref{Pc-1}  give
that
\begin{equation}\label{F2-pc-v}
\begin{array}{ll}
\di (\mb{P}_c F_2)_{\xi_jt}-\mb{N}_{\mb{M}}(\mb{P}_c F_2)_{\xi_j}=-\f1v\big(\xi_1(n_2 v\mb{M})_x\big)_{\xi_j}
+\f{u_1}{v}(n_2 v\mb{M})_{\xi_jx}\\[2mm]
\di\quad-\f1v\big(\xi_1(\mb{P}_c F_2)_x\big)_{\xi_j}+\f{u_1}{v}(\mb{P}_c F_2)_{\xi_jx}
-(n_2 v\mb{M}_{\xi_j})_t-\f{\Phi_x}{v}(\mb{M}_{\xi_j\xi_1}+\bar{\mb{G}}_{\xi_j\xi_1})\\[2mm]
\di\quad -\f{\Phi_x}{v}\widetilde{\mb{G}}_{\xi_j\xi_1}+2\pa_{\xi_j}Q(F_2,\mb{G})+2Q(\mb{P}_c F_2,\mb{M}_{\xi_j}).
\end{array}
\end{equation}
Recalling the following fact
\begin{equation}\label{Fact1}
\partial_{\xi_j}\mb{L}_{\mb{M}}g-\mb{L}_{\mb{M}}(\pa_{\xi_j}g)=2Q(\partial_{\xi_j}\mb{M},g)+2Q(g,\partial_{\xi_j}\mb{M}).
\end{equation}
Multiplying the equations \eqref{Gv} and \eqref{F2-pc-v} by
$\f{\widetilde{\mb{G}}_{\xi_j}}{\mb{M}_\star }$ and $\f{(\mb{P}_cF_2)_{\xi_j}}{\mb{M}_\star }$, respectively, and then integrating with
respect to $x,\xi,t$ imply that
\begin{equation}\label{GvE0}
\begin{array}{ll}
\di \int\int \f{|(\widetilde{\mb{G}}, \mb{P}_c F_2)_{\xi_j}|^2}{2\mb{M}_\star }(x,\xi,t)
d\xi dx +\tilde\sigma\int_0^t\int\int\f{\nu(|\xi|)}{\mb{M}_\star }|(\widetilde{\mb{G}}, \mb{P}_c F_2)_{\xi_j}|^2
d\xi dxds\\[3mm]
 \di
\leq \int\int \f{|(\widetilde{\mb{G}}, \mb{P}_c F_2)_{\xi_j}|^2}{2\mb{M}_\star }(x,\xi,0)
d\xi dx+
\Big|\int_0^t\int\int\Big\{-\f{3}{2\t v}\pa_{\xi_1}\mb{P}_1\Big[\xi_1\big(\f{|\xi-u|^2}{2\t}\zeta_x+\xi\cdot\psi_{x}\big)\mb{M}\Big]\\
\di\quad+\Big(\f{u_1}{v}\mb{G}_{\xi_jx}-\f1v\partial_{\xi_j}\mb{P}_1(\xi_1\mb{G}_x)\Big)
-\Phi_x n_2\Big(\mb{M}_{\xi_j\xi_1}+\sum_{k=0}^4\int\mb{M} \Big(\f{\chi_k}{\mb{M}}\Big)_{\xi_1}d\xi(\chi_k)_{\xi_j}\Big) \\
\di\quad-\f{\Phi_x}{v}\sum_{k=0}^4\int \mb{P}_cF_2 \big(\f{\chi_k}{\mb{M}}\big)_{\xi_1}d\xi(\chi_k)_{\xi_j}
+2\partial_{\xi_j}Q(\mb{G},\mb{G})-\bar{\mb{G}}_{\xi_jt}
+2\big(Q(\mb{M}_{\xi_j},\widetilde{\mb{G}})+Q(\widetilde{\mb{G}},\mb{M}_{\xi_j})\big)\Big\} \f{\widetilde{\mb{G}}_{\xi_j}}{\mb{M}_\star }\\
\di\quad+\Big\{\f1v\big(-\big(\xi_1(n_2 v\mb{M})_x\big)_{\xi_j}+u_1(n_2 v\mb{M})_{\xi_jx}\big)
+\f1v\big(-\big(\xi_1(\mb{P}_cF_2)_x\big)_{\xi_j}+u_1(\mb{P}_cF_2)_{\xi_jx}\big)\\
\di\quad -(n_2 v\mb{M}_{\xi_j})_t-\f{\Phi_x}{v}(\mb{M}_{\xi_j\xi_1}+\bar{\mb{G}}_{\xi_j\xi_1})
+2\pa_{\xi_j}Q(F_2,\mb{G})+2Q(\mb{P}_c F_2,\mb{M}_{\xi_j})\Big\}\f{(\mb{P}_cF_2)_{\xi_j}}{\mb{M}_\star }d\xi dxds\Big|\\
\di\quad +\int_0^t\int\int \Big(-\f{\Phi_x}{v}(\mb{P}_cF_2)_{\xi_j\xi_1}\f{\widetilde{\mb{G}}_{\xi_j}}{\mb{M}_\star }
-\f{\Phi_x}{v}\widetilde{\mb{G}}_{\xi_j\xi_1}\f{(\mb{P}_cF_2)_{\xi_j}}{\mb{M}_\star }\Big)d\xi dxds\\
\di\quad\leq \int\int \f{|(\widetilde{\mb{G}}, \mb{P}_c F_2)_{\xi_j}|^2}{2\mb{M}_\star }(x,\xi,0)
d\xi dx+\sum_{i=1}^{14}T_i.
\end{array}
\end{equation}
Now we estimate the terms on the right hand side of \eqref{GvE0}.
First, by the fact \eqref{Fact1}, it holds that
\begin{equation}\label{T1}
\begin{array}{ll}
\di |T_1|\leq
\f{\tilde\sigma}{16}\int_0^t\int\int\f{\nu(|\xi|)|\widetilde{\mb{G}}_{\xi_j}|^2}{\mb{M}_\star }
d\xi dxds+C\int_0^t\|(\psi_x,\zeta_x)\|^2ds,
\end{array}
\end{equation}
\begin{equation}\label{T2}
\begin{array}{ll}
\di |T_2|=|\int_0^t\int\int
\Big(\f{u_1-\xi_1}{v}\mb{G}_{\xi_jx}-\f{\d_{1j}}{v}\mb{G}_x+\f1v\sum_{k=0}^4<\chi_k,\xi_1\mb{G}_x>(\chi_k)_{\xi_j}\Big)
\f{\widetilde{\mb{G}}_{\xi_j}}{\mb{M}_\star }d\xi dxds|\\
\di\quad =|\int_0^t\int\int
\Big(-\big(\f{u_1-\xi_1}{v}\big)_x\f{|\widetilde{\mb{G}}_{\xi_j}|^2}{2\mb{M}_\star }
+\f{u_1-\xi_1}{v}\bar{\mb{G}}_{\xi_jx}\f{\widetilde{\mb{G}}_{\xi_j}}{\mb{M}_\star }\\
\di\qquad+\f1v\big(-\d_{1j}\mb{G}_x+\sum_{k=0}^4<\chi_k,\xi_1\mb{G}_x>(\chi_k)_{\xi_j}\big)\f{\widetilde{\mb{G}}_{\xi_j}}{\mb{M}_\star }\Big)d\xi dxds|\\
 \di\quad
\leq
(\f{\tilde\sigma}{16}+C(\chi+\d))\int_0^t\int\int\f{\nu(|\xi|)}{\mb{M}_\star }|\widetilde{\mb{G}}_{\xi_j}|^2
d\xi dxds+C\int_0^t\int\int\f{\nu(|\xi|)}{\mb{M}_\star }|\mb{G}_x|^2d\xi dxds\\
\di\qquad +C\int_0^t\int\Big[|(\bar\t_{xx},\bar u_{1xx})|^2+|(\bar\t_x,\bar u_{1x})(v_x,u_x,\t_x)|^2\Big]dxds\\
\di\quad \leq
(\f{\tilde\sigma}{16}+C(\chi+\d))\int_0^t\int\int\f{\nu(|\xi|)}{\mb{M}_\star }|\widetilde{\mb{G}}_{\xi_j}|^2
d\xi dxds+C\int_0^t\int\int\f{\nu(|\xi|)}{\mb{M}_\star }|\mb{G}_x|^2d\xi dxds\\
\di\qquad+C(\chi+\d)\int_0^t\|(\phi_x,\psi_x,\zeta_x)\|^2ds+C\d^{\frac12},
\end{array}
\end{equation}
and
\begin{equation}\label{T4}
\begin{array}{ll}
\di |T_3|+|T_4|\leq
\f{\tilde\sigma}{16}\int_0^t\int\int\f{\nu(|\xi|)|\widetilde{\mb{G}}_{\xi_j}|^2}{\mb{M}_\star }
d\xi dxds +C(\chi+\d)\int_0^t\|(\Phi_x, n_2)\|^2ds\\
\di\qquad\qquad\quad+C(\chi+\d)\int_0^t\int\int\f{\nu(|\xi|)|\mb{P}_c F_2|^2}{\mb{M}_\star }
d\xi dxds.
\end{array}
\end{equation}
Next, it holds that
\begin{equation}\label{T5}
\begin{array}{ll}
\di |T_{5}| \leq
\f{\tilde\sigma}{16}\int_0^t\int\int\f{\nu(|\xi|)|\widetilde{\mb{G}}_{\xi_j}|^2}{\mb{M}_\star }
d\xi dxds+C\int_0^t\int\Big(\int\f{\nu(|\xi|)|\mb{G}_{\xi_j}|^2}{\mb{M}_\star }
d\xi\int\f{|\mb{G}|^2}{\mb{M}_\star } d\xi\Big)dxds\\
\di \qquad\qquad
+C\int_0^t\int\Big(\int\f{\nu(|\xi|)|\mb{G}|^2}{\mb{M}_\star }
d\xi\int\f{|\mb{G}_{\xi_j}|^2}{\mb{M}_\star } d\xi\Big)dxds\\
\di \qquad \leq
(\f{\tilde\sigma}{16}+C(\chi+\d))\int_0^t\int\int\f{\nu(|\xi|)}{\mb{M}_\star }|\widetilde{\mb{G}}_{\xi_j}|^2
d\xi dxds+C(\chi+\delta)
\int_0^t\int\int\f{\nu(|\xi|)}{\mb{M}_\star }|\widetilde{\mb{G}}|^2d\xi dxds+C\delta^{\frac12},
\end{array}
\end{equation}
and
\begin{equation}\label{T6}
\begin{array}{ll}
\di |T_6|\leq \f{\tilde\sigma}{16}\int_0^t\int\int\f{\nu(|\xi|)|\widetilde{\mb{G}}_{\xi_j}|^2}{\mb{M}_\star }
d\xi dxds+C\int_0^t\int\Big[|(\bar\t_{xt},\bar u_{1xt})|^2+|(\bar\t_x,\bar u_{1x})(v_t,u_t,\t_t)|^2\Big]dxds\\
\di\qquad \leq \f{\tilde\sigma}{16}\int_0^t\int\int\f{\nu(|\xi|)}{\mb{M}_\star }|\widetilde{\mb{G}}_{\xi_j}|^2
d\xi dxds+C(\chi+\d)\int_0^t\|(\phi_t,\psi_t,\zeta_t)\|^2ds+C\d^{\frac12}.
\end{array}
\end{equation}
Similarly, we can get
\begin{equation}\label{T7}
\begin{array}{ll}
\di |T_7|\leq \f{\tilde\sigma}{16}\int_0^t\int\int\f{\nu(|\xi|)}{\mb{M}_\star }|\widetilde{\mb{G}}_{\xi_j}|^2
+C\int_0^t\int\int\f{\nu(|\xi|)}{\mb{M}_\star }|\widetilde{\mb{G}}|^2d\xi dxds,
\end{array}
\end{equation}
and
\begin{equation}\label{T8}
\begin{array}{ll}
\di |T_8|\leq \f{\tilde\sigma}{16}\int_0^t\int\int\f{\nu(|\xi|)}{\mb{M}_\star }|\widetilde{\mb{G}}_{\xi_j}|^2d\xi dxds
+C\int_0^t\int\big[|n_{2x}|^2+n_2^2|(v_x,u_x,\t_x)|^2 \big] dxds\\
\di\qquad \leq \f{\tilde\sigma}{16}\int_0^t\int\int\f{\nu(|\xi|)}{\mb{M}_\star }|\widetilde{\mb{G}}_{\xi_j}|^2d\xi dxds
+C\int_0^t\|n_{2x}\|^2ds+C(\chi+\d)\int_0^t\|n_2\|^2ds.
\end{array}
\end{equation}
It holds that
\begin{equation}\label{T9}
\begin{array}{l}
 \di|T_9| =|\int_0^t\int\int \Big[\f{u_1-\xi_1}{v}(\mb{P}_cF_2)_{\xi_jx}-\f{\d_{1j}}{v}(\mb{P}_cF_2)_x\Big]\f{(\mb{P}_cF_2)_{\xi_j}}{\mb{M}_\star } d\xi dxds|\\
 \di\qquad=|\int_0^t\int\int\Big[\big(\f{u_1-\xi_1}{v}\big)_x\f{|(\mb{P}_cF_2)_{\xi_j}|^2}{2\mb{M}_\star }
 +\f{\d_{1j}}{v}(\mb{P}_cF_2)_x\f{(\mb{P}_c F_2)_{\xi_j}}{\mb{M}_\star }\Big] d\xi dxds|\\
 \di\qquad \leq (\f{\tilde\sigma}{16}+C(\chi+\d))\int_0^t\int\int\f{\nu(|\xi|)}{\mb{M}_\star }|(\mb{P}_cF_2)_{\xi_j}|^2
d\xi dxds+C\int_0^t\int\int\f{\nu(|\xi|)}{\mb{M}_\star }|(\mb{P}_c F_2)_x|^2d\xi dxds,
\end{array}
\end{equation}
\begin{equation}\label{T10}
 |T_{10}| \leq \f{\tilde\sigma}{16}
\int_0^t\int\int\f{\nu(|\xi|)}{\mb{M}_\star }|(\mb{P}_cF_2)_{\xi_j}|^2d\xi dxds
+C\int_0^t\|n_{2t}\|^2ds+C(\chi+\delta)\int_0^t\|n_2\|^2ds,
\end{equation}
and
\begin{equation}\label{T11}
 |T_{11}| \leq \f{\tilde\sigma}{16}
\int_0^t\int\int\f{\nu(|\xi|)}{\mb{M}_\star }|(\mb{P}_cF_2)_{\xi_j}|^2d\xi dx
ds+C\int_0^t\|\Phi_x\|^2 ds.
\end{equation}
Furthermore, we have
\begin{equation}\label{T12}
\begin{array}{ll}
\di
 |T_{12}| \leq \f{\tilde\sigma}{16}
\int_0^t\int\int\f{\nu(|\xi|)}{\mb{M}_\star }|(\mb{P}_cF_2)_{\xi_j}|^2d\xi dxds\\
\di\qquad\quad+C\int_0^t\int\int\f{\nu(|\xi|)^{-1}}{\mb{M}_\star }\big(|Q(F_{2\xi_j},\mb{G})|^2+|Q(F_2,\mb{G}_{\xi_j})|^2\big) d\xi dxds
 \\
\di\qquad\leq \f{\tilde\sigma}{16}
\int_0^t\int\int\f{\nu(|\xi|)}{\mb{M}_\star }|(\mb{P}_cF_2)_{\xi_j}|^2d\xi dxds\\
\di\qquad\quad+C\int_0^t\int\big(\int\f{\nu(|\xi|)}{\mb{M}_\star }|F_{2\xi_j}|^2d\xi\cdot\int\f{|\mb{G}|^2}{\mb{M}_\star }d\xi
+\int\f{|F_{2\xi_j}|^2}{\mb{M}_\star }d\xi\cdot\int\f{\nu(|\xi|)}{\mb{M}_\star }|\mb{G}|^2d\xi\big)dxds\\
\di\qquad\quad+C\int_0^t\int\big(\int\f{\nu(|\xi|)}{\mb{M}_\star }|F_{2}|^2d\xi\cdot\int\f{|\mb{G}_{\xi_j}|^2}{\mb{M}_\star }d\xi
+\int\f{|F_{2}|^2}{\mb{M}_\star }d\xi\cdot\int\f{\nu(|\xi|)}{\mb{M}_\star }|\mb{G}_{\xi_j}|^2d\xi\big)dxds\\
\di\qquad \leq \f{\tilde\sigma}{16}
\int_0^t\int\int\f{\nu(|\xi|)}{\mb{M}_\star }|(\mb{P}_cF_2)_{\xi_j}|^2d\xi dxds\\
\di\qquad\quad +C\int_0^t\int\big(|n_2|^2+\int\f{\nu(|\xi|)}{\mb{M}_\star }|(\mb{P}_c F_2)_{\xi_j}|^2d\xi\big)\cdot
\big(|(\bar\t_x,\bar u_{1x})|^2+\int\f{|\widetilde{\mb{G}}|^2}{\mb{M}_\star }d\xi\big)dxds\\
\di\qquad\quad+C\int_0^t\int\big(|n_2|^2+\int\f{|(\mb{P}_cF_2)_{\xi_j}|^2}{\mb{M}_\star }d\xi\big)\cdot
\big(|(\bar\t_x,\bar u_{1x})|^2+\int\f{\nu(|\xi|)}{\mb{M}_\star }|\widetilde{\mb{G}}|^2d\xi\big)dxds\\
\end{array}
\end{equation}
and
\begin{equation}\label{T13}
\begin{array}{ll}
\di |T_{13}|\leq \f{\tilde\sigma}{16}
\int_0^t\int\int\f{\nu(|\xi|)}{\mb{M}_\star }|(\mb{P}_cF_2)_{\xi_j}|^2d\xi dxds+C\int_0^t\int\int\f{\nu(|\xi|)}{\mb{M}_\star }|\mb{P}_c F_2|^2 d\xi dxds.
\end{array}
\end{equation}
 For $T_{14}$, it holds that
\begin{equation}\label{T14}
\begin{array}{ll}
\di |T_{14}|=|\int_0^t\int\int \Big(\f{\Phi_x}{v}(\mb{P}_cF_2)_{\xi_j\xi_1}\f{\widetilde{\mb{G}}_{\xi_j}}{\mb{M}_\star }
+\f{\Phi_x}{v}\widetilde{\mb{G}}_{\xi_j\xi_1}\f{(\mb{P}_cF_2)_{\xi_j}}{\mb{M}_\star }\Big)d\xi dxds|\\
\di ~~\qquad =|\int_0^t\int\int\f{\Phi_x}{v}\widetilde{\mb{G}}_{\xi_j}\cdot(\mb{P}_cF_2)_{\xi_j}\big(\f{1}{\mb{M}_\star }\big)_{\xi_j}d\xi dx ds|\\
\di ~~\qquad \leq C\chi \int_0^t\int\int\f{\nu(|\xi|)}{\mb{M}_\star }|(\widetilde{\mb{G}},\mb{P}_cF_2)_{\xi_j}|^2d\xi dx ds.
\end{array}
\end{equation}
Substituting \eqref{T1}-\eqref{T14} into \eqref{GvE0} and choosing
$\chi,\delta$ sufficiently small gives that
\begin{equation}\label{GvE1}
\begin{array}{ll}
\di \int\int \f{|(\widetilde{\mb{G}}, \mb{P}_c F_2)_{\xi_j}|^2}{\mb{M}_\star }(x,\xi,t)
d\xi dx+\int_0^t\int\int\f{\nu(|\xi|)|(\widetilde{\mb{G}}, \mb{P}_c F_2)_{\xi_j}|^2}{\mb{M}_\star }
d\xi dxds
\\[2mm]
 \di
\leq C\int\int \f{|(\widetilde{\mb{G}}, \mb{P}_c
F_2)_{\xi_j}|^2}{\mb{M}_\star }(x,\xi,0)
d\xi dx +C(\chi+\delta)
\int_0^t\Big[\|n_2\|^2+\sum_{|\alpha'|=1}\|\partial^{\alpha'}(\phi,\psi,\zeta)\|^2\Big]ds+C\delta^{\frac12} \\
\di +C\int_0^t\|(\psi_x,\zeta_x,\Phi_x,n_{2x},n_{2t})\|^2ds+C\int_0^t\int\int\f{\nu(|\xi|)}{\mb{M}_\star }|(\widetilde{\mb{G}},\mb{P}_c F_2,\mb{G}_x,(\mb{P}_c F_2)_x)|^2 d\xi dxds.
\end{array}
\end{equation}
Applying $\partial_{\xi_j}\partial_{\xi_k}~(j,k=1,2,3)$ to the equation \eqref{GE} and \eqref{F2-pc-1},
and using the similar analysis as in obtaining \eqref{GvE1}, one has for $|\beta|=2$
\begin{equation}\label{GvE2}
\begin{array}{ll}
\di \int\int \f{|\partial^\beta(\widetilde{\mb{G}}, \mb{P}_c F_2)|^2}{\mb{M}_\star }(x,\xi,t)
d\xi dx+\int_0^t\int\int\f{\nu(|\xi|)|\partial^\beta(\widetilde{\mb{G}}, \mb{P}_c F_2)|^2}{\mb{M}_\star }
d\xi dxds
\\[2mm]
 \di
\leq C\int\int \f{|\partial^\beta(\widetilde{\mb{G}}, \mb{P}_c F_2)|^2}{\mb{M}_\star }(x,\xi,0)
d\xi dx+C\int_0^t\|(\psi_x,\zeta_x,\Phi_x,n_{2x},n_{2t})\|^2ds\\
\di
+C(\chi+\delta)\int_0^t\Big[\sum_{|\alpha'|=1}\|\pa^{\alpha'}(\phi,\psi,\z)\|^2+\|n_2\|^2\Big]ds+C\delta^{\frac12}\\
\di +C\sum_{0\leq|\beta^\prime|\leq 1}\int_0^t\int\int\f{\nu(|\xi|)|\partial^{\beta^\prime}(\widetilde{\mb{G}},\mb{P}_cF_2, \mb{G}_x, (\mb{P}_cF_2)_x)|^2}{\mb{M}_\star }
d\xi dxds.
\end{array}
\end{equation}
By \eqref{LE-A}, it holds that
\begin{equation}\label{A1}
\begin{array}{ll}
\di \int_0^t\|(\psi_x,\zeta_x,\Phi_x, n_{2x}, n_{2t}\|^2ds+\int_0^t\int\int\f{\nu(|\xi|)|(\widetilde{\mb{G}},\mb{P}_c F_2)|^2}{\mb{M}_\star }d\xi dxds\\
\di \leq C\|(\phi_0,\psi_0,\zeta_0,\phi_{0x},\Phi_{0x},n_{20})\|^2
+C\int\int\f{|(\widetilde{\mb{G}},\mb{P}_c F_{2})|^2}{\mb{M}_\star }(x,\xi,0) d\xi dx+C\delta\\
\di +C\sum_{1\leq|\alpha|\leq 2}\int_0^t\int\int \f{\nu(|\xi|)|\partial^\alpha(\mb{G},\mb{P}_c F_2)|^2}{\mb{M}_\star }d\xi dxds
+C(\chi+\delta)\int_0^t\int\int\f{\nu(|\xi|)|\pa_{\xi_1}(\mb{P}_c F_2)_x|^2}{\mb{M}_\star } d\xi dxds\\
\di+C(\chi+\d)\int_0^t\int\int \f{\nu(|\xi|)|\pa_{\xi_1}(\widetilde{\mb{G}},\mb{P}_c F_2)|^2}{\mb{M}_\star }d\xi dxds,
\end{array}
\end{equation}
which together with \eqref{GvE2} gives
\begin{equation}\label{wt5.4}
\begin{array}{ll}
\di \sum_{1\leq |\beta|\leq 2}\int\int \f{|\partial^\beta(\widetilde{\mb{G}}, \mb{P}_c F_2)|^2}{\mb{M}_\star }(x,\xi,t)
d\xi dx+\sum_{1\leq |\beta|\leq 2}\int_0^t\int\int\f{\nu(|\xi|)|\partial^\beta(\widetilde{\mb{G}}, \mb{P}_c F_2)|^2}{\mb{M}_\star }
d\xi dxds
\\[2mm]
 \di
\leq C\mathcal{N}^2(0)+C\d^{\frac12}+C(\chi+\delta)\int_0^t\Big[\sum_{|\alpha'|=1}\|\pa^{\alpha'}(\phi,\psi,\z)\|^2+\|n_2\|^2\Big]ds\\
\di +C\sum_{1\leq|\alpha|\leq 2}\int_0^t\int\int\f{\nu(|\xi|)|\partial^{\alpha}(\mb{G},\mb{P}_cF_2)|^2}{\mb{M}_\star }
d\xi dxds+C\int_0^t\int\int \f{\nu(|\xi|)|\pa_{\xi_j}(\mb{G},\mb{P}_cF_2)_x|^2}{\mb{M}_\star }d\xi dxds.
\end{array}
\end{equation}
Similarly, one has for $|\alpha^\prime|=1$ and $|\beta^\prime|=1$
\begin{equation}\label{GvE11}
%
\end{equation}
By \eqref{LE-A}, \eqref{55e} and \eqref{n2-he}, it holds that
\begin{equation}\label{A2}
%
\end{equation}
Substituting \eqref{A2} into \eqref{GvE11}, which together with \eqref{xt-e} and \eqref{wt5.4}, one has
\begin{equation}\label{G3}
%
\end{equation}

\

\underline{Step 3. Estimation on the highest order estimates.}

\

Finally, we derive the highest order estimate to control the terms
$\displaystyle \int\phi_{xx}\psi_{1x}dx$ and  \\ $\displaystyle
\sum_{|\alpha|=2}\int\int\frac{\nu(|
\xi|)}{\mb{M}_\star }|\partial^\alpha (\mathbf{G}, \mb{P}_cF_2)|^2d \xi
dx$.  Let $|\alpha|=2$, from \eqref{macro} and \eqref{vpb-l}$_3$, it holds that
$$
\|\partial^\alpha(\rho,m,E,n_2)\|^2\leq
C\int\int\frac{|\partial^\alpha (F_1,F_2)|^2}{\mb{M}_\star }d \xi dx.
$$
So we have for $|\alpha|=2$,
$$
\|\partial^\alpha(\phi,\psi,\zeta,n_2)\|^2\leq
C\int\int\frac{|\partial^\alpha (F_1,F_2)|^2}{\mb{M}_\star }d \xi
dx+C\delta.
$$
Thus in order to estimate $\displaystyle
\int\phi_{xx}\psi_{1x}dx$, it is sufficient to estimate
$\displaystyle \int\int\frac{|\partial^\alpha
(F_1,F_2)|^2}{\mb{M}_\star }d \xi dx$ $(|\alpha|=2)$. From \eqref{vpb-l}, one has
\begin{equation}\label{wt5.5}
\left\{
\begin{array}{ll}
\di
F_{1t}+\frac{\xi_1-u_1}{v} F_{1x}+\frac{\Phi_x}{v}\partial_{\xi_1}F_2= \mb{L}_\mathbf{M}\mathbf{G}+2Q(\mb{G},\mb{G}), \\[2mm]
\di
F_{2t}+\frac{\xi_1-u_1}{v} F_{2x}+\frac{\Phi_x}{v}\partial_{\xi_1}F_1= \mb{N}_{\mb{M}}\mb{P}_c F_2+2Q(F_2,\mb{G}).
\end{array}
\right.
\end{equation}
Applying $\partial^\alpha(|\alpha|=2)$ on the Vlasov-Poisson-Boltzmann equation \eqref{wt5.5}, we have
\begin{equation}\label{F1-a2}
\begin{array}{ll}
\di (\partial^\alpha F_1)_t+ \f{\xi_1-u_1}{v}(\partial^\alpha F_1)_x+\f{\Phi_x}{v}\partial_{\xi_1}(\partial^\alpha F_2)
+\pa^{\alpha}\big(\f{\xi_1-u_1}{v}\big)F_{1x}+\partial^\alpha\big(\f{\Phi_x}{v}\big)\partial_{\xi_1}F_2\\
\di \quad +\sum_{|\alpha^\prime|=1,\alpha^\prime\leq\alpha}C_{\alpha}^{\alpha^\prime}\partial^{\alpha^\prime}
\big(\f{\xi_1-u_1}{v}\big)\partial^{\alpha-\alpha^\prime}F_{1x} +\sum_{|\alpha^\prime|=1,\alpha^\prime\leq\alpha}C_{\alpha}^{\alpha^\prime}\partial^{\alpha^\prime}
\big(\f{\Phi_x}{v}\big)\partial_{\xi_1}\partial^{\alpha-\alpha^\prime}F_2\\
\di=\mb{L}_\mathbf{M}\partial^\alpha\mathbf{G}+\sum_{|\alpha^\prime|=1,\alpha^\prime\leq\alpha}2C_{\alpha}^{\alpha^\prime}
\big(Q(\pa^{\alpha'}\mb{M},\pa^{\alpha-\alpha'}\mb{G})+Q(\pa^{\alpha'}\mb{G},\pa^{\alpha-\alpha'}\mb{M})\big)\\
\di\quad+2Q(\mb{G},\pa^{\alpha}\mb{M})+2Q(\pa^{\alpha}\mb{M},\mb{G})+2 \partial^\alpha Q(\mathbf{G}, \mathbf{G}).
\end{array}
\end{equation}
and
\begin{equation}\label{F2-a2}
\begin{array}{ll}
\di (\partial^\alpha F_2)_t+ \f{\xi_1-u_1}{v}(\partial^\alpha F_2)_x+\f{\Phi_x}{v}\partial_{\xi_1}(\partial^\alpha F_1)
+\pa^{\alpha}\big(\f{\xi_1-u_1}{v}\big)F_{2x}+\partial^\alpha\big(\f{\Phi_x}{v}\big)\partial_{\xi_1}F_1\\
\di \quad +\sum_{|\alpha^\prime|=1,\alpha^\prime\leq\alpha}C_{\alpha}^{\alpha^\prime}\partial^{\alpha^\prime}
\big(\f{\xi_1-u_1}{v}\big)\partial^{\alpha-\alpha^\prime}F_{2x} +\sum_{|\alpha^\prime|=1,\alpha^\prime\leq\alpha}C_{\alpha}^{\alpha^\prime}\partial^{\alpha^\prime}
\big(\f{\Phi_x}{v}\big)\partial_{\xi_1}\partial^{\alpha-\alpha^\prime}F_1\\
\di=\mb{N}_\mathbf{M}\partial^\alpha(\mathbf{P}_cF_2)+\sum_{|\alpha^\prime|=1,\alpha^\prime\leq\alpha}2C_{\alpha}^{\alpha^\prime}
Q(\pa^{\alpha'}(\mb{P}_cF_2),\pa^{\alpha-\alpha'}\mb{M})+2Q(\mb{P}_c F_2,\pa^{\alpha}\mb{M})+2 \partial^\alpha Q(F_2, \mathbf{G}).
\end{array}
\end{equation}
Multiplying \eqref{F1-a2} by $\frac{\partial^\alpha F_1}{\mb{M}_\star }=\frac{\partial^\alpha
\mathbf{M}}{\mb{M}_\star }+\frac{\partial^\alpha
\mathbf{G}}{\mb{M}_\star }$, and \eqref{F2-a2} by $\frac{\partial^\alpha
F_2}{\mb{M}_\star }=\frac{\partial^\alpha
(n_2 v\mathbf{M})}{\mb{M}_\star }+\frac{\partial^\alpha
(\mathbf{P}_c F_2)}{\mb{M}_\star }$, respectively,  gives
\begin{equation}\label{h1}
\begin{array}{l}
\quad\displaystyle (\frac{|\partial^\alpha
F_1|^2}{2\mb{M}_\star })_t-\frac{\partial^\alpha
\mathbf{G}}{\mb{M}_\star }\mb{L}_\mathbf{M}\partial^\alpha \mathbf{G}
=-\f{\Phi_x}{v}\frac{\partial^\alpha
F_1}{\mb{M}_\star }\partial_{\xi_1}(\partial^\alpha F_2)+\frac{\partial^\alpha
F_1}{\mb{M}_\star }\Big[-\pa^{\alpha}(\f{\xi_1-u_1}{v})F_{1x}\\[2mm]
\di -\partial^\alpha(\f{\Phi_x}{v})\partial_{\xi_1}F_2
-\sum_{|\alpha^\prime|=1,\alpha^\prime\leq\alpha}C_{\alpha}^{\alpha^\prime}\partial^{\alpha^\prime}(\f{\xi_1-u_1}{v})\pa^{\alpha-\alpha'}F_{1x}
-\sum_{|\alpha^\prime|=1,\alpha^\prime\leq\alpha}C_{\alpha}^{\alpha^\prime}\partial^{\alpha^\prime}(\f{\Phi_x}{v})\partial_{\xi_1}\partial^{\alpha-\alpha^\prime}F_2\\
\di +\sum_{|\alpha^\prime|=1,\alpha^\prime\leq\alpha}2C_{\alpha}^{\alpha^\prime}
\big(Q(\pa^{\alpha'}\mb{M},\pa^{\alpha-\alpha'}\mb{G})+Q(\pa^{\alpha'}\mb{G},\pa^{\alpha-\alpha'}\mb{M})\big)
+2Q(\mb{G},\pa^{\alpha}\mb{M})+2Q(\pa^{\alpha}\mb{M},\mb{G})\\
\di+2 \partial^\alpha Q(\mathbf{G}, \mathbf{G})\Big]+\frac{\partial^\alpha
\mathbf{M}}{\mb{M}_\star }\mb{L}_\mathbf{M}\partial^\alpha\mathbf{G}-\big(\f{\xi_1-u_1}{v}\f{|\pa^{\alpha }F_1|^2}{2\mb{M}_\star }\big)_x
+\big(\f{\xi_1-u_1}{v}\big)_x\f{|\pa^{\alpha}F_1|^2}{2\mb{M}_\star }
\end{array}
\end{equation}
and
\begin{equation}\label{h2}
\begin{array}{l}
\quad\displaystyle (\frac{|\partial^\alpha
F_2|^2}{2\mb{M}_\star })_t-\frac{\partial^\alpha (\mathbf{P}_cF_2)}{\mb{M}_\star }\mb{N}_\mathbf{M}\partial^\alpha(\mathbf{P}_cF_2)
=-\f{\Phi_x}{v}\frac{\partial^\alpha F_2}{\mb{M}_\star }\partial_{\xi_1}(\partial^\alpha F_1)
+\frac{\partial^\alpha F_2}{\mb{M}_\star }\Big[-\pa^{\alpha}\big(\f{\xi_1-u_1}{v}\big)F_{2x}\\[2mm]
\di -\partial^\alpha(\f{\Phi_x}{v})\partial_{\xi_1}F_1
-\sum_{|\alpha^\prime|=1,\alpha^\prime\leq\alpha}C_{\alpha}^{\alpha^\prime}\partial^{\alpha^\prime}
\big(\f{\xi_1-u_1}{v}\big)\partial^{\alpha-\alpha^\prime}F_{2x} -\sum_{|\alpha^\prime|=1,\alpha^\prime\leq\alpha}C_{\alpha}^{\alpha^\prime}\partial^{\alpha^\prime}
\big(\f{\Phi_x}{v}\big)\partial_{\xi_1}\partial^{\alpha-\alpha^\prime}F_1\\
\di +\sum_{|\alpha^\prime|=1,\alpha^\prime\leq\alpha}2C_{\alpha}^{\alpha^\prime}
Q(\pa^{\alpha'}(\mb{P}_cF_2),\pa^{\alpha-\alpha'}\mb{M})+2Q(\mb{P}_c F_2,\pa^{\alpha}\mb{M})+2 \partial^\alpha Q(F_2, \mathbf{G})\Big]\\
\di +\f{\pa^{\alpha}(n_2 v\mb{M})}{\mb{M}_\star }\mb{N}_{\mb{M}}\pa^{\alpha}(\mb{P}_c F_2)-\big(\f{\xi_1-u_1}{v}\f{|\pa^{\alpha}F_2|^2}{2\mb{M}_\star }\big)_x
+(\f{\xi_1-u_1}{v})_x\f{|\pa^{\alpha}F_2|^2}{2\mb{M}_\star }.
\end{array}
\end{equation}
Adding \eqref{h1} and \eqref{h2} together and noting that
$$
-\f{\Phi_x}{v}\frac{\partial^\alpha F_1}{\mb{M}_\star }\partial_{\xi_1}(\partial^\alpha F_2)
-\f{\Phi_x}{v}\frac{\partial^\alpha F_2}{\mb{M}_\star }\partial_{\xi_1}(\partial^\alpha F_1)
=-(\f{\Phi_x}{v}\frac{\partial^\alpha F_1\partial^\alpha F_2}{\mb{M}_\star })_{\xi_1}
+\f{\Phi_x}{v}\partial^\alpha F_1\partial^\alpha F_2\big(\f{1}{\mb{M}_\star }\big)_{\xi_1},
$$
and then integrating the resulting equation over $x, \xi, t$ imply that
\begin{equation}\label{h3}
\begin{array}{ll}
\displaystyle\quad \int\int \frac{|\partial^\alpha
F_1|^2+|\partial^\alpha F_2|^2}{2\mb{M}_\star }(x,\xi,t)d\xi dx- \int\int \frac{|\partial^\alpha
F_{10}|^2+|\partial^\alpha F_{20}|^2}{2\mb{M}_\star }d\xi dx\\
\di\quad -\int_0^t\int\int\Big[ \frac{\partial^\alpha
\mathbf{G}}{\mb{M}_\star }\mb{L}_\mathbf{M}\partial^\alpha \mathbf{G}+\frac{\partial^\alpha
(\mathbf{P}_cF_2)}{\mb{M}_\star }\mb{N}_\mathbf{M}\partial^\alpha(\mathbf{P}_cF_2) \Big]d\xi dxds \\
\di =\int_0^t\int\int\Big\{
\f{\Phi_x}{v}\partial^\alpha F_1\partial^\alpha F_2\big(\f{1}{\mb{M}_\star }\big)_{\xi_1}
+\frac{\partial^\alpha F_1}{\mb{M}_\star }\Big[-\pa^{\alpha}(\f{\xi_1-u_1}{v})F_{1x}
-\partial^\alpha(\f{\Phi_x}{v})\partial_{\xi_1}F_2\\
\di\quad -\sum_{|\alpha^\prime|=1,\alpha^\prime\leq\alpha}C_{\alpha}^{\alpha^\prime}\partial^{\alpha^\prime}(\f{\xi_1-u_1}{v})\pa^{\alpha-\alpha'}F_{1x}
-\sum_{|\alpha^\prime|=1,\alpha^\prime\leq\alpha}C_{\alpha}^{\alpha^\prime}\partial^{\alpha^\prime}(\f{\Phi_x}{v})\partial_{\xi_1}\partial^{\alpha-\alpha^\prime}F_2\\
\di\quad +\sum_{|\alpha^\prime|=1,\alpha^\prime\leq\alpha}2C_{\alpha}^{\alpha^\prime}
\big(Q(\pa^{\alpha'}\mb{M},\pa^{\alpha-\alpha'}\mb{G})+Q(\pa^{\alpha'}\mb{G},\pa^{\alpha-\alpha'}\mb{M})\big)
+2\big(Q(\mb{G},\pa^{\alpha}\mb{M})+Q(\pa^{\alpha}\mb{M},\mb{G})\big)\\
\di\quad+2 \partial^\alpha Q(\mathbf{G}, \mathbf{G})\Big]
+\frac{\partial^\alpha F_2}{\mb{M}_\star }\Big[-\pa^{\alpha}\big(\f{\xi_1-u_1}{v}\big)F_{2x}
-\partial^\alpha(\f{\Phi_x}{v})\partial_{\xi_1}F_1
-\sum_{|\alpha^\prime|=1,\alpha^\prime\leq\alpha}C_{\alpha}^{\alpha^\prime}\partial^{\alpha^\prime}
\big(\f{\xi_1-u_1}{v}\big)\partial^{\alpha-\alpha^\prime}F_{2x}\\
\di\quad
-\sum_{|\alpha^\prime|=1,\alpha^\prime\leq\alpha}C_{\alpha}^{\alpha^\prime}\partial^{\alpha^\prime}
\big(\f{\Phi_x}{v}\big)\partial_{\xi_1}\partial^{\alpha-\alpha^\prime}F_1
+\sum_{|\alpha^\prime|=1,\alpha^\prime\leq\alpha}2C_{\alpha}^{\alpha^\prime}
Q(\pa^{\alpha'}(\mb{P}_cF_2),\pa^{\alpha-\alpha'}\mb{M})+2Q(\mb{P}_c F_2,\pa^{\alpha}\mb{M})\\
\di\quad+2 \partial^\alpha Q(F_2, \mathbf{G})\Big]
+\frac{\partial^\alpha \mathbf{M}}{\mb{M}_\star }\mb{L}_\mathbf{M}\partial^\alpha\mathbf{G}
+\f{\pa^{\alpha}(n_2 v\mb{M})}{\mb{M}_\star }\mb{N}_{\mb{M}}\pa^{\alpha}(\mb{P}_c F_2)
+\big(\f{\xi_1-u_1}{v}\big)_x\f{|\pa^{\alpha}(F_1,F_2)|^2}{2\mb{M}_\star }
\Big\} d\xi dxds\\
\di\quad:=\sum_{i=1}^{18} W_i.
\end{array}
\end{equation}
We only estimate the terms $W_1, W_2, W_3, W_{10}, W_{16}, W_{17}$ on the right hand side of \eqref{h3}, while the other terms can be estimated similarly as in \cite{LWYZ}.
First, it holds that
\begin{equation}\label{W1}
\begin{array}{ll}
\di W_1=\int_0^t \int\int \f{\Phi_x}{v}\f{\partial^\alpha(\mb{M}+\mb{G})\cdot \partial^\alpha(n_2 v\mb{M}+\mb{P}_cF_2)}{\mb{M}_\star }\cdot\f{\xi_1-u_{1*}}{R\t} d\xi dxds\\
\di\qquad\leq C\int_0^t\int|\Phi_x|\Big[|\partial^\alpha(v,u,\t)|+\sum_{|\alpha^\prime|=1}|\partial^{\alpha^\prime}(v,u,\t)|^2
+\Big(\int\f{\nu(|\xi|)|\partial^\alpha\mb{G}|^2}{\mb{M}_\star }d\xi\Big)^{\f12}\Big]\\
\di \qquad\qquad \cdot\Big[|\partial^\alpha n_2|+|n_2||\partial^\alpha(v,u,\t)|+\sum_{|\alpha^\prime|=1}\big(|\partial^{\alpha^\prime} n_2||\partial^{\alpha^\prime}(v,u,\t)|+|n_2||\partial^{\alpha^\prime}(v,u,\t)|^2\big)\\
\di\qquad\qquad+\Big(\int\f{\nu(|\xi|)|\partial^\alpha(\mb{P}_cF_2)|^2}{\mb{M}_\star }d\xi\Big)^{\f12}\Big]dxds\\
\di\qquad \leq C(\chi+\delta)\int_0^t\Big[\|\partial^\alpha (\phi,\psi,\zeta,n_2)\|^2
+\sum_{|\alpha^\prime|=1}\|\partial^{\alpha^\prime}(\phi,\psi,\zeta,n_2)\|^2+\|(\Phi_x,n_2)\|^2\\
\di\qquad\qquad +\int\int\f{\nu(|\xi|)|\partial^\alpha(\mb{G},\mb{P}_cF_2)|^2}{\mb{M}_\star }d\xi dx\Big]ds+C\d^{\frac12}.
\end{array}
\end{equation}
Then, we can compute $W_2$ as
\begin{equation}\label{W2}
\begin{array}{ll}
\di W_2=-\int_0^t\int\int\pa^{\alpha}\big(\f{\xi_1-u_1}{v}\big)\f{\pa^{\alpha}\mb{M}+\pa^{\alpha}\mb{G}}{\mb{M}_\star }
(\mb{M}_x+\mb{G}_x) d\xi dxds\\
\di\qquad\leq C\int_0^t\int\big(\pa^{\alpha}(v,u_1)+\sum_{|\alpha'|=1}|\pa^{\alpha'}(v,u_1)|^2\big)\cdot\Big(
|\pa^{\alpha}(v,u,\t)|+\sum_{|\alpha'|=1}|\pa^{\alpha'}(v,u,\t)|^2\\
\di\qquad\qquad+\Big(\int\f{\nu(|\xi|)|\pa^{\alpha}\mb{G}|^2}{\mb{M}_\star }d\xi\Big)^{\f12}\Big)
\cdot\Big(|(v_x,u_x,\t_x)|+\Big(\int\f{|\mb{G}_x|^2}{\mb{M}_\star }d\xi\Big)^{\f12}\Big)dxds\\
\di\qquad \leq C(\chi+\d)\int_0^t\Big[\|\pa^{\alpha}(\phi,\psi,\zeta)\|^2
+\sum_{|\alpha'|=1}\|\pa^{\alpha'}(\phi,\psi,\z)\|^2+\int\int\f{\nu(|\xi|)|\mb{G}_x|^2}{\mb{M}_\star } d\xi dx\\
\di\qquad\qquad +\int\int\f{\nu(|\xi|)|\pa^{\alpha}\mb{G}|^2}{\mb{M}_\star } d\xi dx\Big]ds+C\d^{\frac12}
\end{array}
\end{equation}
and
\begin{equation}\label{W3}
\begin{array}{ll}
\di W_3 =-\int_0^t\int\int \pa^{\alpha}\big(\f{\Phi_x}{v}\big)\f{\pa^{\alpha}\mb{M}+\pa^{\alpha}\mb{G}}{\mb{M}_\star }
\big(n_2 v\mb{M}_{\xi_1}+(\mb{P}_cF_2)_{\xi_1}\big) d\xi dxds\\
\di\qquad\leq C\int_0^t\int\big(|\pa^{\alpha}(v,\Phi_x)|+\sum_{|\alpha'|=1}|\pa^{\alpha'}(v,\Phi_x)|^2\big)\cdot\Big(
|\pa^{\alpha}(v,u,\t)|+\sum_{|\alpha'|=1}|\pa^{\alpha'}(v,u,\t)|^2\\
\di\qquad\qquad+\Big(\int\f{\nu(|\xi|)|\pa^{\alpha}\mb{G}|^2}{\mb{M}_\star }d\xi\Big)^{\f12}\Big)
\cdot\Big(|n_2|+\Big(\int\f{|(\mb{P}_cF_2)_{\xi_1}|^2}{\mb{M}_\star }d\xi\Big)^{\f12}\Big)dxds\\
\di\qquad \leq C(\chi+\d)\int_0^t\Big[\|\pa^{\alpha}(\phi,\psi,\zeta,\Phi_x)\|^2
+\sum_{|\alpha'|=1}\|\pa^{\alpha'}(\phi,\psi,\z,\Phi_x)\|^2+\|n_2\|^2\Big]ds+C\d^{\frac12}\\
\di\qquad+C(\chi+\d)\int_0^t\Big[\int\int\f{\nu(|\xi|)|\pa^{\alpha}\mb{G}|^2}{\mb{M}_\star } d\xi dx
 +\int\int\f{\nu(|\xi|)|\pa_{\xi_1}(\mb{P}_c F_2)|^2}{\mb{M}_\star } d\xi dx\Big]ds.
\end{array}
\end{equation}
Then it holds that,
\begin{equation}\label{W10}
\begin{array}{ll}
\di W_{10} =-\int_0^t\int\int\pa^{\alpha}\big(\f{\Phi_x}{v}\big)\f{\pa^{\alpha}(n_2 v\mb{M})+\pa^{\alpha}(\mb{P}_cF_2)}{\mb{M}_\star }\pa_{\xi_1}F_{1}d\xi dxds\\
\di\qquad=-\int_0^t\int\int\pa^{\alpha}\big(\f{\Phi_x}{v}\big)v\pa^{\alpha}n_2\f{\mb{M}}{\mb{M}_\star }\pa_{\xi_1}F_1d\xi dxds
-\int_0^t\int\int\pa^{\alpha}\big(\f{\Phi_x}{v}\big)\f{\pa^{\alpha}(n_2 v\mb{M})-v\mb{M}\pa^{\alpha}n_2}{\mb{M}_\star }\pa_{\xi_1}F_1d\xi dxds\\
\di\qquad\quad-\int_0^t\int\int\pa^{\alpha}\big(\f{\Phi_x}{v}\big)\f{\pa^{\alpha}(\mb{P}_cF_2)}{\mb{M}_\star }\pa_{\xi_1}F_1d\xi dxds\\
\di\qquad=-\int_0^t\int\int\pa^{\alpha}\big(\f{\Phi_x}{v}\big)v\pa^{\alpha}n_2\big(\f{\mb{M}}{\mb{M}_\star }-1\big)
(\mb{M}_{\xi_1}+\widetilde{\mb{G}}_{\xi_1}+\bar{\mb{G}}_{\xi_1})d\xi dxds\\
\di\qquad\quad-\int_0^t\int\int\pa^{\alpha}\big(\f{\Phi_x}{v}\big)\f{\pa^{\alpha}(n_2 v\mb{M})-v\mb{M}\pa^{\alpha}n_2}{\mb{M}_\star }
(\mb{M}_{\xi_1}+\widetilde{\mb{G}}_{\xi_1}+\bar{\mb{G}}_{\xi_1})d\xi dxds\\
\di\qquad\quad+\int_0^t\int\int\pa^{\alpha}\big(\f{\Phi_x}{v}\big)\f{\xi_1-u_1}{R\t}\big(\f{\mb{M}}{\mb{M}_\star }-1\big)\pa^{\alpha}(\mb{P}_cF_2)d\xi dxds\\
\di\qquad\quad-\int_0^t\int\int\pa^{\alpha}\big(\f{\Phi_x}{v}\big)\f{\pa^{\alpha}(\mb{P}_cF_2)}{\mb{M}_\star }
(\widetilde{\mb{G}}_{\xi_1}+\bar{\mb{G}}_{\xi_1})d\xi dxds
+\int_0^t\int\pa^{\alpha}\big(\f{\Phi_x}{v}\big)\f{\xi_1-u_1}{R\t}\pa^{\alpha}(\mb{P}_cF_2) d\xi dxds,
\end{array}
\end{equation}
where
$$
\begin{array}{l}
\di\int_0^t\int\pa^{\alpha}\big(\f{\Phi_x}{v}\big)\f{\xi_1-u_1}{R\t}\pa^{\alpha}(\mb{P}_cF_2) d\xi dxds
=\int_0^t\int\pa^{\alpha}\big(\f{\Phi_x}{v}\big)\f{1}{R\t}\pa^{\alpha}\big(\int\xi_1\mb{P}_cF_2 d\xi\big)dxds\\
\di=-\int_0^t\int \f{1}{R\t}\pa^{\alpha}\big(\f{\Phi_x}{v}\big)\pa^{\alpha}\big(\f{\Phi_x}{2v}\big)_tdxds
=-\int_0^t\int \f{1}{4R\t}\Big(|\pa^{\alpha}\big(\f{\Phi_x}{v}\big)|^2\Big)_t dxds\\
\di=-\int_0^t\int \f{1}{4R\t}|\pa^{\alpha}\big(\f{\Phi_x}{v}\big)|^2(x,t)dx
+\int_0^t\int \f{1}{4R\t}|\pa^{\alpha}\big(\f{\Phi_x}{v}\big)|^2(x,0)dx
+\int_0^t\int \big(\f{1}{4R\t}\big)_t |\pa^{\alpha}\big(\f{\Phi_x}{v}\big)|^2dxds\\
\di \leq -\int_0^t\int \f{1}{4R\t}|\pa^{\alpha}\big(\f{\Phi_x}{v}\big)|^2(x,t)dx
+C\|\pa^{\alpha}(\Phi_{0x},v_0)\|^2+C\sum_{|\alpha'|=1}\|\pa^{\alpha'}(\Phi_{0x},v_0)\|^2\\
\di\quad +C(\chi+\d)\int_0^t\Big[\|\pa^{\alpha}(\Phi_x,\phi)\|^2+\sum_{|\alpha'|=1}\|\pa^{\alpha'}(\Phi_x,\phi)\|^2\Big]ds+C\d^{\frac12},
\end{array}
$$
which together with \eqref{W10}, we can get
\begin{equation}\label{w10-1}
\begin{array}{l}
\di W_{10}\leq -\int_0^t\int \f{1}{4R\t}|\pa^{\alpha}\big(\f{\Phi_x}{v}\big)|^2(x,t)dx
+C\|\pa^{\alpha}(\Phi_{0x},v_0)\|^2+C\sum_{|\alpha'|=1}\|\pa^{\alpha'}(\Phi_{0x},v_0)\|^2\\
\di\qquad\quad +C(\chi+\d+\eta_0)\int_0^t\Big[\|\pa^{\alpha}(\phi,\psi,\z,\Phi_x,n_2)\|^2
+\int\int\f{\nu(|\xi|)|\pa^{\alpha}(\mb{P}_cF_2)|^2}{\mb{M}_\star } d\xi dx\Big]ds\\
\di\qquad\quad+C(\chi+\d)\sum_{|\alpha'|=1}\int_0^t\|\pa^{\alpha'}(\phi,\psi,\z,\Phi_x,n_2)\|^2ds+C\d^{\frac12}.
\end{array}
\end{equation}
Since $\mathbf{P}_1(\partial^\alpha \mathbf{M})$ does not contain
$\partial^\alpha(v,u,\theta)$. Thus, we have
\begin{equation}\label{hq1}
\begin{array}{l}
\di W_{16}=\int_0^t\int\int\frac{\partial^\alpha \mathbf{M}}{\mb{M}_\star }\mb{L}_\mathbf{M}\partial^\alpha\mathbf{G}d \xi
dxds
=\int_0^t\int\int\frac{\mb{P}_1(\partial^\alpha \mathbf{M})}{\mb{M}_\star }\mb{L}_\mathbf{M}\partial^\alpha\mathbf{G}d\xi dxds\\
\displaystyle\le
\frac{\widetilde{\sigma}}{16}
\int_0^t\int\int\frac{\nu(| \xi|)}{\mb{M}_\star }|\partial^\alpha
\mathbf{G}|^2d \xi dxds+C(\chi+\delta)\sum_{|\alpha^\prime|=1}\int_0^t\|\partial^{\alpha^\prime}(\phi,\psi,\zeta)\|^2ds+C\delta^{\frac12}
\end{array}
\end{equation}
and
\begin{equation}\label{hq2}
\begin{array}{l}
\di W_{17}=\int_0^t\int\int\f{v\mb{M}\pa^{\alpha}n_2}{\mb{M}_\star }\mb{N}_{\mb{M}}\pa^{\alpha}(\mb{P}_c F_2)d\xi dxds
+\int_0^t\int\int\f{\pa^{\alpha}(n_2 v\mb{M})-v\mb{M}\pa^{\alpha}n_2}{\mb{M}_\star }\mb{N}_{\mb{M}}\pa^{\alpha}(\mb{P}_c F_2)d\xi dxds\\
\di\leq |\int_0^t\int\int v\pa^{\alpha}n_2\big(\f{\mb{M}}{\mb{M}_\star }-1\big)\mb{N}_{\mb{M}}\pa^{\alpha}(\mb{P}_c F_2) d\xi dxds|
+C\int_0^t\int\big(|n_2||\pa^{\alpha}(v,u,\t)|\\
\di\quad+|n_2|\sum_{|\alpha'|=1}|\pa^{\alpha'}(v,u,\t)|^2
+\sum_{|\alpha'|=1}|\pa^{\alpha'}n_2||\pa^{\alpha'}(v,u,\t)|\big)
\big(\int\f{\nu(|\xi|)|\pa^{\alpha}(\mb{P}_cF_2)|^2}{\mb{M}_\star }d\xi\big)^{\f12}dxds\\
\di\leq C\eta_0\int_0^t\|\pa^{\alpha}n_2\|^2 ds+C(\chi+\d)\int_0^t\Big[\|\pa^{\alpha}(\phi,\psi,\z)\|^2
+\sum_{|\alpha'|=1}\|\pa^{\alpha'}(\phi,\psi,\z,n_2)\|^2\Big] ds\\
\di\quad+C(\chi+\d+\eta_0)\int_0^t\int\int\f{\nu(|\xi|)|\pa^{\alpha}(\mb{P}_cF_2)|^2}{\mb{M}_\star } d\xi dxds+C\d^{\frac12},
\end{array}
\end{equation}
where we have used the fact that
$$
\di \int_0^t\int v\pa^{\alpha}n_2\big(\int \mb{N}_{\mb{M}}(\pa^{\alpha}\mb{P}_c F_2)d\xi\big)dxds=0,
$$
and the small constant $\eta_0$ is defined in Lemma \ref{Lemma 4.2}. Substituting the above estimates into \eqref{h3} and choosing $\chi,\delta,\eta_0$ suitably small yield that
\begin{equation}
\begin{array}{l}
 \displaystyle \sum_{|\alpha|=2}\Big[\|\partial^\alpha\big(\f{\Phi_x}{v}\big)(\cdot,t)\|^2
 +\int\int\frac{|\partial^\alpha (F_1,F_2)|^2}{\mb{M}_\star }(x,\xi,t)d \xi
dx\Big]+\sum_{|\alpha|=2}\int_0^t\int\int\frac{\nu(|\xi|)}{\mb{M}_\star }|\partial^\alpha
(\mathbf{G},\mb{P}_c F_2)|^2d \xi dxds\\
\displaystyle\leq C \sum_{|\alpha|=2}\Big[\|\partial^\alpha(\Phi_{0x},v_0)\|^2
+\int\int\frac{|\partial^\alpha (F_{10},F_{20})|^2}{\mb{M}_\star }d \xi dx\Big]
+C\sum_{|\alpha'|=1}\|\pa^{\alpha'}(\Phi_{0x},v_0)\|^2+C\d^{\frac12}\\
\di +C(\chi+\delta+\eta_0)\sum_{|\alpha|=2}\int_0^t\|\partial^\alpha(\phi,\psi,\zeta,n_2,\Phi_x)\|^2ds
+C(\chi+\d)\sum_{|\alpha^\prime|=1}\int_0^t\|\partial^{\alpha^\prime}(\phi,\psi,\zeta,n_2,\Phi_x)\|^2ds\\
\di+C(\chi+\delta)\int_0^t\Big[\sum_{|\alpha^\prime|=1}\int\int\frac{\nu(| \xi|)}{\mb{M}_\star }|\partial^{\alpha^\prime}
(\mathbf{G}, \mb{P}_c F_2)|^2d \xi dx
+\sum_{0\leq |\beta'|\leq 1}\int\int\frac{\nu(| \xi|)}{\mb{M}_\star }|\partial^{\beta'}
(\widetilde{\mathbf{G}}, \mb{P}_c F_2)|^2d \xi dx\Big]ds\\
\displaystyle +C(\chi+\delta)\sum_{|\alpha^\prime|=1,|\beta^\prime|=1}\int_0^t\int\int
\f{\nu(|\xi|)|\pa^{\alpha'}\pa^{\beta'}(\mb{G}, \mb{P}_c F_2)|^2}{\mb{M}_\star }d\xi dxds
+C(\chi+\d)\int_0^t\|(\Phi_x,n_2)\|^2ds.
\end{array}
\end{equation}
From \eqref{phixtt-ap}, for $|\alpha|=2$, it holds that
\begin{equation}\label{phi-x-a}
\begin{array}{ll}
\di \int_0^t\|\partial^\alpha \Phi_x\|^2ds\leq
C \int_0^t\|n_{2x}\|^2 ds+C(\chi+\d)\int_0^t\Big[\|\pa^{\alpha}\phi\|^2+\|(\Phi_x,\Phi_{xt},n_2)\|^2\Big]ds\\
\di+C\sum_{|\alpha^\prime|=1}\int_0^t\int\int\f{\nu(|\xi|)|\partial^{\alpha^\prime}(\mb{P}_cF_2)|^2}{\mb{M}_\star }
d\xi dxds
 +C(\chi+\delta)\int_0^t\int\int\f{\nu(|\xi|)|\mb{P}_cF_2|^2}{\mb{M}_\star } d\xi dxds.
\end{array}
\end{equation}
The proof of Proposition \ref{Prop3.2} is completed.

%
%
\section{The Proof of Theorem \ref{thm}}
\setcounter{equation}{0}
Finally, by the Proposition \ref{Prop3.1} and \ref{Prop3.2}, it holds that
\begin{equation}\label{Final-E}
\begin{array}{ll}
\di\sup_{0\leq t<\i} \Big[\|(\phi,\psi,\zeta)(\cdot,t)\|^2_{H^1}+\|(\Phi_x,n_2,n_{2x})(\cdot,t)\|^2+\sum_{0\leq|\beta|\leq2}\int\int\f{|\partial^\beta(\widetilde{\mb{G}},\mb{P}_c F_2)|^2}{\mb{M}_\star }(x,\xi,t) d\xi dx\\
\di +\sum_{|\alpha^\prime|=1,0\leq|\beta^\prime|\leq1}\int\int \f{|\partial^{\alpha^\prime}\partial^{\beta^\prime}(\mb{G},\mb{P}_c F_2)|^2}{\mb{M}_\star }(x,\xi,t)d\xi dx+\sum_{|\alpha|=2}\int\int \f{|\partial^\alpha(F_1,F_2)|^2}{\mb{M}_\star }(x,\xi,t) d\xi dx\Big] \\
\di +\d\int_0^\i\||\hat w|(\phi,\psi,\z)\|^2dxds
+\sum_{1\leq|\alpha|\leq2}\int_0^\i\|\partial^{\alpha}(\phi,\psi,\zeta,n_2)\|^2ds+\int_0^\i\|(\Phi_x,\Phi_{xt},n_2)\|^2ds\\
\di +\sum_{1\leq|\alpha|\leq2}\int_0^\i\int\int \f{\nu(|\xi|)|\partial^\alpha(\mb{G},\mb{P}_c F_2)|^2}{\mb{M}_\star }d\xi dxds+\sum_{0\leq|\beta|\leq2}\int_0^\i\int\int\f{\nu(|\xi|)|\partial^\beta(\widetilde{\mb{G}},\mb{P}_c F_2)|^2}{\mb{M}_\star }d\xi dxds\\
\di +\sum_{|\alpha^\prime|=1,|\beta^\prime|=1}\int_0^\i\int\int \f{\nu(|\xi|)|\partial^{\alpha^\prime}\partial^{\beta^\prime}(\mb{G},\mb{P}_c F_2)|^2}{\mb{M}_\star }d\xi dxds\leq C(\mathcal{N}(0)^2+\delta^{\frac12})\leq C(\mathcal{E}(0)^2+\d^{\frac12}).
\end{array}
\end{equation}

Therefore, we close the a priori assumption by choosing suitably small positive constants $\varepsilon_0, \d$  and one has
\begin{equation*}
\begin{array}{ll}\label{aa}
\qquad \di \int_0^{+\infty}\int\int\frac{|(F_1-\mb{M}_{[\bar v,\bar u,\bar\t]}, F_2)_x|^2}{\mb{M}_\star } d\xi dxds\\[2mm]
\qquad\di\leq \int_0^{+\infty}\int\int\frac{|\big(\mb{M}_x-(\mb{M}_{[\bar v,\bar u,\bar\t]})_x,(n_2 v \mb{M})_x\big)|^2}{\mb{M}_\star } d\xi dxds
+\int_0^{+\infty}\int\int\frac{|(\mb{G}_x,(\mb{P}_cF_2)_x)|^2}{\mb{M}_\star } d\xi dxds\\[2mm]
\qquad\di\leq C\int_0^{+\infty}\|(\phi,\psi,\zeta,n_2)_x\|^2 ds
+C\d\int_0^{+\infty}\||\hat w|(\phi,\psi,\zeta)\|^2 ds+C\d^2\\
\di\qquad+\int_0^{+\infty}\int\int\frac{|(\mb{G}_x,(\mb{P}_cF_2)_x)|^2}{\mb{M}_\star } d\xi dxds
\leq C(\mathcal{N}(0)^2+\d^{\frac12})\leq C(\mathcal{E}(0)^2+\d^{\frac12}).
\end{array}
\end{equation*}
From the Vlasov-Poisson-Boltzmann system \eqref{vpb-l}, we  can obtain
\begin{equation*}
\begin{array}{ll}\label{bb}
&\di \int_0^{+\infty}\Big|\frac{d}{dt}\int\int\frac{|(F_1-\mb{M}_{[\bar v,\bar u,\bar\t]}, F_2)_x|^2}{\mb{M}_\star }d\xi dx\Big|ds
\leq C(\mathcal{E}(0)^2+\d^{\frac12}).
\end{array}
\end{equation*}
Therefore, one has
\begin{equation*}
\begin{array}{ll}
\di\int_0^{+\infty}\Big(\int\int\frac{|(F_1-\mb{M}_{[\bar v,\bar u,\bar\t]}, F_2)_x|^2}{\mb{M}_\star } d\xi dx
+\Big|\frac{d}{dt}\int\int\frac{|(F_1-\mb{M}_{[\bar v,\bar u,\bar\t]}, F_2)_x|^2}{\mb{M}_\star } d\xi dx\Big|
\Big)ds<\infty,
\end{array}
\end{equation*}
which implies that
$$
\lim_{t\rightarrow+\infty} \int\int\frac{|(F_1-\mb{M}_{[\bar v,\bar u,\bar\t]}, F_2)_x|^2}{\mb{M}_\star } d\xi dx=0.
$$
By Sobolev inequality
$$
\begin{array}{ll}
\di \|\int\frac{|(F_1-\mb{M}_{[\bar v,\bar u,\bar\t]}, F_2)|^2}{\mb{M}_\star } d\xi\|^2_{L^{\infty}_x}\\
\di \leq C\Big(\int\int\frac{|(F_1-\mb{M}_{[\bar v,\bar u,\bar\t]}, F_2)|^2}{\mb{M}_\star } d\xi dx\Big)
\cdot\Big(\int\int\frac{|(F_1-\mb{M}_{[\bar v,\bar u,\bar\t]}, F_2)_x|^2}{\mb{M}_\star } d\xi dx\Big),
\end{array}
$$
we can prove that
\begin{equation}\label{pr1}
\lim_{t\rightarrow+\infty} \sup_x\int\frac{|(F_1-\mb{M}_{[\bar v,\bar u,\bar\t]}, F_2)|^2}{\mb{M}_\star } d\xi =0.
\end{equation}
Similarly, one can prove that
$$
\lim_{t\rightarrow+\infty} \|(\Phi_x,n_2)\|=0.
$$
Thus the time-asymptotic convergence of the solutions $(F_1,F_2)$ to the viscous contact wave $\mb{M}_{[\bar v, \bar u, \bar \t]}$
can be derived directly from \eqref{pr1}, it holds that
 $$
 \mathcal{E}(t)\leq C(\mathcal{N}(t)+\delta^{\frac12})\leq C(\mathcal{E}(0)+\d^{\frac12}),\qquad \forall t\in [0,+\infty],
 $$
which proves \eqref{Et}. And the proof of Theorem \ref{thm} is completed.

%
\end{document}